\UseRawInputEncoding

\documentclass[final,3p,times]{elsarticle}

\usepackage{placeins}
\usepackage[normalem]{ulem}

\usepackage{makeidx}
\usepackage{amsfonts}
\usepackage{verbatim}
\usepackage[ansinew]{inputenc}
\usepackage[dvipsnames]{pstricks}
\usepackage{epstopdf}
\usepackage{bbm}
\usepackage{pdfsync}
\usepackage[mathscr]{euscript}
\usepackage{appendix}

\usepackage{amssymb}
\usepackage{amsmath}
\usepackage{tabularx}
\usepackage{graphicx}
\usepackage{subfigure}
\usepackage{xcolor}
\usepackage{csquotes}
\usepackage{algorithm}
\usepackage{algpseudocode}

\newcommand{\ba}{{\mathbf{a}}}

\newcommand{\bb}{{\mathbf{b}}}

\newcommand{\bd}{{\mathbf{d}}}

\newcommand{\bD}{{\mathbf{D}}}

\newcommand{\be}{{\mathbf{e}}}

\newcommand{\bF}{{\mathbf{F}}}

\newcommand{\bH}{{\mathbf{H}}}

\newcommand{\bI}{{\mathbf{I}}}
\newcommand{\bK}{{\mathbf{K}}}

\newcommand{\bL}{{\mathbf{L}}}

\newcommand{\bM}{{\mathbf{M}}}
\newcommand{\bn}{{\mathbf{n}}}
\newcommand{\bN}{{\mathbf{N}}}

\newcommand{\bp}{{\mathbf{p}}}

\newcommand{\prm}{{\prime}}

\newcommand{\bt}{{\mathbf{t}}}

\newcommand{\bT}{{\mathbf{T}}}

\newcommand{\bu}{{\mathbf{u}}}

\newcommand{\bX}{{\mathbf{X}}}

\newcommand{\bx}{{\mathbf{x}}}

\newcommand{\bzero}{{{\bf{0}}}}

\newcommand{\bchi}{{\boldsymbol{\chi}}}

\newcommand{\bmu}{{\boldsymbol{\mu}}}

\newcommand{\beps}{{\boldsymbol{\varepsilon}}}

\newcommand{\boldeta}{{\boldsymbol{\eta}}}

\newcommand{\br}{{\boldsymbol{r}}}
\newcommand{\bJ}{{\boldsymbol{J}}}

\newcommand{\std}[1]{{#1}_\mathrm{macro}}

\newcommand{\der}[2]{{\frac{\partial #1}{\partial #2}}}

\journal{Computer Methods in Applied Mechanics and Engineering}

\begin{document}

\begin{frontmatter}





\title{Strain localization in softening plasticity without modifying standard constitutive models: a deformable Cosserat approach}


\author[1]{Andrea Panteghini}

\address[1]{Department of Civil, Environmental, Architectural Engineering and Mathematics (DICATAM), University of Brescia, Italy, Email: andrea.panteghini@unibs.it}

\author[2]{MB Rubin}

\address[2]{Faculty of Mechanical Engineering, Technion-Israel Institute of Technology, 32000 Haifa, Israel, Email: mbrubin@tx.technion.ac.il}

\begin{abstract}
This paper presents a formulation for strain localization in softening plasticity based on a deformable Cosserat model. The approach enables the direct use of standard elastoplastic constitutive models formulated for a classical Cauchy continuum, without any modification of the stress update algorithm or consistent tangent operator.
The key feature of the proposed framework is a strict separation between dissipative and energetic mechanisms: all dissipation is confined to the macro-continuum, while the micro-continuum contributes exclusively through linear elastic terms associated with the director field. As a result, the constitutive structure of the underlying elastoplastic model is preserved, and existing models can be employed as true black-box components.
The internal length scale arises naturally from the micro-continuum and governs the development, interaction, and selection of localization patterns, rather than acting as a diffusive or artificial parameter.
The formulation is straightforward to implement within standard finite element frameworks, requiring only additional linear contributions to the residual and tangent operators, without altering the structure of the constitutive integration.
The performance of the proposed approach is assessed through benchmark problems involving shallow foundations on soil, which represent a particularly demanding test due to the onset of complex, interacting, and unstable localization mechanisms. Both Tresca and Matsuoka-Nakai plasticity models are considered, including cases exhibiting highly unstable post-peak responses.
Numerical results show that load-displacement responses, dissipated energy, and shear-band patterns converge upon mesh refinement, even in the presence of strongly nonlinear and interacting localization processes. These findings demonstrate that the proposed framework provides a robust and physically consistent approach for the analysis of strain localization in softening plasticity.
\end{abstract}





\end{frontmatter}




\section{Introduction}
\label{intro}
Elasto-plastic constitutive models within the classical Cauchy continuum are widely used for the analysis of structural and geotechnical problems. 
However, capturing key features of material behavior -- particularly for geomaterials -- often requires the introduction of strain softening, which leads to severe numerical and theoretical difficulties.

It is well established that, in the presence of softening, boundary value problems exhibit non-uniqueness of the solution, which manifests in numerical simulations as pathological mesh sensitivity. In particular, computed responses do not converge upon mesh refinement, and the width of localized deformation bands becomes dependent on the discretization size (see, e.g., \cite{de1991simulation,DeBorst1991,sabet2019mesh}). 
From a practical standpoint, these issues severely affect the robustness of numerical analyses. In many cases, the introduction of softening leads to loss of convergence, with simulations failing to complete even during the initial loading steps. 
These difficulties are now understood as a direct consequence of the loss of well-posedness of the governing equations.  In quasi-static problems, the onset of localization is associated with a loss of ellipticity, which underlies both the observed mesh dependence and the severe numerical instabilities \cite{Sabet2019}.

Several approaches have been proposed in the literature to overcome the pathological mesh dependence associated with strain localization in softening materials. 

One class of approaches is based on mechanical enrichments of the Cauchy continuum. For example, the Cosserat continuum, originally introduced by the Cosserat brothers \cite{cosserat1909theorie}, augments the classical kinematics by introducing independent rotational degrees of freedom governed by additional balance equations, thereby incorporating intrinsic length scales into the formulation.
In this work, this formulation is referred to as a \emph{rigid Cosserat model}, since the enriched kinematics can be characterized by the rotation of a rigid triad of director vectors.
In particular, rigid Cosserat models have long been regarded as a natural framework for the analysis of shear banding, starting from the seminal work of M\"uhlhaus and Vardoulakis \cite{muhlhaus1987thickness}, which showed that the introduction of an intrinsic length scale leads to localization bands of finite thickness. Starting from the pivotal computational studies of de Borst \cite{de1991simulation} and Peri\'c \cite{peric1994error}, the rigid Cosserat model has been widely employed in numerical analyses and benchmark problems \cite{sharbati2006computational, khoei2008enriched, khoei20103d, sabet2019mesh, PL2022a, PL2022b}, as well as in practical applications involving localization phenomena \cite{russo2020thermomechanics, ebrahimian2012cosserat}. However, it is now well recognized that this regularization is only partial: while objective results can be obtained in problems characterized by a single dominant shear band, mesh dependence may reappear when multiple bands nucleate, interact, or evolve sequentially \cite{duretz2023comparison, PR2026}.

This class further includes micromorphic models \cite{eringen1964nonlinear, suhubl1964nonlinear, forest2009micromorphic, forest2020micromorphic}, that are closely related but conceptually distinct from the rigid Cosserat. They introduce a fully independent micro-deformation tensor at each material point, thereby allowing for a very general and systematic description of microstructural effects and internal length scales.

A second class of approaches aims instead at restoring well-posedness through explicit regularization of the governing equations. These include gradient-enhanced plasticity and damage models \citep{peerlings1996gradient, peerlings1998gradient, jirasek2003comparison, lorentz2011gradient, anand2026fracture}, viscoplastic formulations \citep{needleman1992analyses}, and nonlocal approaches based on spatial averaging \citep{bazant1988nonlocal, bazant2002nonlocal}. More recently, phase-field formulations \citep{miehe2010phasefield_cma, miehe2010phasefield_nme, duda2015phase, miehe2017phase} have gained increasing attention, introducing auxiliary field variables governed by additional differential equations to regularize localization phenomena.
While these approaches have proven effective in controlling mesh dependence, they are typically introduced as mathematical regularization techniques, rather than as mechanically enriched continuum models.

A common limitation of both classes of approaches is that they require the formulation of  constitutive models, in which the additional fields are constitutively coupled to the primary mechanical response. As a consequence, the resulting formulations cannot be readily combined with existing elastoplastic models based on the Cauchy continuum.

The framework adopted in this work is based on the deformable Cosserat model recently proposed in \citep{rubin2025thermomechanical, rubin2026correction}, which enriches the continuum with a triad of deformable director vectors. 
It was previously shown in \cite{PR2026} that, in the context of elasticity with damage, the original formulation is able to predict mesh- and rate-independent shear band evolution, since elastic and total deformations essentially coincide. 
However, preliminary investigations in the context of elastoplasticity revealed that this property is lost when elastic and plastic deformations are no longer equivalent, and that the original formulation does not provide a consistent description of localization in this case. 
For this reason, a modified version of the deformable Cosserat model has been developed for elastoplastic response involving volumetric plastic deformations, and is presented in the companion theoretical paper \citep{miles}.
It is also noted that, as discussed in \citep{rubin2026AM}, in the small deformation setting adopted here, the resulting formulation exhibits strong similarities with classical microstructural models developed by Mindlin \citep{mindlin1964microstructure}. 

In the deformable Cosserat model, the material is described as two interacting continua: a macro-continuum, which follows a standard Cauchy continuum description, and a micro-continuum, associated with the evolution of a deformable director triad. The two continua are constitutively uncoupled: the macro-continuum governs the entire dissipative response through standard elastoplasticity, while the micro-continuum  introduces additional energetic terms associated with the \emph{relative} deformation of the directors with respect to the macro-continuum deformation. These energetic terms become significant in regions of strong strain localization and govern the structure of the emerging deformation patterns.

This structure leads to a key computational advantage. The constitutive response of the macro-continuum can be obtained directly from any standard elastoplastic model developed for classical Cauchy continua, without modifying either the stress update algorithm or the consistent tangent operator. In other words, the constitutive law can be employed as a true black-box component. Once this contribution is computed, the additional terms associated with the micro-continuum are linear and can be assembled straightforwardly within the finite element formulation.

The objective of the present work is to assess the capability of the deformable Cosserat model to predict strain localization in softening plasticity within a physically meaningful setting. In particular, the internal length $\ell$ is treated as a genuine material parameter that governs the selection of the deformation mechanisms, rather than as a numerical regularization parameter. 

For a fixed value of $\ell$, the formulation yields mesh-independent solutions that converge towards well-defined localization patterns. At the same time, varying $\ell$ can lead to qualitatively different responses, both in terms of global behavior and the associated shear band structures. The proposed framework therefore does not merely regularize the solution, but provides a physically grounded mechanism for the selection and evolution of localized deformation patterns, while enabling the direct use of standard constitutive models as true black-box components.


\section{The deformable Cosserat model}

The formulation adopted in this work is the small-deformation version of the deformable Cosserat model developed in \citep{miles}. The general framework, which is thermodynamically consistent, is derived in an Eulerian large-strain setting. Here, it is specialized to small strains and rotations.

The model consists of two interacting continua: a \emph{macro}-continuum and a \emph{micro}-continuum. 
The macro-continuum governs the entire dissipative response and is described by a standard elastoplastic constitutive model. The micro-continuum is non-dissipative and contributes only through additional energetic terms associated with the mismatch between the directors $\bd_i$
and material line elements $\ba_i$, as well as through the corresponding director curvatures.

 The two continua are constitutively uncoupled: the macro-response remains entirely standard, whereas the micro-response is governed by its own elastic energy. Their interaction is therefore not introduced through constitutive coupling terms, but arises from the kinematics together with the balance equations and boundary conditions.

This separation between macro- and micro-responses is a key feature of the formulation. All dissipation is confined to the macro-continuum, while the micro-continuum acts as an energetic mechanism that becomes relevant only when the deformation field develops sufficiently strong spatial variations. In this way, the model can influence localization patterns without modifying the standard elastoplastic constitutive structure of the macro-continuum.

\subsection{Kinematics}

\begin{figure}[t]
    \centering
    \includegraphics[width=0.62\textwidth]{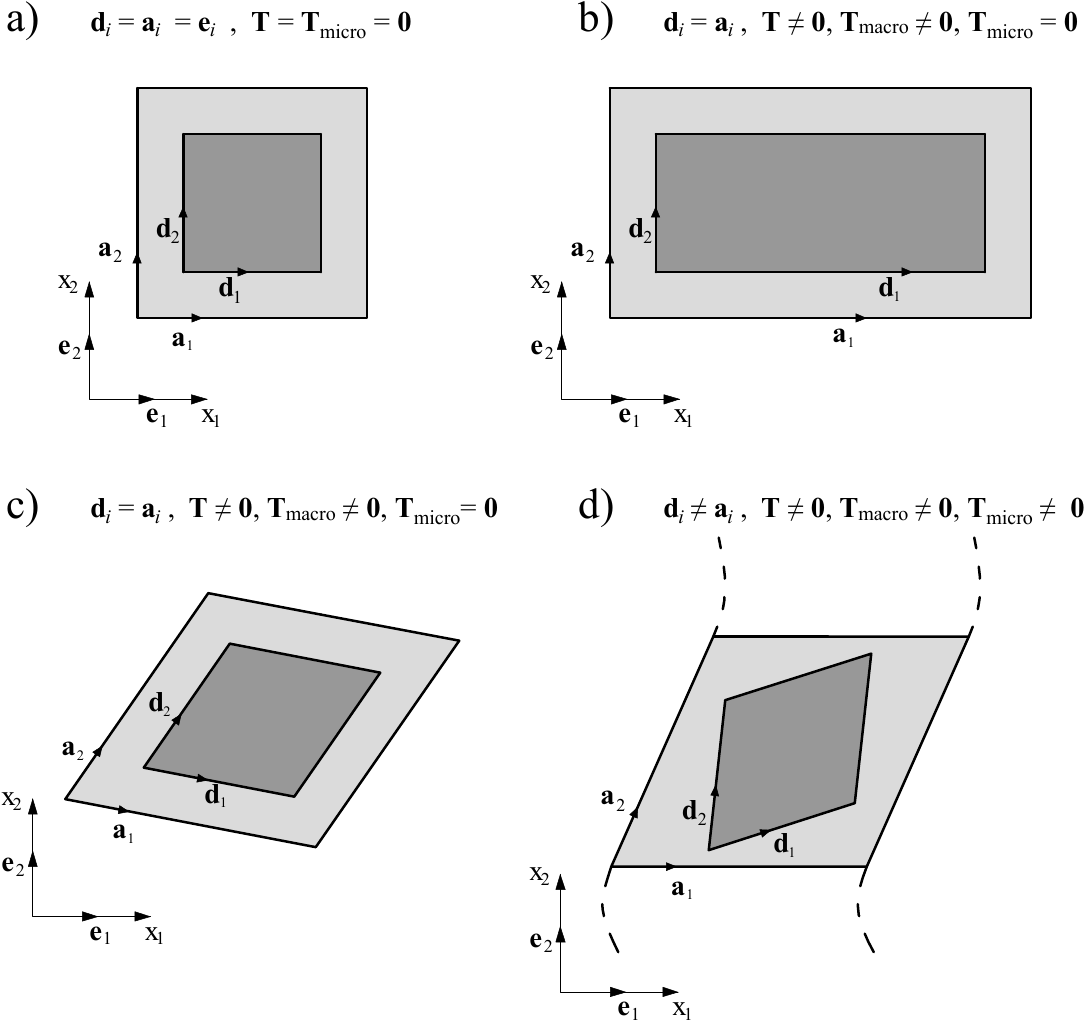}
   \caption{Schematic representation of the kinematics of the deformable Cosserat model in a 2D setting. The vectors $\ba_i$ denote the material line elements associated with the displacement gradient, whereas $\bd_i$ denote the directors describing the micro-continuum. (a) Undeformed configuration, for which $\bd_i=\ba_i=\be_i$ and therefore $\bT=\bT_\mathrm{micro}=\bzero$. (b) Homogeneous stretch: although $\bT\neq \bzero$ and $\bT_\mathrm{macro}\neq \bzero$, the directors still follow the material line elements, so that $\bd_i=\ba_i$ and $\bT_\mathrm{micro}=\bzero$. (c) Homogeneous shear/distortion: again $\bd_i=\ba_i$, hence $\bT_\mathrm{micro}=\bzero$ and the response is entirely governed by the macro-continuum. (d) Non-homogeneous deformation: the directors no longer follow the material line elements, i.e. $\bd_i\neq\ba_i$, so that a mismatch develops and $\bT_\mathrm{micro}\neq \bzero$. The dashed lines schematically represent non-vanishing director curvatures, which generate corresponding micro-couples.}
    
    \label{fig:kinematics}
\end{figure}

With reference to Fig.~\ref{fig:kinematics}, the kinematics is described in terms of a triad of vectors $\ba_i$, representing material line elements, and a triad of directors $\bd_i$, describing the micro-continuum. Under the assumption of small strains and rotations, the vectors $\ba_i$ are obtained from the displacement gradient as
\begin{equation}
    \ba_i=\be_i+\left(\frac{\partial \bu}{\partial\bx}\right)\be_i\,,
\end{equation}
where $\bu$ is the displacement field and $\be_i$ is a fixed orthonormal triad.

The directors $\bd_i$ are determined from a non-symmetric second-order tensor $\boldeta$, 
as
\begin{equation}
    \bd_i=\be_i+\boldeta^T\be_i\,.
\end{equation}
The tensor $\boldeta$ is treated as an independent kinematic variable, whose evolution is governed by the balance equations.
This allows the deformation of the director triad to differ from the macroscopic deformation, introducing additional kinematic freedom with respect to classical continuum models. Accordingly, each material point is characterized by $12$ degrees of freedom in a general 3D setting (three displacement components and nine components of $\boldeta$), which reduce to $6$ degrees of freedom in a 2D plane strain setting (two displacement components and four components of $\boldeta$).

\subsection{Balance laws and boundary conditions}

Under the assumption of no external micro-couples per unit of mass, and neglecting all the inertia terms, the balances of linear momentum and  director momentum respectively read
\begin{equation}
    \mathrm{div}\,\bT+ \rho_0\bb=\bzero\,,
\end{equation}
\begin{equation}
    \bT_\mathrm{micro}+\mathrm{div} \,\bM^i \otimes\be_i=\bzero\,,
    \label{eq:dirmomentum}
\end{equation}
where $\bT$ is a second-order (non-symmetric) stress tensor, $\bb$ are the body forces (per unit mass), $\rho_0$ is the zero-stress density, $\bM^i$
are second-order director couple tensors\footnote{Following \cite{miles}, $\bM^i=M^{i}_{jk}\be_j\otimes \be_k$} having the dimensions of $[F][L]^{-1}$ and $\bT_\mathrm{micro}$ denotes the non-symmetric stress associated with the micro-continuum. Here, $\mathrm{div(\,\cdot\,)}$ is the divergence operator with respect to position $\bx$, $\ba\otimes\bb$ denotes the tensor product between two vectors $\ba,\,\bb$, and the usual summation convention is used for repeated indices.
The director momentum equation \eqref{eq:dirmomentum}  governs the evolution of the director field, relating the micro-stress to the divergence of the micro-couples.
Since the micro-stress and the micro-couples are assumed to be linearly elastic functions of $\bchi$ and the director gradients, Eq.~\eqref{eq:dirmomentum} is linear in the enhanced kinematic tensor $\boldeta$ for a given displacement field. This property will play a key role in the simplicity of the resulting finite element formulation.

The traction vector $\bt$ and the contact director couples $\bmu^i$, applied to a material surface $\partial P$ with unit outward normal $\bn$ are given by
\begin{equation}
    \bt=\bT\cdot\bn\,,\quad \bmu^i=\bM^i\cdot \bn\,,
\end{equation}
so that their external power on $\partial P$ is
\begin{equation}
    \int_{\partial P}\left(\bt \cdot \dot \bu +\bmu^i\cdot\dot{\bd}_i\right)\,\mathrm{d}S\,,
\end{equation}
where $\mathrm{d}S$ is the element of area on $\partial P$.

\subsection{Specific strain energy density and constitutive assumptions}
It is assumed that the specific (per unit mass) strain energy density $\Sigma$ can be additively decomposed as
\begin{equation}   \Sigma=\Sigma_\mathrm{macro}+\Sigma_\mathrm{micro}\,.
\label{eq:decomp}
\end{equation}
Accordingly, the stress tensor is written as
\begin{equation}
    \bT=\bT_\mathrm{macro}+\bT_\mathrm{micro}\,,
\label{eq:decompstress}
\end{equation}
here $\bT_\mathrm{macro}$ is a \emph{symmetric} second-order tensor describing the stress in the macro-continuum, consistently with the assumption that the latter obeys a general standard isotropic elastoplastic constitutive law.

As shown in Fig.\ref{fig:kinematics}, the micro-continuum stress $\bT_\mathrm{micro}$ is related to the \emph{mismatch} between the directors $\bd_i$ and material line elements described by $\ba_i$, which is measured by the non-symmetric second-order tensor\footnote{Using the definition of $\chi^j_i$ in \cite{miles}, $\bchi=\chi^j_i\,\be_j\otimes\be_i$}
\begin{equation}
\bchi=\partial \bu/\partial\bx-\boldeta^T.
\label{eq:defchi}
\end{equation}
Thus, the micro-stress is driven by the tensor $\bchi$, which accounts for both relative deformation and relative rotation.

Since the micro-stress and micro-couples are assumed to be linearly elastic, $\Sigma_\mathrm{micro}$ is chosen as a quadratic form. The micro-curvature measure $\zeta^k_{ij}$ (which corresponds to the symbol $\chi^k_{ij}$ in \cite{miles}) is introduced as
\begin{equation}
    \zeta^k_{ij}=  \left(\eta_{ij}+\eta_{ji}\right)_{,k}\,,
\end{equation}
where $(\,\cdot\,)_{,k}$ denotes partial differentiation with respect to the $k$-th spatial coordinate $x_k$.  The strain energy density of the micro-continuum is taken as:
\begin{equation}
    \rho_0\Sigma_\mathrm{micro}=\frac{G}{2} \left\{ 
    k_1 \left(\bchi : \bI \right)^2+k_2\, \bchi^\prm:\bchi^\prm+\,\ell^2 \zeta^k_{ij} \zeta^k_{ij}
    \right\}\,,
    \label{eq:freemicro}
\end{equation}
where $G$ is the shear modulus, $k_1$ and $k_2$ are non-negative dimensionless material constants, $\ell$ is a material length, and $\bI$ is the second-order identity tensor. Moreover, $(\cdot)^\prm$ denotes the deviatoric part of a second-order tensor and $:$ denotes the inner product between two second-order tensors.

Consistently with the small-strain framework here adopted, by assuming a constant density $\rho_0$,  following \cite{miles}, the rate of dissipation can be written as:
\begin{equation}
   \mathcal{D}=\bT : \bL-\bT_\mathrm{micro}: \dot \boldeta^T+M^i_{jk} \dot{\eta}_{ij,k} -\rho_0\dot \Sigma \geq0\,,   \label{eq:dissgen}
\end{equation}
where $\bL=\partial\dot\bu/\partial\bx$ is the velocity gradient. 
By considering the stress decomposition \eqref{eq:decompstress}, recalling that $\bT_\mathrm{macro}$ is symmetric, and using the definition of $\bchi$ given in Eq. \eqref{eq:defchi}, the rate of dissipation can be rewritten as
\begin{equation}
\mathcal{D}=\bT_\mathrm{macro}:\bD+\bT_\mathrm{micro}:\dot{\bchi}
+M^i_{jk} \dot{\eta}_{ij,k}-\rho_0 \dot \Sigma \geq0.
\label{eq:dissgen1}
\end{equation}
where  $\bD$ is the symmetric part of the velocity gradient.
It is assumed that the constitutive response of the micro-continuum is elastic, so that the mechanical power of the micro-continuum is balanced by the derivative of the micro-strain energy:
\begin{equation}
\bT_\mathrm{micro}:\dot{\bchi} +M^i_{jk} \dot{\eta}_{ij,k}-\rho_0 \dot \Sigma_\mathrm{micro} = 0 \,.
\label{eq:dissgen_micro}
\end{equation}
Since this identity must hold for arbitrary admissible rates, the micro-stress and micro-couples are given by
\begin{equation}
\begin{gathered}    \bT_\mathrm{micro}=\rho_0\frac{\partial \Sigma_\mathrm{micro}}{\partial \bchi}=G\left(k_1 \bI\otimes \bI+k_2\, \mathcal{I}_d \right)\bchi\,,\\
    M^i_{jk}=\rho_0\frac{\partial \Sigma_\mathrm{micro}}{\partial \eta_{ij,k}}=2G \ell^2 \zeta^k_{ij}\,,
\end{gathered}
\label{eq:tmicro}
\end{equation}
where $\mathcal{I}_d= \mathcal{I}-(1/3)\bI\otimes\bI$ is the deviatoric projector and $\mathcal{I}$ is the fourth-order identity tensor. Then, with the help of \eqref{eq:dissgen_micro}, the rate of dissipation \eqref{eq:dissgen1} requires
\begin{equation}
    \mathcal{D}=\bT_\mathrm{macro}:\bD-\rho_0\dot\Sigma_\mathrm{macro}\geq0\,.
    \label{eq:dissmacro}
\end{equation}
The macro strain energy $\Sigma_\mathrm{macro}$ is assumed to be an isotropic function of the \emph{symmetric} elastic strain $\beps_e$, which is governed by the following evolution equation:
\begin{equation}
    \dot \beps_e = \bD-\bD_p \,,
    \label{eq:evoleps}
\end{equation}
where $\bD_p$ is the symmetric plastic strain rate of the macro-continuum. Moreover, the constitutive equation for the macro-stress is specified by
\begin{equation}
    \bT_\mathrm{macro} = \rho_0 \frac{\partial \Sigma_\mathrm{macro}}{\partial \beps_e}\,,
\end{equation}
so that the rate of dissipation \eqref{eq:dissmacro} reduces to the  requirement for standard plasticity in the macro-continuum:
\begin{equation}
    \mathcal{D}= \bT_\mathrm{macro}: \bD_p \geq0\,.
\end{equation}
It is emphasized that all dissipation is confined to the macro-continuum, while the micro-continuum provides a purely energetic contribution.

\section{Finite element formulation}
\label{sec5}

The internal virtual power $\delta W^\mathrm{int}$ of a continuum occupying a region $P$ is expressed as\footnote{
It should be noted that Eq. \eqref{eq:VirtPower} is 
a special case, for the simple \emph{constitutive} assumptions here adopted. For a more general expression, the reader may refer to \cite{miles}.}
\begin{equation}
\delta W^\mathrm{int}=\int_P \big[T_{ij}  \delta u_{i,j}-\left(\bT_\mathrm{micro} \right)_{ij} \delta \eta_{ji}+\frac{1}{2}M^i_{jk} \left(\delta\eta_{ij,k} +\delta\eta_{ji,k}\right)\big]\,
\mathrm{d}V\,,
\label{eq:VirtPower}
\end{equation}
where $\delta\bu$ and $\delta\boldeta$ denote arbitrary kinematically admissible variations, and $\mathrm{d}V$ is the volume element.

The adopted FE discretization is isoparametric, such that the same shape functions are employed to interpolate both the reference geometry and the displacement field, i.e.
\begin{equation}
\nonumber
x_{j} (\br)= \sum_{k=1}^nN^{(k)}(\br) \hat x_{j}^{(k)}, \quad
u_{j} (\br)= \sum_{k=1}^nN^{(k)}(\br) \hat u_j^{(k)}\,,
\end{equation}
where $x_{j}$ and $u_j$ are the $j$-th components of the coordinates and of the displacement, respectively, $n$ is the total number of nodes, 
$N^{(k)}(\br)$ is a standard polynomial shape function referred to the $k$-th node, $\br$ is the vector containing the intrinsic coordinates of the parent element, and $\hat x_{j}^{(k)}$ and $\hat u_{j}^{(k)}$ are, respectively, the $j$-th coordinate and the $j$-th displacement component of the $k$-th node\footnote{In this work the hat symbol, $\hat{(\,\cdot\,)}$ denotes the nodal values of the corresponding variables.}.
The tensor components $\eta_{ij}$ are discretized as:
\begin{equation}
\eta_{ij}(\br)= \sum_{k=1}^{n_\eta} N_\eta^{(k)}(\br) \hat \eta_{ij}^{(k)} \,,
\end{equation}
where $n_{\eta}\leq n$ is the number of nodes in which $\eta_{ij}$ is discretized, $N_\eta^{(k)}(\br)$ is  a standard polynomial shape function (that may be of different order with respect to $N^{(k)}(\br)$) referred to the $k$-th node, and $\hat \eta_{ij}^{(k)}$ is the $ij$ component of $\boldeta$ at the node $k$.

The matrices $\bH$ and $\bH_\eta$, respectively  containing the derivatives of the shape functions $N^{(k)}$ and $N_\eta^{(k)}$ with respect to the  coordinate system,  can be computed as
\begin{equation}
\bH 
= \bJ^{-1} \bH_{r},\quad \bH_\eta=\bJ^{-1} \bH_{\eta r}, \quad \bJ=\bH_{r} \hat{\bX} \,,
\end{equation}
where $\bH_r$ and $\bH_{\eta r}$ are, respectively, the matrices that contain the derivatives of the shape functions $N^{(k)}$ and $N^{(k)}_\eta$ with respect to $\br$ and $\bJ$ is the Jacobian matrix, and $\hat{\bX}$ is a matrix containing the nodal coordinates.

By assuming 2D plane strain, and adopting an array notation\footnote{
	At the implementation level, the nonzero components of each second-order 2D tensor (symmetric or non-symmetric) are stored in a 5-component array. For instance,
	\begin{equation}
\nonumber
	[\bT] = \begin{bmatrix}
	T_{11} & T_{12} & 0\\
	T_{21} & 
	T_{22} & 0\\
	0 & 0 & T_{33}
	\end{bmatrix}\rightarrow
	\bT=\begin{bmatrix}
	T_{11}\\
	T_{22}\\
	T_{33}\\
	T_{12}\\
	T_{21}
	\end{bmatrix},
    \quad
    \boldeta=\begin{bmatrix}
	\eta_{11}\\
	\eta_{22}\\
	0\\
	\eta_{12}\\
	\eta_{21}
    \end{bmatrix} \,.
	\end{equation}

Moreover,
\begin{equation}
\nonumber
	[\boldeta]_{ij,k} \rightarrow\boldeta_{,k}=
	\begin{bmatrix}
	\eta_{11,1}\\
	\eta_{22,1}\\
	\eta_{12,1}\\
	\eta_{21,1}\\
    \eta_{11,2}\\
	\eta_{22,2}\\
	\eta_{12,2}\\
	\eta_{21,2}\\
	\end{bmatrix}, \quad
    [\bM]^i_{jk}\rightarrow \bM=\begin{bmatrix}
    M^1_{11}\\
    M^2_{21}\\
    M^1_{21}\\
    M^2_{11}\\
    M^1_{12}\\
    M^2_{22}\\
    M^1_{22}\\
    M^2_{12}\\
    \end{bmatrix} \,,
	\end{equation}
} the fields $\delta u_{i,j}$, $\delta \eta_{ji}$ and $\delta \eta_{ij,k}$ can be discretized as  functions of the  nodal displacements and tensor variations $\delta \hat \bu, \delta \hat \boldeta$  as:
\begin{equation}
\delta u_{i,j}=\mathcal{A} \delta \hat \bu, \quad
\delta \eta_{ji}\equiv \delta \boldeta^T= \mathcal{N}_\tau \delta \hat \boldeta, \quad
\frac{1}{2}\left(\delta \eta_{ij,k}+
\delta \eta_{ji,k}\right)
= \mathcal{G}_\mathrm{sym} \delta \hat \boldeta \,,
\end{equation}
where $\mathcal{A}(\bH)$ is a linear function\footnote{All the FE operators,  are reported in \ref{Apn_FEOper}.} that maps the components of $\bH$  into a $5\times  2n$ matrix, $\mathcal{N}_\tau (\bN_{\eta})$ is a linear function that maps the shape functions $\bN_{\eta}$ into a $5\times  4n_\eta$ matrix, and  $\mathcal{G}_\mathrm{sym}(\bH_{\eta})$, here defined in its symmetrized form, consistently with the adopted evolution equation for $M^i_{jk}$, maps the components $\bH_{\eta}$ into a $8\times  4n_\eta$ matrix.

Under these assumptions, the internal virtual power \eqref{eq:VirtPower} can be discretized as
\begin{equation}
\delta W^{\mathrm{int}}=\int_P \big[ \delta
\hat \bu^T  \mathcal{A}^T \bT
+ \delta
\hat \boldeta^T \left(  \mathcal{G}_\mathrm{sym}^T  \bM - \mathcal{N}_\tau^T \bT_\mathrm{micro} \right) \big] 
\mbox{d}V\,,
\end{equation}
which, being valid for any kinematically admissible value of $\delta \hat \bu, \delta \hat \boldeta$, gives the following internal forces:
\begin{equation}
    \bF^{\mathrm{int}}_\mathrm{u}=\int_P   \mathcal{A}^T \bT\,\mathrm{d}V\,,\quad
    \bF^{\mathrm{int}}_{\eta}=\int_P \left(\mathcal{G}_\mathrm{sym}^T  \bM - \mathcal{N}_\tau^T \bT_\mathrm{micro}\right)
\mathrm{d}V \,.
\label{eq:residuals}
\end{equation}
The non-standard micro-stress $\bT_\mathrm{micro}$ and the micro-couples $M^i_{jk}$, here collected in the vector $\bM$, are computed directly from the nodal variables through the linear elastic relations of the micro-continuum.
In particular, from Eq. \eqref{eq:tmicro} it results:
\begin{equation}
\begin{gathered}
    \bT_\mathrm{micro}=G\left(k_1 \bI\otimes \bI+k_2\, \mathcal{I}_d \right) 
 \bchi\,=G\left(k_1 \bI\otimes \bI+k_2\, \mathcal{I}_d \right) \left(\mathcal{A}\,\hat \bu -\mathcal{N}_\tau\, \hat \boldeta\right),\\
    M^i_{jk}= 2 G \ell^2 \left(\eta_{ij,k}+\eta_{ji,k}\right)\rightarrow    
    \bM=4G\ell^2 \,\mathcal{G}_\mathrm{sym}\, \hat \boldeta\,.
\end{gathered}
\label{eq:nonstdstress}
\end{equation}
The macro-stress $\bT_\mathrm{macro}$ is the only part affected by elastoplasticity. It is a symmetric stress that depends on the symmetric elastic strain of the macro-continuum, on the increment of the symmetric part of the displacement gradient,
\begin{equation}
    \Delta \beps_\mathrm{macro} = \mathcal{A}_\mathrm{sym}\,\Delta \hat \bu\,,
    \label{eq:epsstd}
\end{equation}
and on a set of internal variables $\std{\bp}$ associated exclusively with the standard plastic response.  Again, $\mathcal{A}_\mathrm{sym}(\bH)$, reported in \ref{Apn_FEOper} for a 2D plane strain element, is a linear function that maps the components of $\bH$  into a $5\times  2n$ matrix.

Accordingly, the macro stress $ \bT_\mathrm{macro}$ is obtained from a standard elastoplastic constitutive update as
\footnote{
For implementation purposes, standard constitutive routines formulated in compact Voigt notation can be readily employed. To this end, we introduce the compact (symmetric) representations $\std{\underline{\bT}}$ and $\std{\underline{\Delta \beps}}$, collecting the independent components of second-order symmetric tensors.

The relation between the full tensorial representation and the compact one is established through linear mapping operators. In particular, the symmetric strain increment entering the constitutive update is obtained as
\[
\underline{\Delta \beps}_\mathrm{macro} = \mathcal{S}\,\Delta \beps_\mathrm{macro}\,,
\]
where $\mathcal{S}$ computes the Voigt representation of a  symmetric strain in full array notation.

The corresponding stress returned by the standard constitutive routine, $\std{\underline{\bT}}$, is then embedded into the full (generally non-symmetric) representation through the lifting operator
\[
\std{\bT} = \mathcal{P}\,\std{\underline{\bT}}.
\]

Consistently, the tangent operator in full form is obtained as
\[
\frac{\partial \std{\bT}}{\partial \std{\beps}}
=
\mathcal{P}\,
\frac{\partial \std{\underline{\bT}}}{\partial \std{\underline{\beps}}}\,
\mathcal{S}.
\]
The operators $\mathcal{S}$ and $\mathcal{P}$ are reported in \ref{Apn_FEOper} for the plane strain formulation.
Therefore, any standard constitutive integration algorithm, together with its consistent tangent, can be employed without modification within the present framework.
}.
:
\begin{equation}
    \bT_\mathrm{macro}
    =
\std{\bT}\bigl(\Delta \beps_\mathrm{macro},\,
\beps^\mathrm{e}_{\mathrm{macro},n},\,\std{\bp}\bigr),
\end{equation}
where $\beps^\mathrm{e}_{\mathrm{macro},n}$ are the (standard) elastic strains of the macro-continuum at the beginning of the increment, i.e., at time $t_n$.
This is the central feature of the formulation: both the macro-stress and the corresponding consistent elastoplastic tangent can be obtained directly from \emph{any} standard elastoplastic constitutive routine, which can therefore be employed as a black box without modifying either the stress update algorithm or the tangent operator.
 Hence, the stress $\bT$
can be computed from 
Eq. \eqref{eq:decompstress} as:
\begin{equation}
    \bT = \bT_\mathrm{macro}+G\left(k_1 \bI\otimes \bI+k_2\, \mathcal{I}_d \right) \left(\mathcal{A}\,\hat \bu -\mathcal{N}_\tau\, \hat \boldeta\right).
\label{eq:TImpl}
\end{equation}

The FE stiffness matrix is obtained by differentiating the internal forces $\bF^\mathrm{int}_\mathrm{u}$ and $\bF^\mathrm{int}_\eta$ with respect to the nodal displacements $\hat{\bu}$ and the nodal tensor $\hat \boldeta$, i.e.,
\begin{equation}
    \bK=\begin{bmatrix}
        \bK_\mathrm{uu} & \bK_{\mathrm{u}\eta}\\
        \bK_{\eta \mathrm{u}} & \bK_{\eta \eta} 
    \end{bmatrix}
    \label{eq:defstiff}
\end{equation}
where:
\begin{equation}
\begin{aligned}
   \bK_\mathrm{uu}= \frac{\partial \bF^\mathrm{int}_\mathrm{u}}{\partial \hat \bu}  =\int_P  \mathcal{A}^T  \left[
   \frac{\partial \std{\bT}}
   {\partial \std{\beps}}
\mathcal{A}_\mathrm{sym}
+G\left(k_1 \bI\otimes \bI+k_2\, \mathcal{I}_d \right) \mathcal{A}
\right]
    \,
   \mathrm{d}V\,,
   \end{aligned}
   \label{eq:stiff1}
\end{equation}

\begin{equation}
\begin{aligned}
   \bK_{\mathrm{u}\eta}= \frac{\partial \bF^\mathrm{int}_\mathrm{u}}{\partial \hat \boldeta}  =\int_P   - G\mathcal{A}^T 
   \left(k_1 \bI\otimes \bI+k_2\, \mathcal{I}_d \right)
    \mathcal{N}_\tau \,
   \mbox{d}V\,,
   \end{aligned}
   \label{eq:stiff2}
\end{equation}

\begin{equation}
\begin{aligned}
   \bK_{\eta \mathrm{u}}= \frac{\partial \bF^\mathrm{int}_\eta}{\partial \hat \bu}  =\int_P   -  G\,\mathcal{N}_\tau^T \left(k_1 \bI\otimes \bI+k_2\, \mathcal{I}_d \right)
     \mathcal{A} \,
   \mbox{d}V\,,
   \end{aligned}
   \label{eq:stiff3}
\end{equation}

\begin{equation}
\begin{aligned}
   \bK_{\eta \eta}= \frac{\partial \bF^\mathrm{int}_\eta}{\partial \hat \boldeta}  =\int_P   G \left[\,\mathcal{N}_\tau^T\left(k_1 \bI\otimes \bI+k_2\, \mathcal{I}_d \right) \, \mathcal{N}_\tau+
   4  \ell^2 \left(\mathcal{G}_\mathrm{sym}^T \,\mathcal{G}_\mathrm{sym} \right)
   \right]\,
   \mathrm{d}V\,.
   \end{aligned}
   \label{eq:stiff4}
\end{equation}
It should be noted that, once $\std{\bT}$ and the consistent elastoplastic tangent 
$\partial \std{\bT}/\partial \std{\beps}$ are obtained from a standard constitutive update, 
the entire formulation reduces to the assembly of linear contributions. The complete FE algorithm is summarized in Algorithm \ref{alg:FE}.

\begin{algorithm}
\caption{Algorithm for the FE formulation}
\label{alg:FE}
\begin{algorithmic}
\Statex At the element level, the formulation is expressed in terms of the nodal variables $\hat{\bu}$ and $\hat{\boldeta}$, together with the standard constitutive variables of the macro-continuum at the beginning of the increment.

\State ASSEMBLE the FE operators $\mathcal{A}$, $\mathcal{A}_\mathrm{sym}$, $\mathcal{N}_\tau$, and $\mathcal{G}_\mathrm{sym}$ from the shape functions and their derivatives;

\State COMPUTE the macro (standard) strain increment  $\Delta \beps_\mathrm{macro}$ from the nodal displacement increment using Eq. \eqref{eq:epsstd}:
\begin{equation}
\Delta \beps_\mathrm{macro}= \mathcal{A}_\mathrm{sym}\Delta \hat \bu
\,;\nonumber
\end{equation}
\State CALL the standard elastoplastic routine to compute the response of the macro-continuum: 
\Statex \hspace{\algorithmicindent}GIVEN
\Statex \hspace{\algorithmicindent}\hspace{\algorithmicindent} the (standard) elastic strain $\beps^\mathrm{e}_{\mathrm{macro},n}$ (at the beginning of the increment, i.e., at time $t_n$), 
\Statex \hspace{\algorithmicindent}\hspace{\algorithmicindent} the (standard) state variables $\bp_{\mathrm{macro},n}$ (at the beginning of the increment, i.e., at time $t_n$),  
\Statex \hspace{\algorithmicindent}\hspace{\algorithmicindent} the macro strain increment $\Delta \beps_\mathrm{macro}$, 
\Statex \hspace{\algorithmicindent}COMPUTE
\Statex \hspace{\algorithmicindent}\hspace{\algorithmicindent}  the (standard) macro-stress $\std{\bT}$,
\Statex \hspace{\algorithmicindent}\hspace{\algorithmicindent} the (standard) elastic strain of the macro-continuum $\beps^\mathrm{e}_\mathrm{macro}$,
\Statex \hspace{\algorithmicindent}\hspace{\algorithmicindent}  the (standard) consistent elastoplastic tangent 
      $\partial \std{\bT}/\partial\beps_\mathrm{macro}$,
\Statex \hspace{\algorithmicindent}\hspace{\algorithmicindent} the updated (standard) state variables $\std{\bp}$ of the macro-continuum.

\State COMPUTE the (non-symmetric) director stress $\bT_\mathrm{micro}$ and the micro-couples $\bM$ using Eq. \eqref{eq:nonstdstress}:
\begin{equation}
\bT_\mathrm{micro}=G\left(k_1 \bI\otimes \bI+k_2\, \mathcal{I}_d \right) 
\left( \mathcal{A} \hat \bu-\mathcal{N}_\tau \hat \boldeta\right),\quad
\bM=4 G \ell^2\, \mathcal{G}_\mathrm{sym}\, \hat \boldeta\,;
\nonumber
\end{equation}

\State COMPUTE the internal forces $\bF^{\mathrm{int}}_\mathrm{u}$ and $ \bF^{\mathrm{int}}_{\eta}$ using Eq. \eqref{eq:residuals}:
\begin{equation}
    \bF^{\mathrm{int}}_\mathrm{u}=\int_P   \mathcal{A}^T \left( \std{\bT}+\bT_\mathrm{micro}\right)\,\mathrm{d}V\,,\quad
    \bF^{\mathrm{int}}_{\eta}=\int_P \left(\mathcal{G}_\mathrm{sym}^T  \bM - \mathcal{N}_\tau^T \bT_\mathrm{micro}\right)
\mathrm{d}V \,;
\nonumber\end{equation}

\State COMPUTE the FE stiffness matrix $\bK$ using Eqs. \eqref{eq:defstiff}-\eqref{eq:stiff4}:
\begin{equation}
    \bK=\begin{bmatrix}
        \bK_\mathrm{uu} & \bK_{\mathrm{u}\eta}\\
        \bK_{\eta \mathrm{u}} & \bK_{\eta \eta} 
    \end{bmatrix}\,,
\nonumber
\end{equation}
where:
\begin{equation}
\begin{aligned}
   \bK_\mathrm{uu}=\int_P  \mathcal{A}^T  \left[
   \frac{\partial \std{\bT}}
   {\partial \std{\beps}}
\mathcal{A}_\mathrm{sym}
+G\left(k_1 \bI\otimes \bI+k_2\, \mathcal{I}_d \right) \mathcal{A}
\right]
    \,
   \mathrm{d}V\,,
   \end{aligned}
   \quad
\begin{aligned}
   \bK_{\mathrm{u}\eta}= \int_P   - G\mathcal{A}^T 
   \left(k_1 \bI\otimes \bI+k_2\, \mathcal{I}_d \right)
    \mathcal{N}_\tau \,
   \mbox{d}V\,,
   \end{aligned}
 \nonumber
\end{equation}
\begin{equation}
\begin{aligned}
   \bK_{\eta \mathrm{u}} =\int_P   -  G\,\mathcal{N}_\tau^T \left(k_1 \bI\otimes \bI+k_2\, \mathcal{I}_d \right)
     \mathcal{A} \,
   \mbox{d}V\,,
   \end{aligned}
  \quad
\begin{aligned}
   \bK_{\eta \eta}=\int_P   G \left[\,\mathcal{N}_\tau^T\left(k_1 \bI\otimes \bI+k_2\, \mathcal{I}_d \right) \, \mathcal{N}_\tau+
   4  \ell^2 \left(\mathcal{G}_\mathrm{sym}^T \,\mathcal{G}_\mathrm{sym} \right)
   \right]\,
   \mathrm{d}V\,.
   \end{aligned}
\nonumber
\end{equation}

\Statex The above algorithm is written in a form valid for a general 3D stress state. The FE operators $\mathcal{A}$, $\mathcal{A}_\mathrm{sym}$, $\mathcal{N}_\tau$, and $\mathcal{G}_\mathrm{sym}$ are detailed in  \ref{Apn_FEOper} for the 2D plane strain  case.

\end{algorithmic}
\end{algorithm}

\section{Benchmark problem: shallow strip footing}

A shallow rigid strip footing has been analyzed under plane strain conditions by employing either the Tresca or the Matsuoka--Nakai \cite{matsuoka1974stress,PL2014} plasticity models. Whilst the former is typically adopted to reproduce the behaviour of soils in undrained conditions, the latter class of models is more appropriate for frictional geomaterials in drained regimes. 

The considered benchmark is well known to be particularly demanding from a numerical standpoint \cite{DeBorst1984}, as it involves the development of highly localized deformation patterns beneath the footing, often associated with the interaction of multiple failure mechanisms. Such features make it especially suitable for assessing the ability of a numerical formulation to capture strain localization in a mesh-objective manner.

The footing has a total width $2B=2\,\mathrm{m}$ and is subjected to vertical and centered loading conditions. The out-of-plane thickness $h$ is equal to $1\,\mathrm{m}$. Owing to symmetry, only half of the domain is modeled, corresponding to a footing width $B=1\,\mathrm{m}$. The soil domain extends over a $50\,\mathrm{m}\times50\,\mathrm{m}$ region, which is sufficiently large to avoid any spurious influence of the boundary conditions on the computed response. 

Boundary conditions are prescribed as follows. Along the axis of symmetry, the horizontal displacement is constrained, $u_1=0$, and the director component $\eta_{21}$ is set equal to zero consistently with the symmetry of the director field. No other essential boundary conditions are imposed on the director field. Along the right boundary, only the horizontal displacement is constrained, whereas along the bottom boundary only the vertical displacement is constrained.
The rigid strip footing is modeled by prescribing the vertical displacement of the soil nodes located beneath the footing. The horizontal displacement of the same nodes is also constrained, $u_1=0$, so as to prevent tangential relative motion between the soil and the footing.

Depending on the constitutive model and loading scenario, the analyses are performed either starting from an unstressed initial configuration or from an initial geostatic state induced by self-weight. The latter case is particularly relevant for frictional soils, where the stress field generated by gravity plays a crucial role in the subsequent development of failure mechanisms. The specific assumptions adopted in each case are detailed in the following subsections.

The mechanical behaviour of the macro-continuum is described by an isotropic elastoplastic model with linear elasticity and associated plasticity, based on the General Yield Criterion proposed in \cite{LP2016}. The yield function is expressed as
\begin{equation}
   f= q\,\Gamma_f(\theta)-M(\bar \varepsilon_p)\,p-\kappa(\bar \varepsilon_p),
\end{equation}
where $q=\sqrt{\frac{3}{2}\bT_\mathrm{macro}^\prm:\bT_\mathrm{macro}^\prm}$ denotes the equivalent von Mises macro-stress, $p=-(\bT_\mathrm{macro}:\bI)/3$ is the mean macro-stress, and
\begin{equation}
    \theta=\frac{1}{3}\arcsin\left(-\frac{27}{2}\frac{\det \bT_\mathrm{macro}^\prm }{q^3}\right)
\end{equation}
is the Lode angle. The scalar variable $\bar \varepsilon_p$ represents a measure of accumulated plastic strain.

The function $\Gamma_f(\theta)$ controls the shape of the yield surface in the octahedral plane and is defined as
\begin{equation}
    \Gamma_f(\theta)=a_f \cos \left[\frac{\arccos\left(-b_f\sin3\theta \right)}{3}-c_f\frac{\pi}{6}\right].
\end{equation}
By appropriate selection of the parameters $a_f$, $b_f$, and $c_f$, this formulation can reproduce a variety of classical yield criteria, including  rounded Tresca and Matsuoka--Nakai \cite{matsuoka1974stress,PL2014}, as well as Lade--Duncan \cite{Lade1974} and von Mises.

\begin{figure}
\centering
\begin{tabular}{cc}
		\includegraphics[trim={0cm 0cm 12cm 0cm},clip,width=0.25\textwidth]{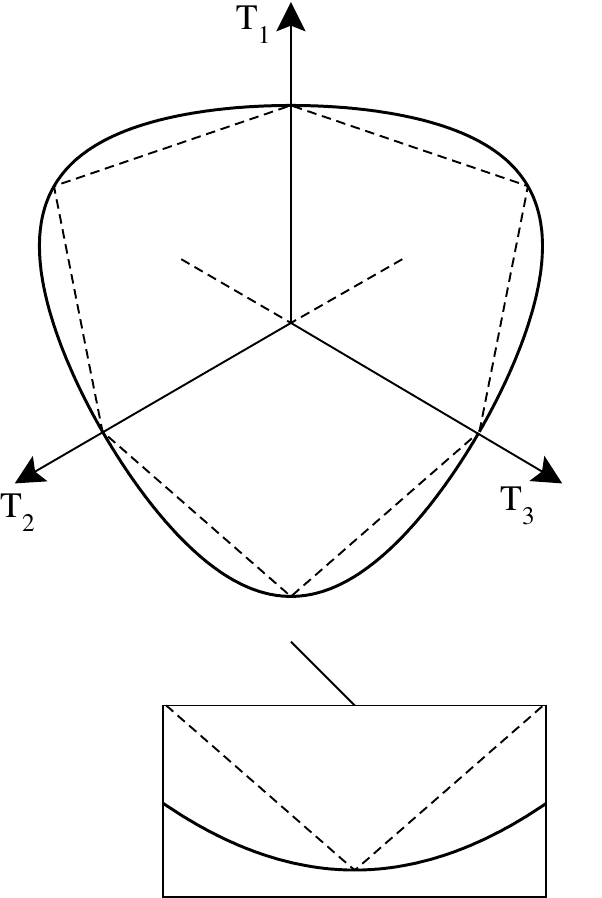}
        \hspace{0.5cm}
        &
        \hspace{0.5cm}
        \includegraphics[trim={12cm 0cm 0cm 0cm},clip,width=0.25\textwidth]{figures/mnakai_phi20.pdf}
        \\
        (a) & (b)
        \end{tabular}
    \caption{Plot of the adopted yield criteria in the octahedral plane.  (a) The solid line is the Matsuoka-Nakai criterion. For comparison the dashed line is the Mohr-Coulomb yield criterion for the same shear resistance angle $\phi=20^\circ$; (b) The solid line is the Rounded Tresca criterion based on the parameters reported in Table \ref{tab:shape}. For comparison the dashed line is the sharp Tresca criterion and the dotted line is the von Mises circle}.
    \label{fig:shapefunction}
\end{figure}

In the case of criteria exhibiting corners in the octahedral plane, the adopted parameters smooth the apexes while preserving convexity of the yield surface. In the present study, the shape of the octahedral section is kept constant throughout the analysis, i.e. the parameters $a_f$, $b_f$, and $c_f$ do not evolve with plastic deformation. The values adopted to reproduce the rounded Tresca and Matsuoka--Nakai criteria are reported in Table~\ref{tab:shape}, and the corresponding shape functions are shown in Fig.~\ref{fig:shapefunction}.

\begin{table}[t]
 	\centering
 	\begin{tabular}{c|c c c}
        \hline
 		\hline  		
       Model    &   $a_f$ &    $b_f$ & $c_f$ \\ 
       \hline
        Outer smooth Tresca  &$1.151579$& $0.9999$& $1.0$\\
         Matsuoka-Nakai, $\phi=30^\circ$ &  $1.442221$ &$0.746712$& $0.0$\\
        Matsuoka-Nakai, $\phi=20^\circ$ &  $1.328450$ &$0.552093$& $0.0$\\
        \hline
        \end{tabular}
        \caption{Material parameters adopted for the shape function $\Gamma_f$, with $\phi$ being the shear resistance angle.}
        \label{tab:shape}
\end{table}

The material parameter $M(\bar \varepsilon_p)$, defined as a function of the shear resistance angle $\phi$ as
\begin{equation}
    M(\bar \varepsilon_p)=\frac{6 \sin \phi(\bar \varepsilon_p)}{3-\sin\phi(\bar \varepsilon_p)}\,,
    \label{eq:defM}
\end{equation}
represents the slope of the meridional section in triaxial compression conditions ($\theta=\pi/6$), whereas $\kappa(\bar \varepsilon_p)$ defines the intercept of the meridional section with the deviatoric stress axis in triaxial compression. Their evolution laws, which are responsible for the onset of softening in the considered problems, are specified separately for each constitutive model in the following subsections.

The numerical integration of the elastoplastic constitutive model is fully described in \cite{PL2014, PL2018}, and is employed here without any modification, consistently with the black-box nature of the proposed formulation.

\begin{figure}[t]
	\centering
		\includegraphics[width=0.9\textwidth]{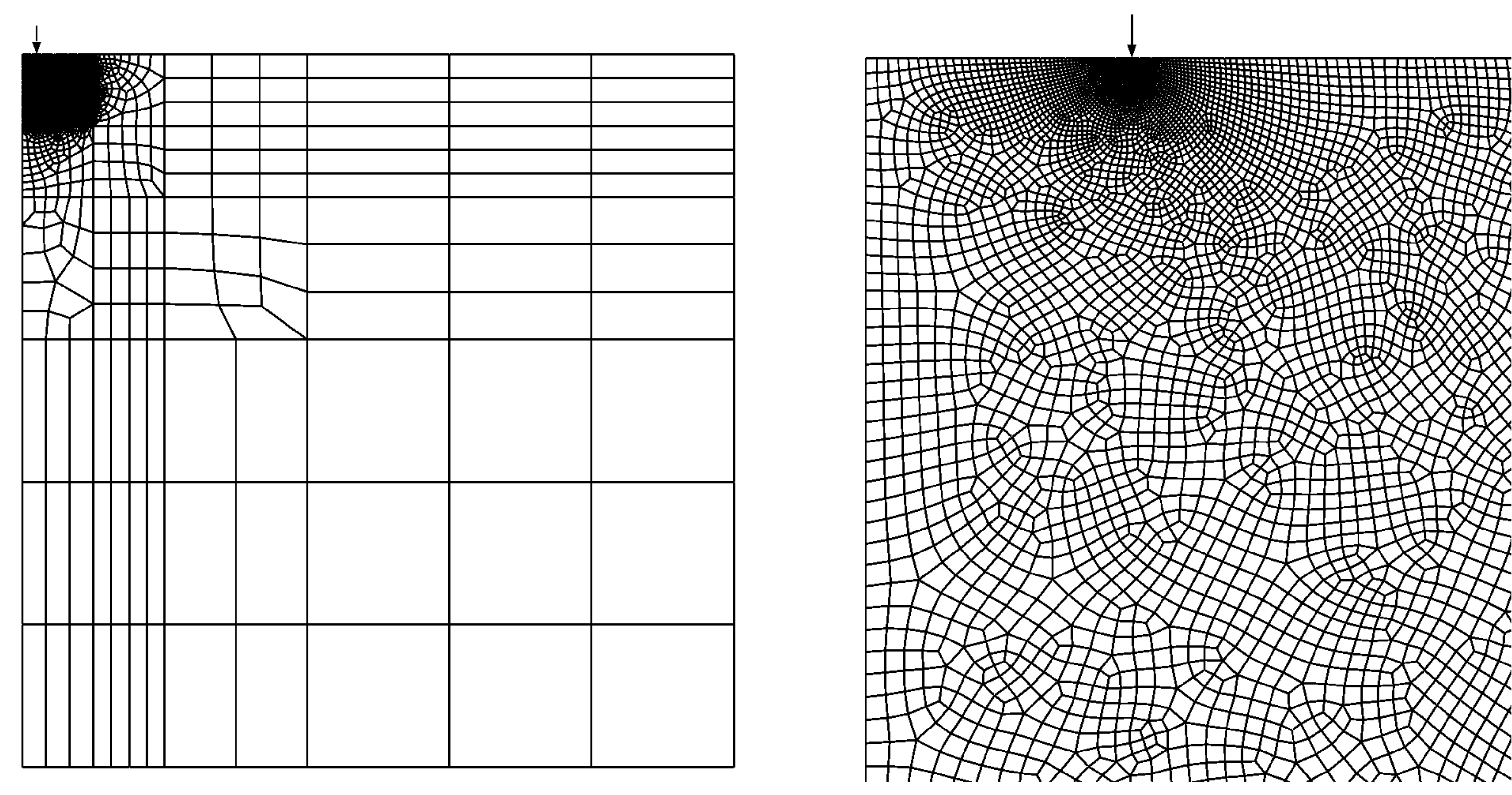} 
\caption{Mesh 1 (coarse mesh). The arrow indicates the edge of the shallow footing.}
\label{fig:mesh}
\end{figure}

Four different meshes have been employed. The coarsest one is depicted in Fig.~\ref{fig:mesh}. Each subsequent mesh (Mesh 2, Mesh 3, Mesh 4) has been obtained by subdividing each element of the previous mesh into four elements. This deliberate choice avoids introducing preferred mesh orientations and allows for a systematic assessment of the mesh objectivity of the formulation.

The adopted finite elements are quadratic in the displacement field $\bu$ and linear in the director field $\boldeta$ (i.e., $n=8$, $n_\eta=4$). The internal force vector and the FE stiffness matrix are evaluated using four Gauss integration points. The main properties of the adopted meshes are summarized in Table~\ref{tab:meshfooting}.

\begin{table}[t]
    \centering
    \begin{tabular}{c c c c c}
    \hline\hline
        Mesh     &  Number of & Number of  & Characteristic                             & Deformable\\
                 &  Elements  &   Nodes    & Element Length $\ell_e$          &   Cosserat DoFs\\
         \hline
         $1$   &   $7,916$    &  $24,028$    &  $13.4 \cdot 10^{-2}$ m &                  $80,728$        \\
         $2$  &    $31,664$   &  $95,550$    & $6.715 \cdot 10^{-2}$ m &                 $318,870$       \\
         $3$  &    $126,656$   &  $381,082$   & $3.358 \cdot 10^{-2}$ m &                 $1,271,014$      \\       
         $4$  &    $506,624$   &  $1,522,098$   & $1.679 \cdot 10^{-2}$ m &             $5,075,142$\\       
         \hline
         \hline
    \end{tabular}
    \caption{Meshes for the strip footing problem}
    \label{tab:meshfooting}
\end{table}

All analyses have been performed using the commercial finite element code Abaqus \cite{ABA24}. The proposed finite element formulation has been implemented through a user-defined element (\texttt{UEL}) subroutine.

\subsection{Tresca cohesive soil}

The  elastic response of the soil under undrained conditions is modeled by assuming a shear modulus $G=416.7\,\mathrm{MPa}$ and a bulk modulus $K_v=55560\,\mathrm{MPa}$, corresponding to an elastic modulus $E=1247\,\mathrm{MPa}$ and a Poisson ratio $\nu=0.4963$. Since the undrained plastic response is pressure-independent, the parameter $M$ is set equal to zero.

The accumulated plastic strain is defined as
\begin{equation}
 \bar \varepsilon_p=\int \sqrt{\frac{2}{3}\bD_p^\prm:\bD_p^\prm}\,\mathrm{d}t.
\end{equation}

Two constitutive responses are considered for the macro-continuum, namely perfect plasticity, $\kappa=\kappa_0$, and isotropic softening, defined as
\begin{equation}
    \kappa=\kappa_\infty+(\kappa_0-\kappa_\infty)\exp(-a_h \bar \varepsilon_p).
    \nonumber
\end{equation}
The softening parameters are chosen as $\kappa_0=980\,\mathrm{kPa}$, $\kappa_\infty=9.8\,\mathrm{kPa}$ and $a_h=10$.

All analyses are carried out under displacement control using the standard nonlinear static solver of Abaqus~\cite{ABA24}. The initial and maximum time increments are set equal to $10^{-4}\,t_0$ and $5\cdot10^{-3}\,t_0$, respectively, where $t_0$ is the (dummy) load time. The imposed normalized vertical displacement is $\bar u/B=0.1$. No geostatic initial state and no lateral surcharge are considered, and a rough footing--soil interface is assumed. Unless otherwise specified, the micro-continuum parameters are set to $k_1=k_2=0.1$.

\begin{figure}[t]
	\centering
		\includegraphics[width=0.47\textwidth]{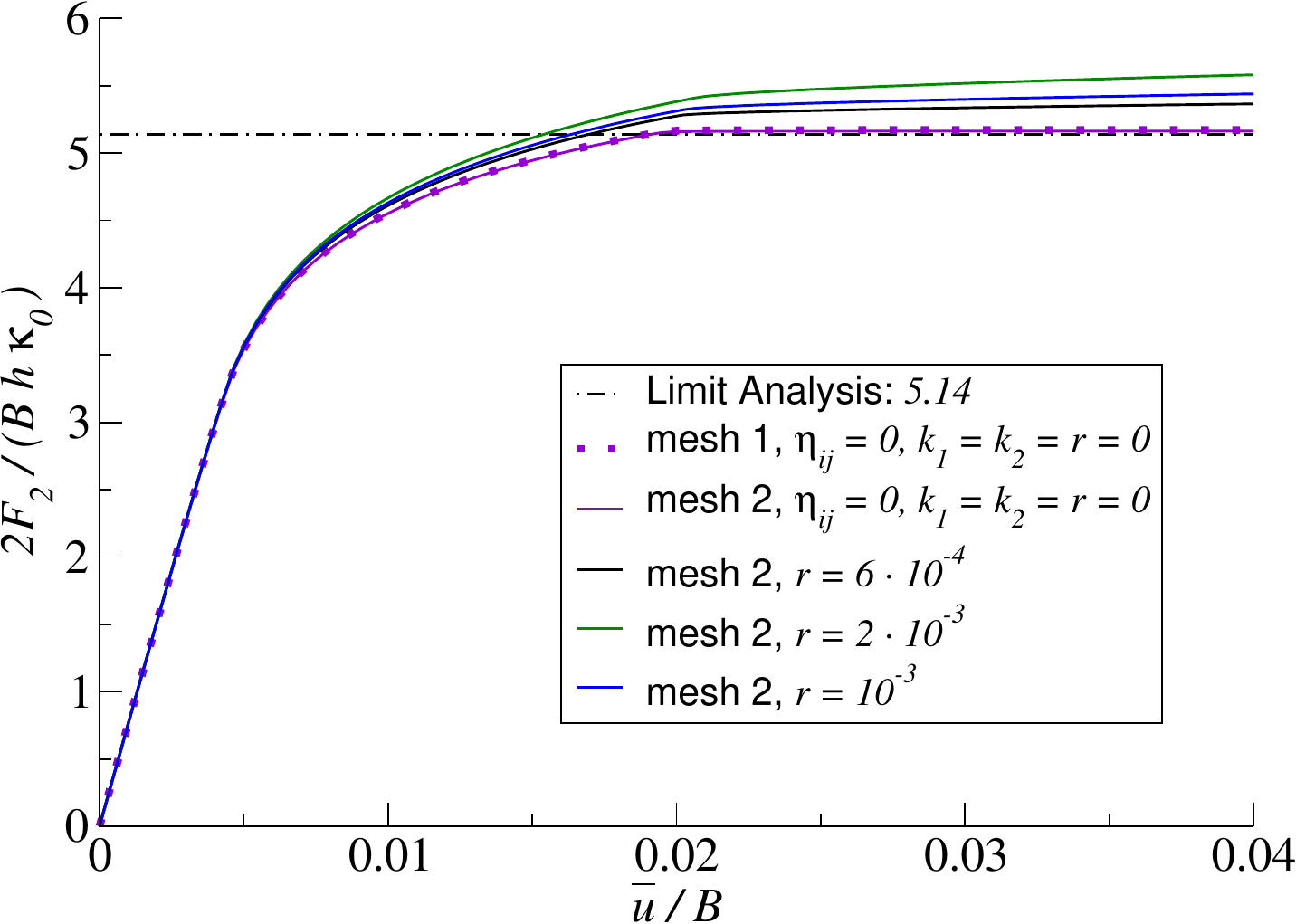} 
\caption{Strip footing on a Tresca soil under perfect plasticity: normalized load--displacement curves. The Cauchy response (obtained by enforcing $\boldeta=\bzero$ and $k_1=k_2=\ell=0$) coincides with Prandtl's solution. Results for the deformable Cosserat continuum are shown for $r=\ell/B=6\cdot10^{-4}$, $10^{-3}$, and $2\cdot10^{-3}$, highlighting the recovery of the classical solution and the emergence of a mild size effect as $r$ increases.}
\label{fig:footing_tresca_pp}
\end{figure}

Fig.~\ref{fig:footing_tresca_pp} shows the normalized load--displacement curves obtained under perfect plasticity. In this case, only Mesh~1 and Mesh~2 are considered, since no mesh dependence is expected in the absence of softening. The classical Cauchy solution is recovered by enforcing $\boldeta=\bzero$ at all nodes and setting $k_1=k_2=\ell=0$, and coincides with Prandtl's solution~\cite{Prandtl1920}. The deformable Cosserat results are obtained for $r=\ell/B=6\cdot10^{-4}$, $10^{-3}$, and $2\cdot10^{-3}$, showing that the proposed formulation correctly recovers the classical response in the limit case. Increasing $r$ leads to a slightly stiffer response and a modest increase in the peak load, indicating the emergence of a size effect when the internal length becomes non-negligible with respect to the footing width.

\begin{figure}[t]
	\centering
	\begin{tabular}{cc}
		\includegraphics[width=0.47\textwidth]{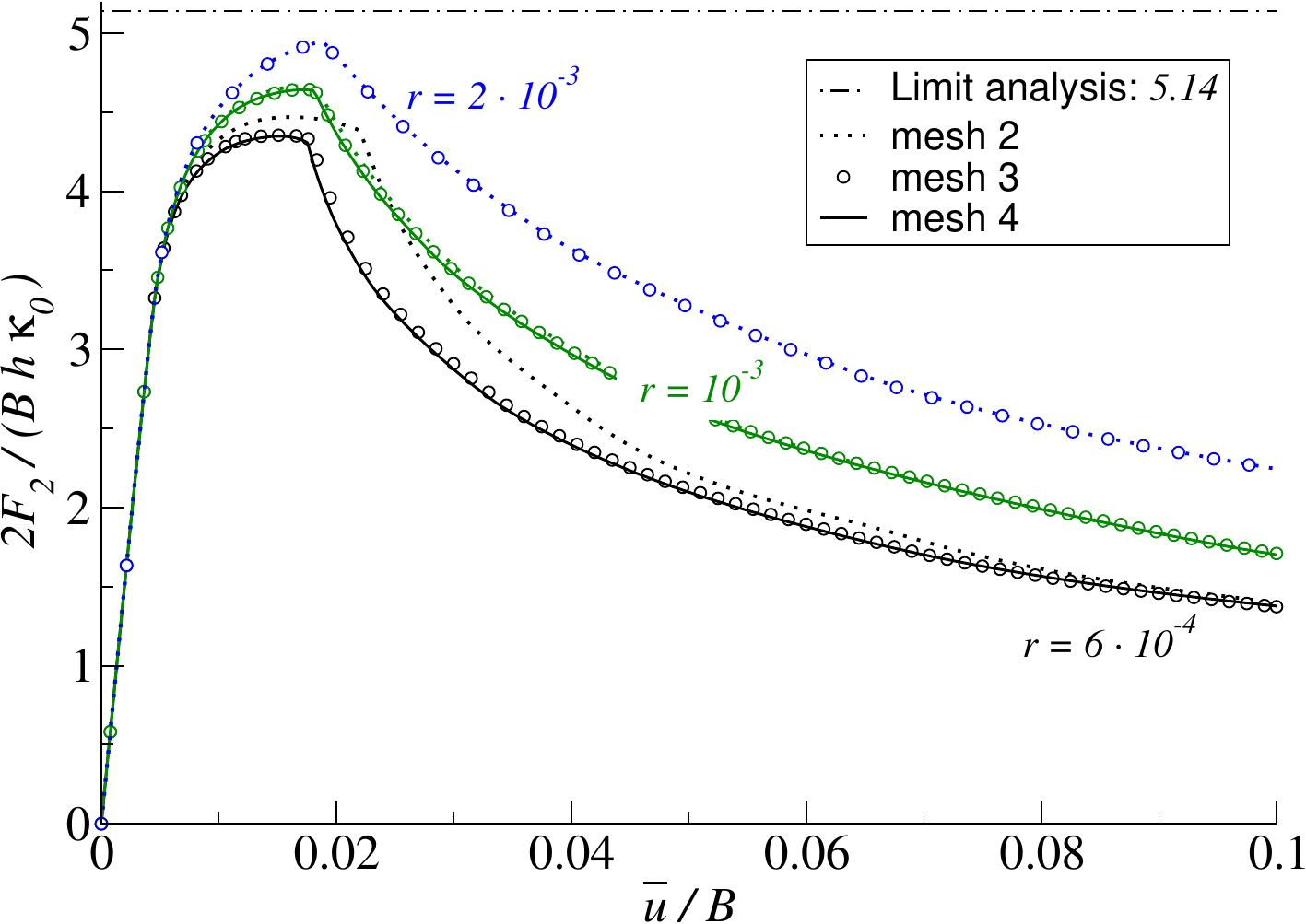} 
		&\includegraphics[width=0.47\textwidth]{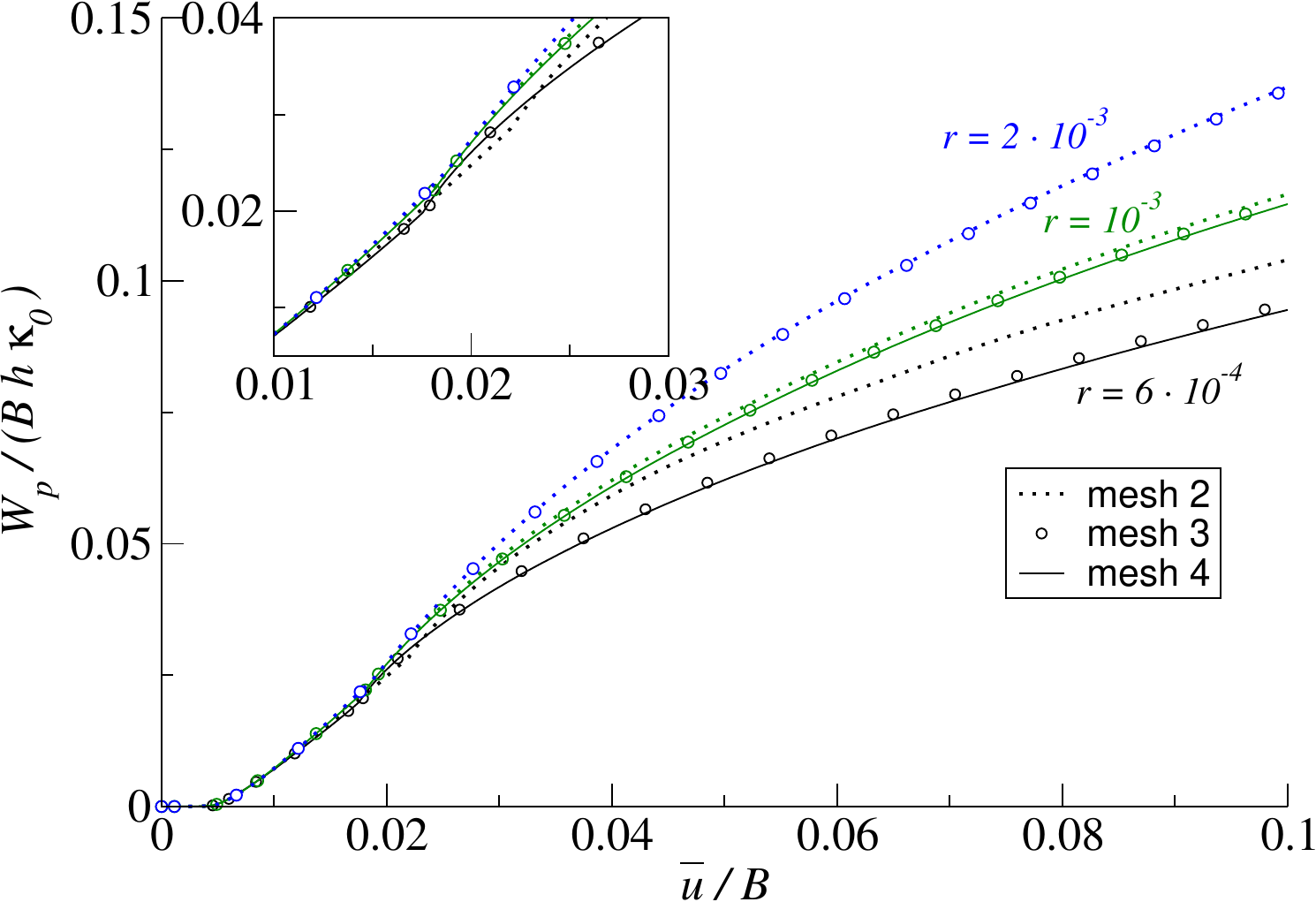}\\
		(a) & (b)
	\end{tabular}
\caption{Strip footing on a Tresca soil with exponential isotropic softening: (a) normalized load--displacement curves and (b) total dissipated energy as functions of the applied displacement. Results are reported for different values of $r=\ell/B$ and for multiple meshes, showing convergence of the both global response and the dissipation history upon mesh refinement.}
\label{fig:footing_tresca_softening}
\end{figure}

The results obtained with exponential softening are reported in Fig.~\ref{fig:footing_tresca_softening}, in terms of load--displacement curves and total dissipated energy as functions of the applied displacement. A clear convergence of the global response is observed upon mesh refinement, both for the load--displacement curves and for the total dissipated energy. The convergence of the entire dissipation history (and not only of its final value) further supports the consistency of the numerical results. The rate of convergence depends on the ratio $r=\ell/B$: for $r=2\cdot10^{-3}$, Mesh~2 and Mesh~3 already provide essentially coincident results, and Mesh~4 is therefore not considered; for $r=10^{-3}$, convergence is very good for all the three meshes, while for $r=6\cdot10^{-4}$, Mesh~2 is not yet converged and Mesh~3 and Mesh~4 are required to approach convergence. This behavior reflects the fact that smaller values of $r$ lead to finer localization patterns, which require a higher spatial resolution.

\begin{figure}[t]
	\centering
	\begin{tabular}{cc}
		\includegraphics[trim={2cm 2cm 10cm 0.5cm},clip,width=0.4\textwidth]{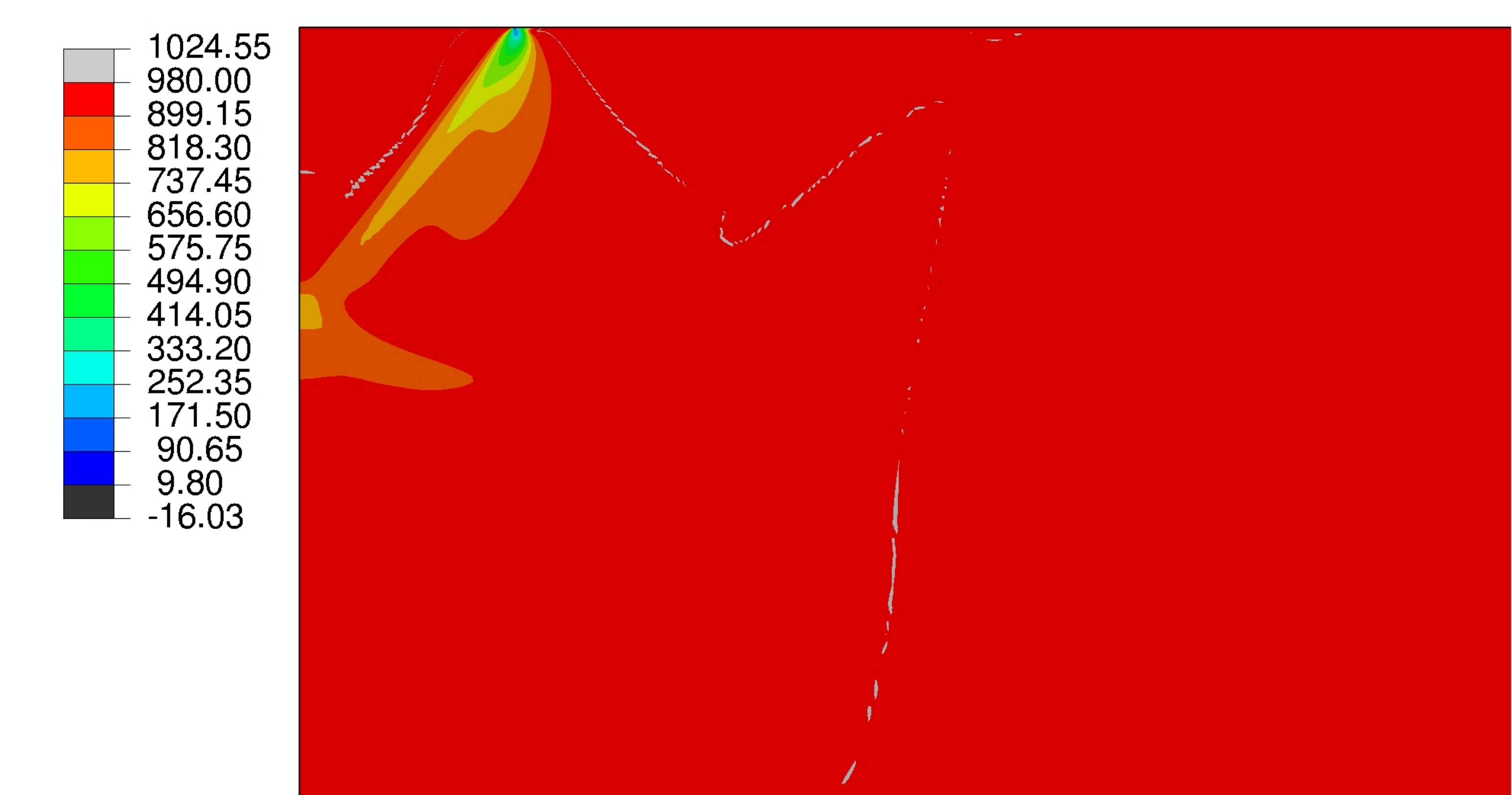} 
		&\includegraphics[trim={2cm 2cm 10cm 0.5cm},clip,width=0.4\textwidth]{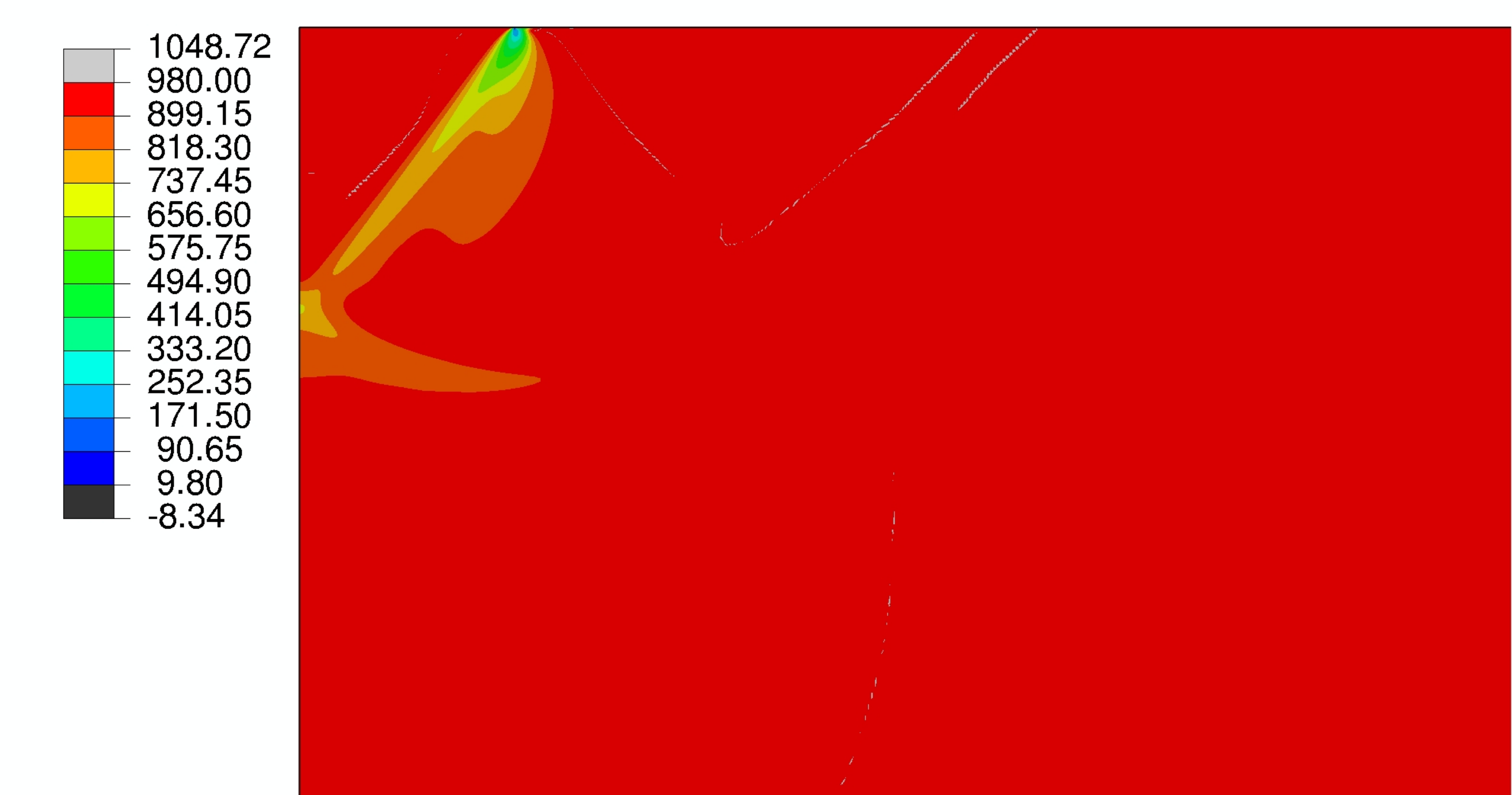} \\
		(a) & (b)\\
		\includegraphics[trim={2cm 2cm 10cm 0.5cm},clip,width=0.4\textwidth]{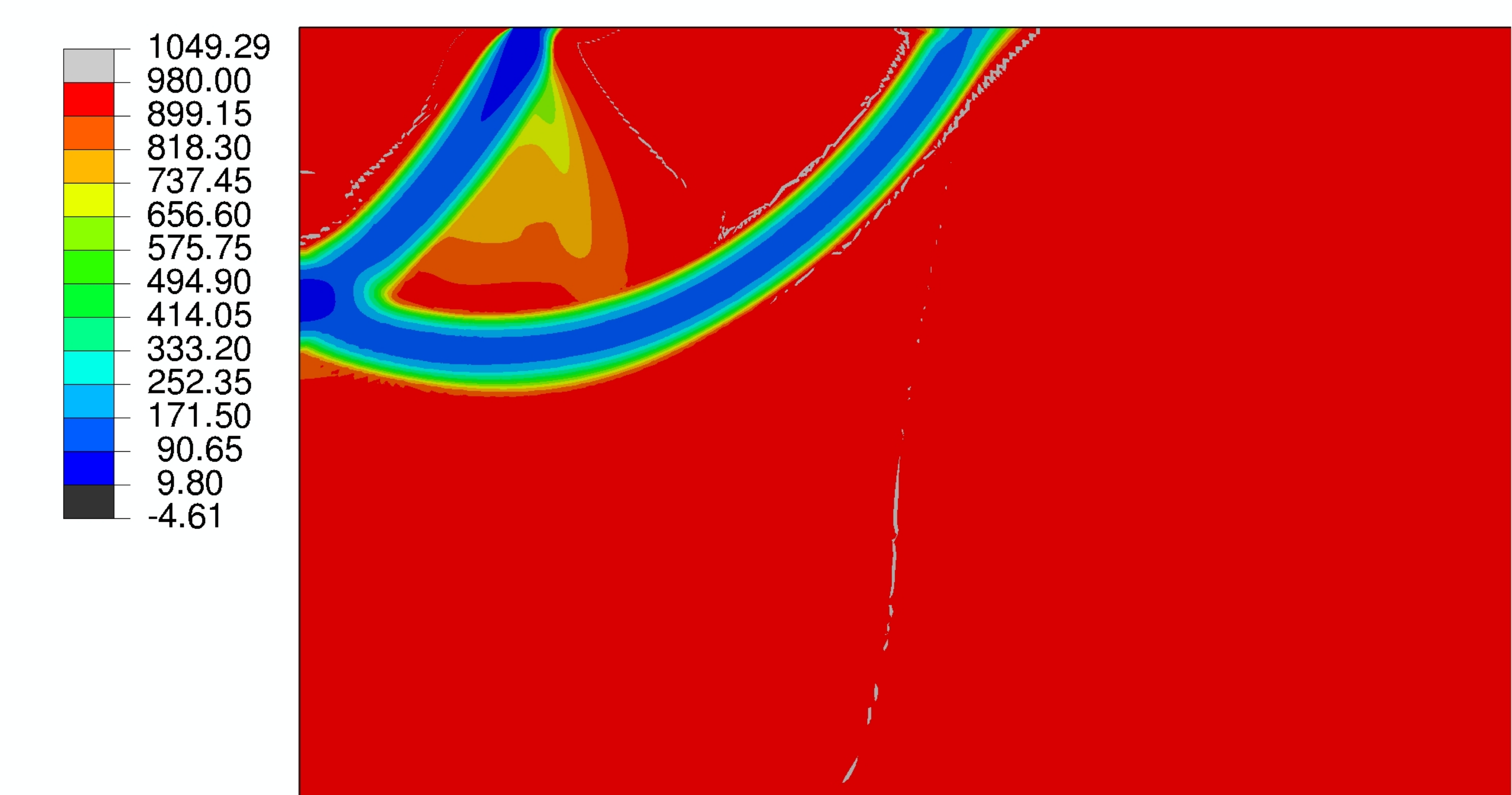} 
		&\includegraphics[trim={2cm 2cm 10cm 0.5cm},clip,width=0.4\textwidth]{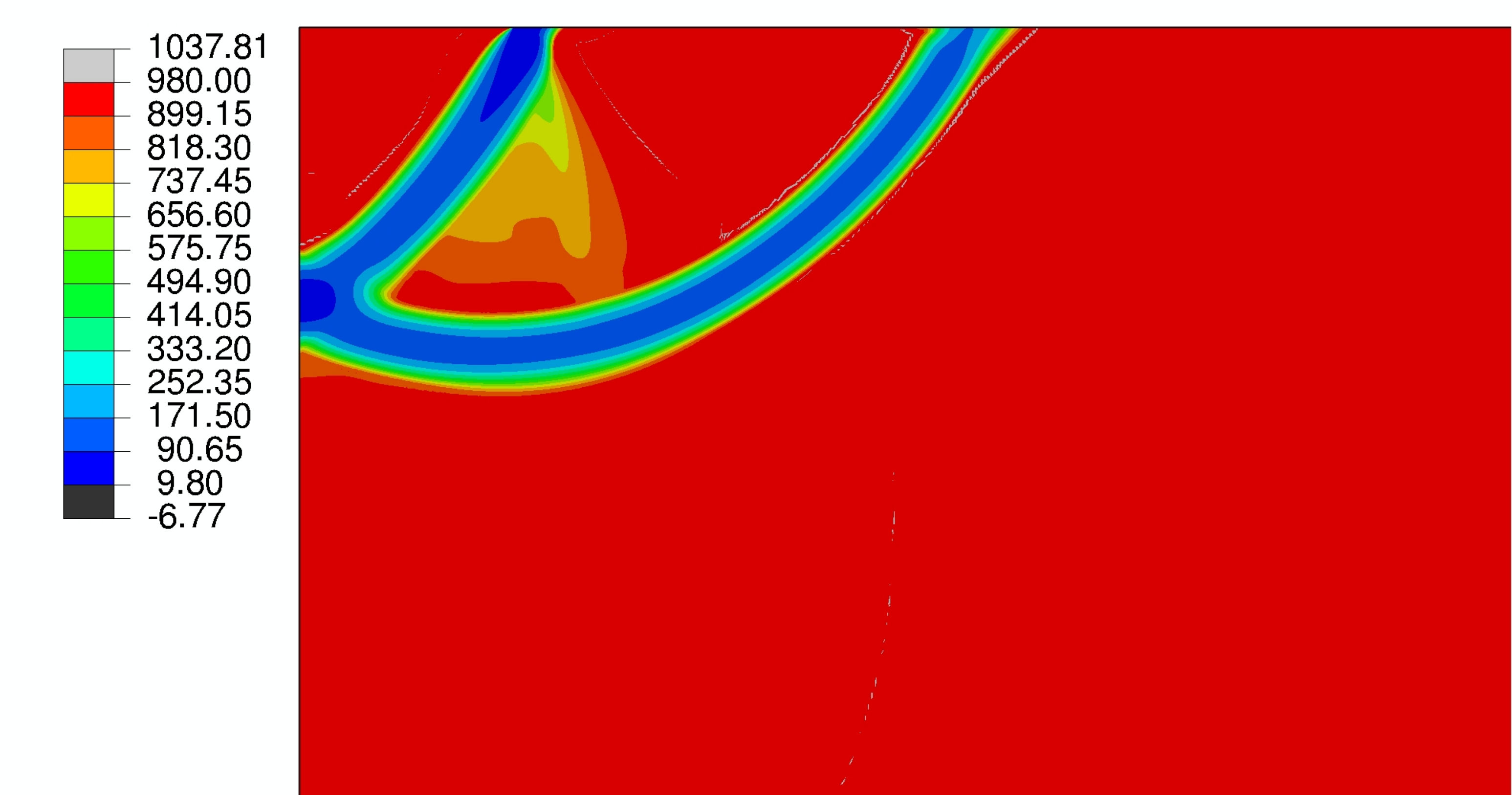} \\
		(c) & (d)
        \end{tabular}
\caption{Spatial distribution of $\kappa$ for the strip footing on Tresca soil with exponential softening, $r=2\cdot10^{-3}$. (a)--(b) Mesh~2 and Mesh~3 at peak load; (c)--(d) Mesh~2 and Mesh~3 at the end of the analysis ($\bar u/B=0.1$). The localization pattern, consisting of two dominant shear bands with a slightly diffuse intermediate zone, is essentially coincident across meshes, indicating full mesh objectivity.}
\label{fig:tresca_conv_r2e-3}
\end{figure}

\begin{figure}[t]
	\centering
	\begin{tabular}{cc}
		\includegraphics[trim={2cm 2cm 10cm 0.5cm},clip,width=0.4\textwidth]{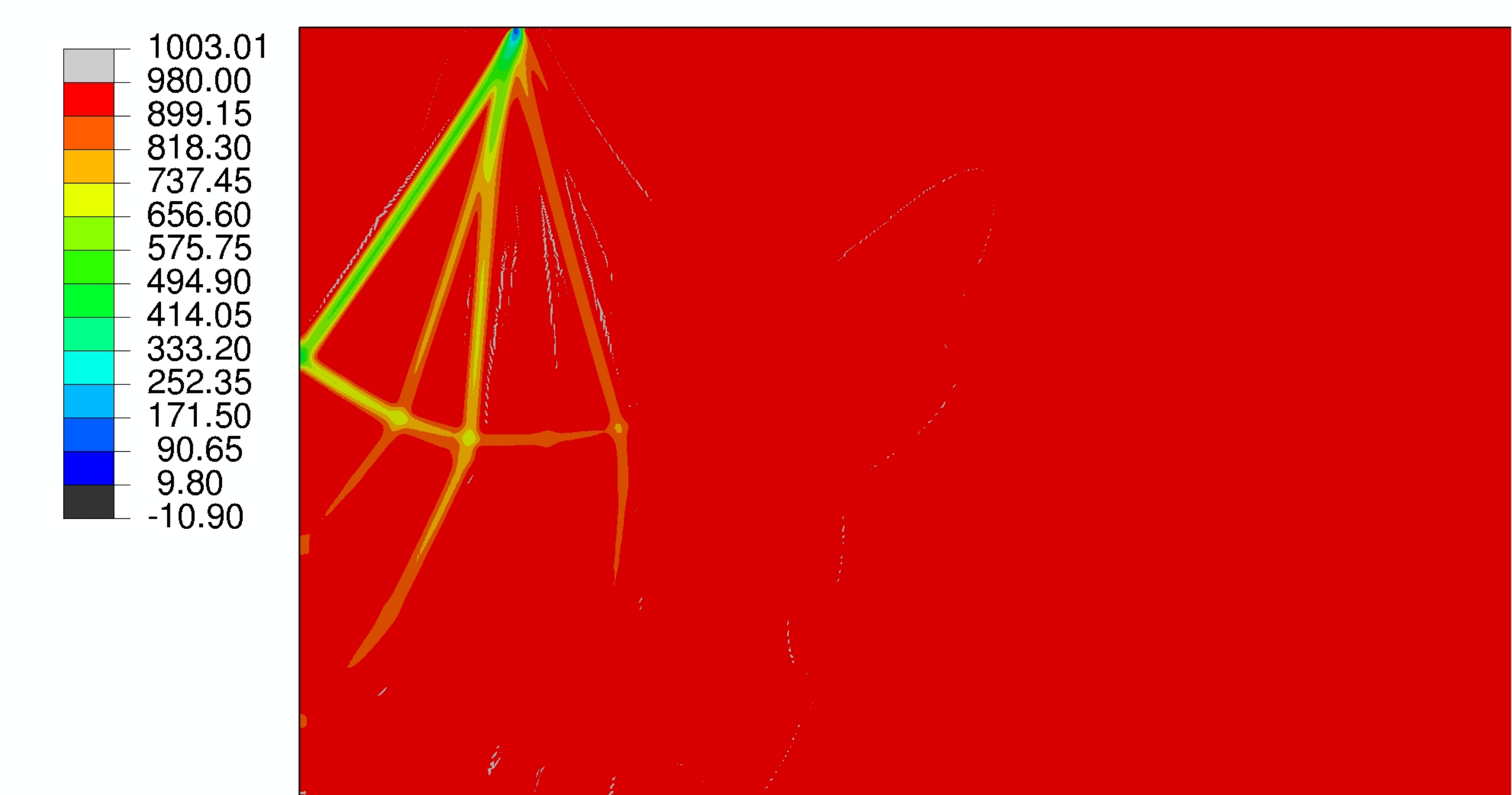} 
		&\includegraphics[trim={2cm 2cm 10cm 0.5cm},clip,width=0.4\textwidth]{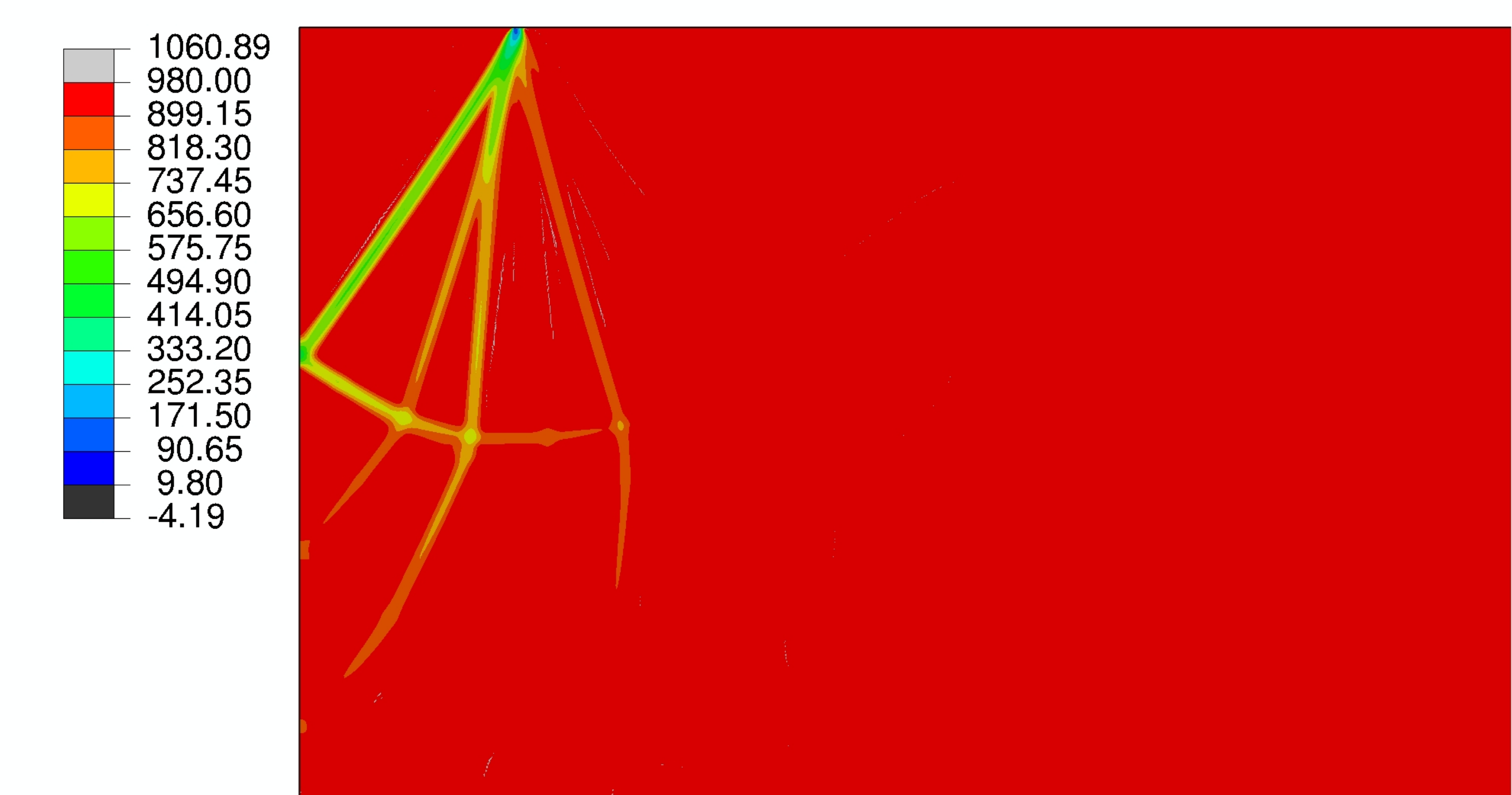} \\
		(a) & (b)\\
		\includegraphics[trim={2cm 2cm 10cm 0.5cm},clip,width=0.4\textwidth]{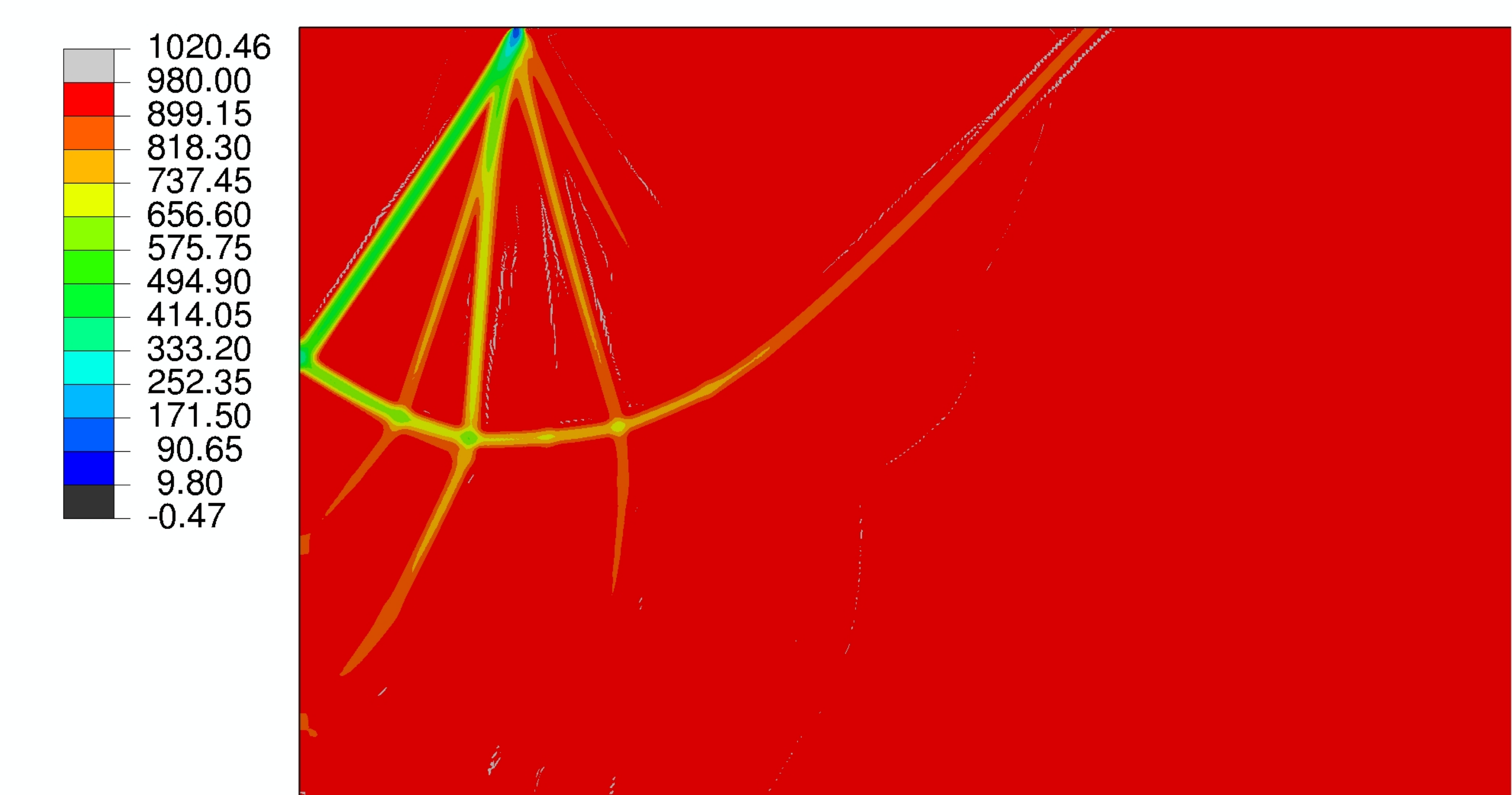} 
		&\includegraphics[trim={2cm 2cm 10cm 0.5cm},clip,width=0.4\textwidth]{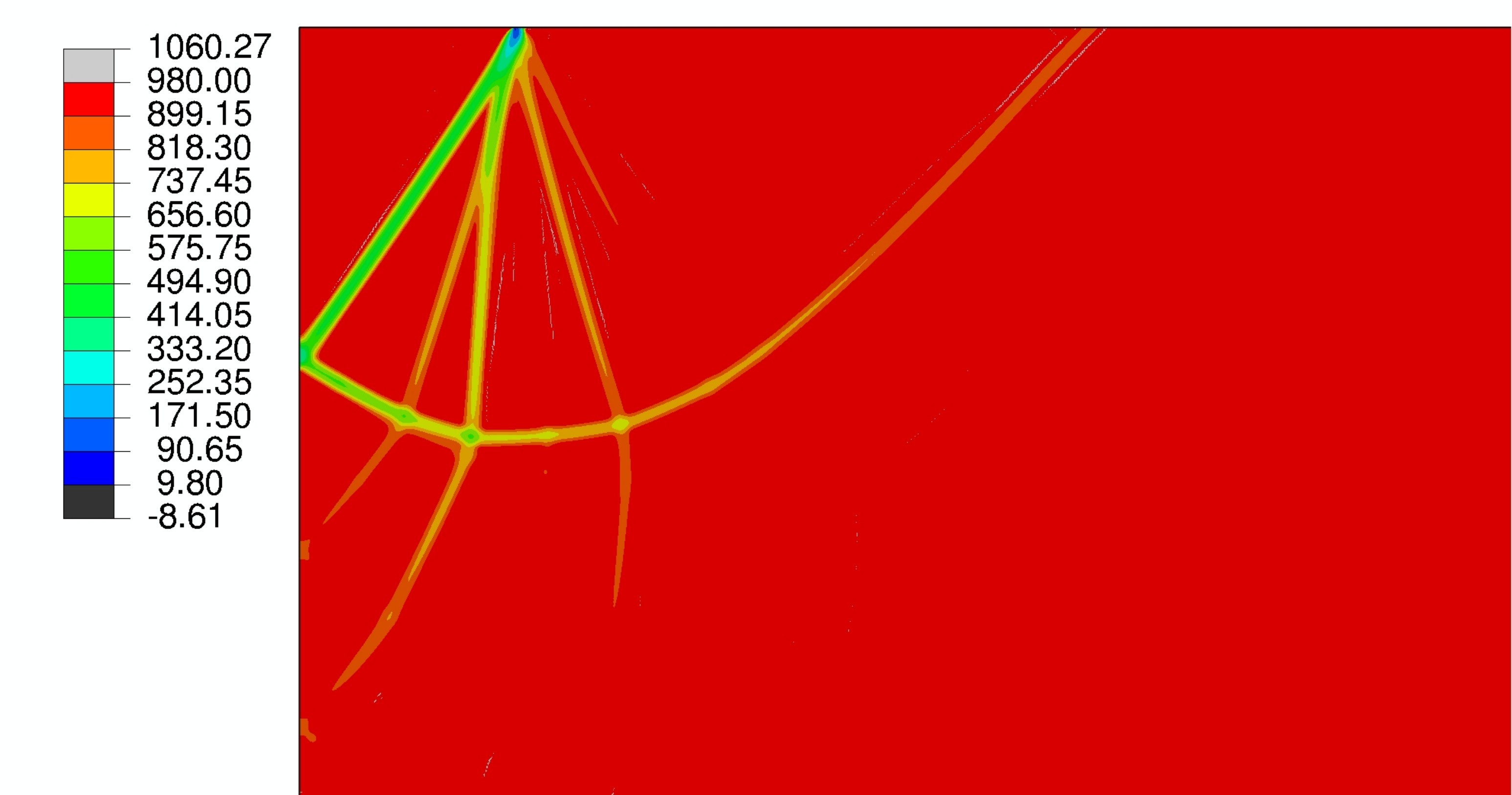} \\
		(c) & (d) \\
		\includegraphics[trim={2cm 2cm 10cm 0.5cm},clip,width=0.4\textwidth]{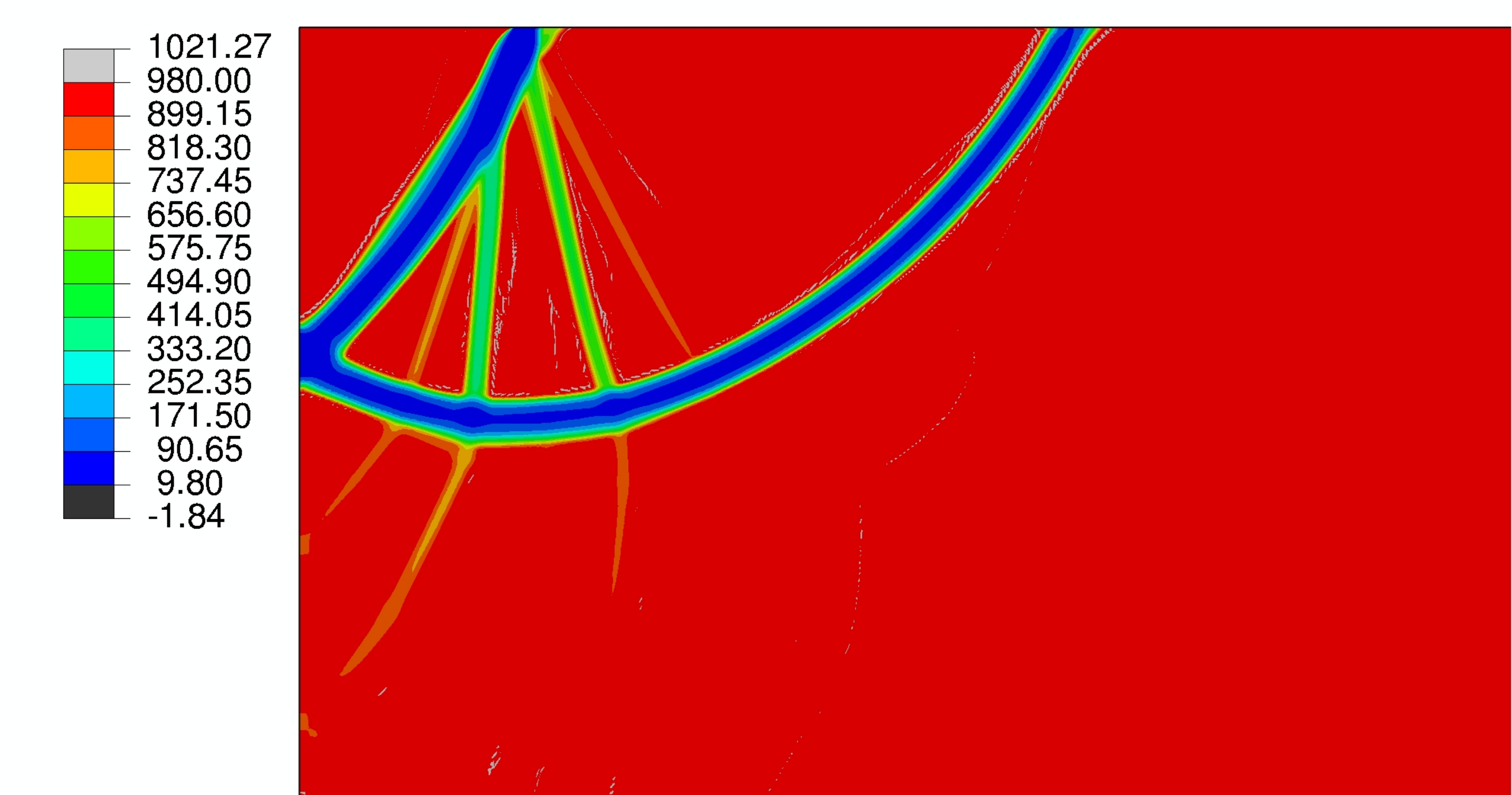} 
		&\includegraphics[trim={2cm 2cm 10cm 0.5cm},clip,width=0.4\textwidth]{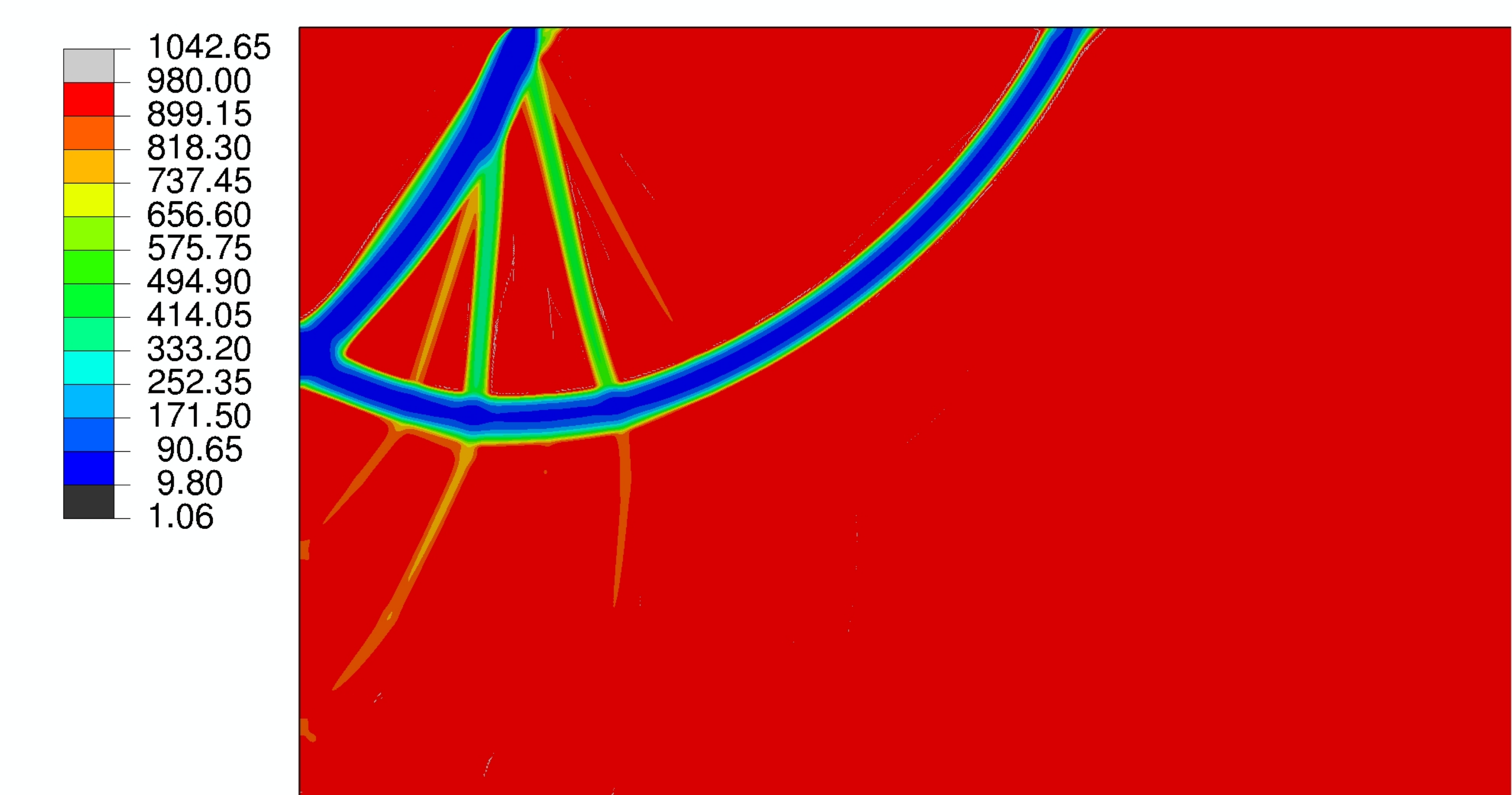} \\
		(e) & (f)
        \end{tabular}
\caption{Spatial distribution of $\kappa$ for the strip footing on Tresca soil with exponential softening, $r=10^{-3}$. (a)--(b) Mesh~3 and Mesh~4 at peak load; (c)--(d) at $\bar u/B=0.02$; (e)--(f) at the end of the analysis ($\bar u/B=0.1$). In addition to the main shear bands, a network of secondary bands develops. Both primary and secondary structures are nearly coincident across meshes, demonstrating mesh-independent prediction of complex localization patterns.}
\label{fig:tresca_conv_r1e-3}
\end{figure}

\begin{figure}[t]
	\centering
	\begin{tabular}{cc}
		\includegraphics[trim={2cm 2cm 10cm 0.5cm},clip,width=0.4\textwidth]{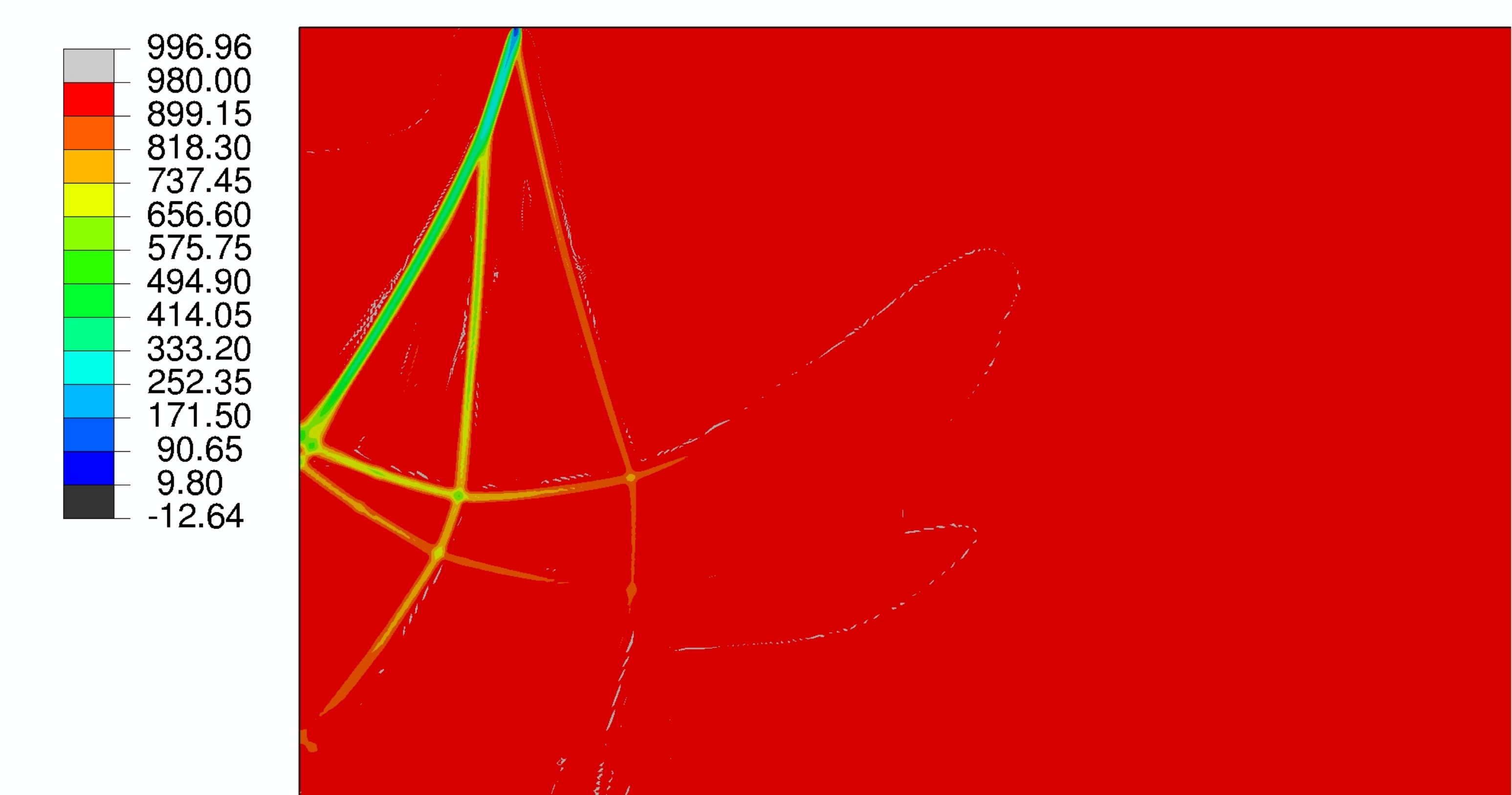} 
		&\includegraphics[trim={2cm 2cm 10cm 0.5cm},clip,width=0.4\textwidth]{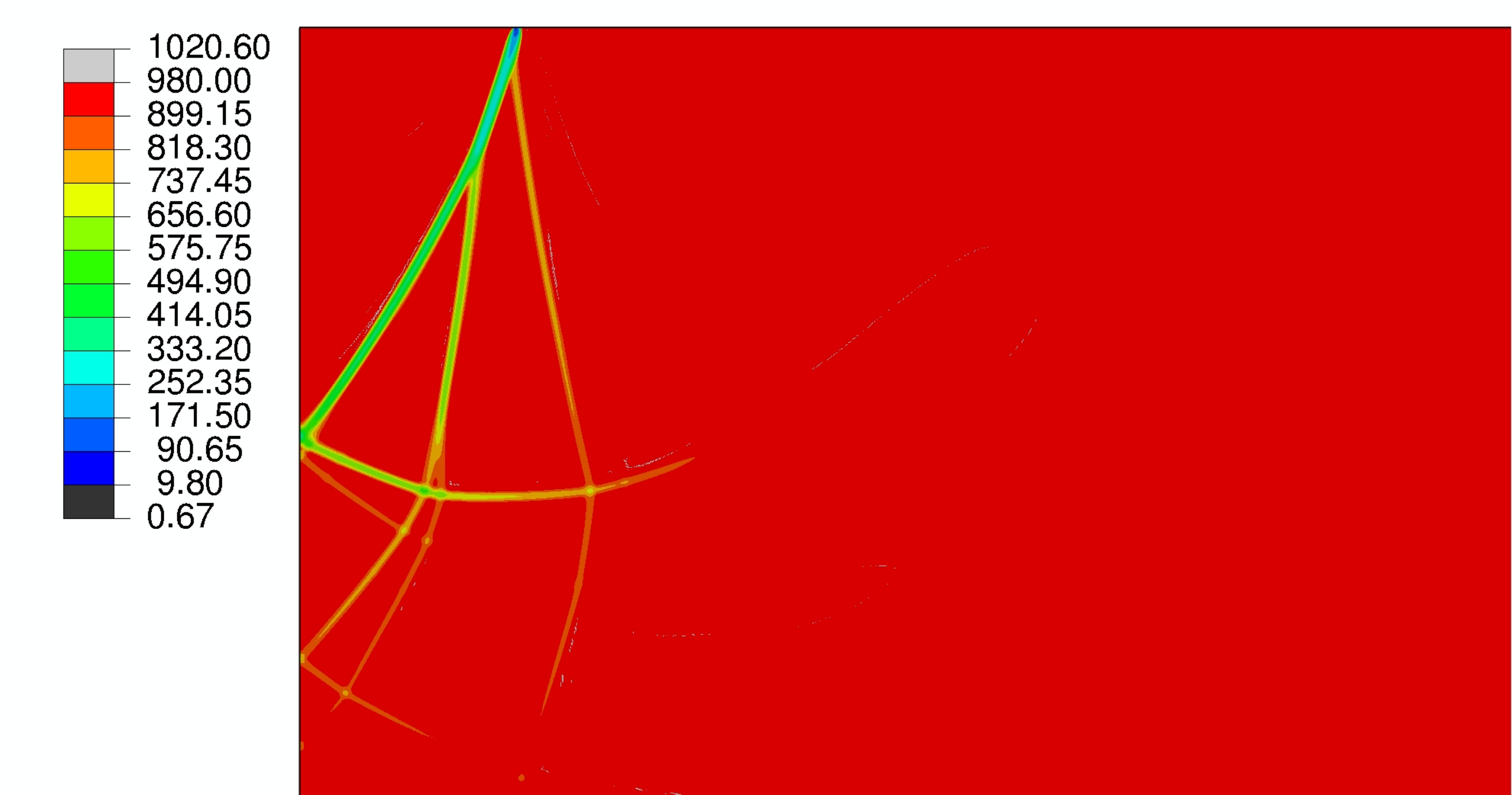} \\
		(a) & (b)\\
		\includegraphics[trim={2cm 2cm 10cm 0.5cm},clip,width=0.4\textwidth]{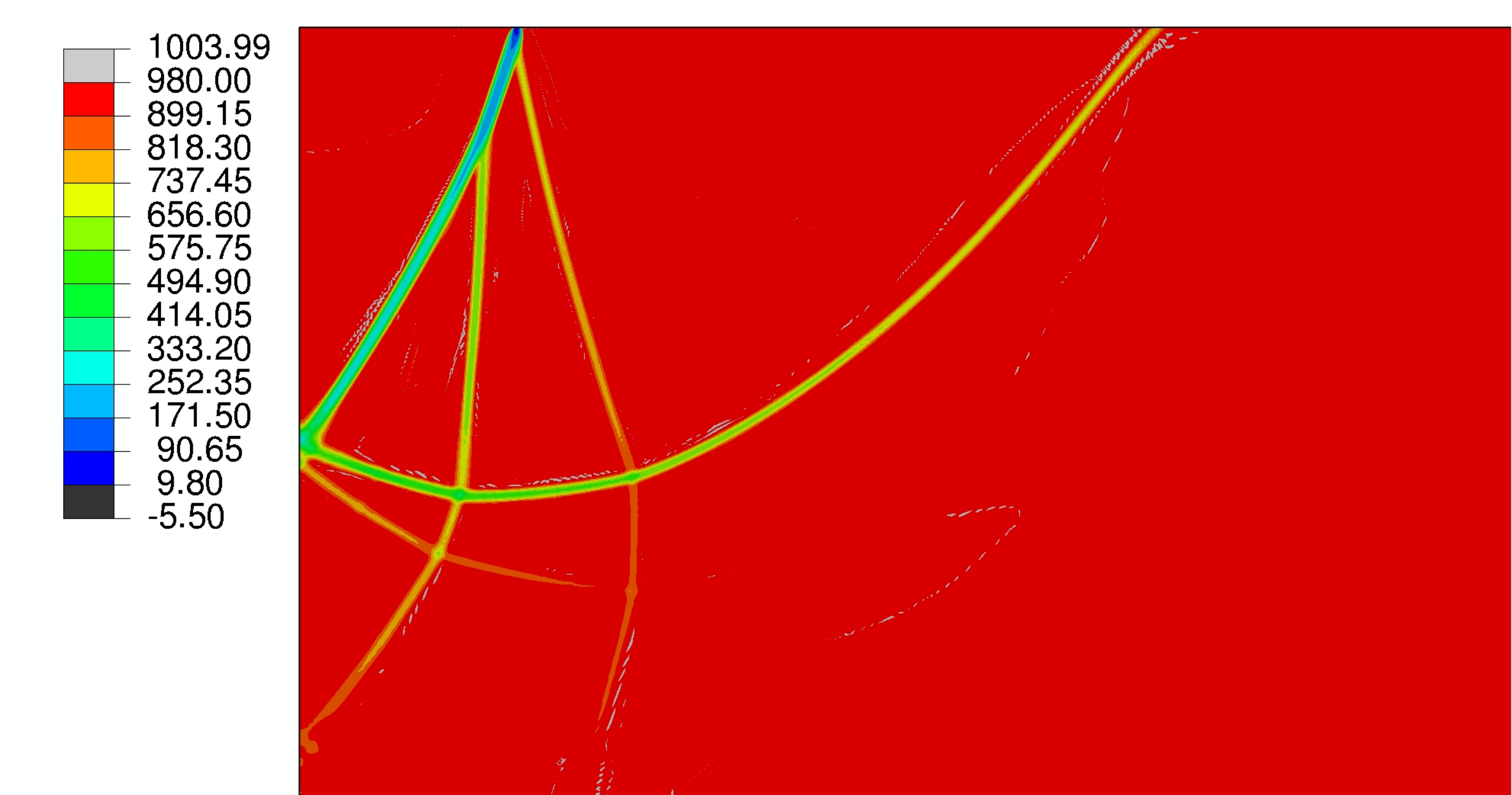} 
		&\includegraphics[trim={2cm 2cm 10cm 0.5cm},clip,width=0.4\textwidth]{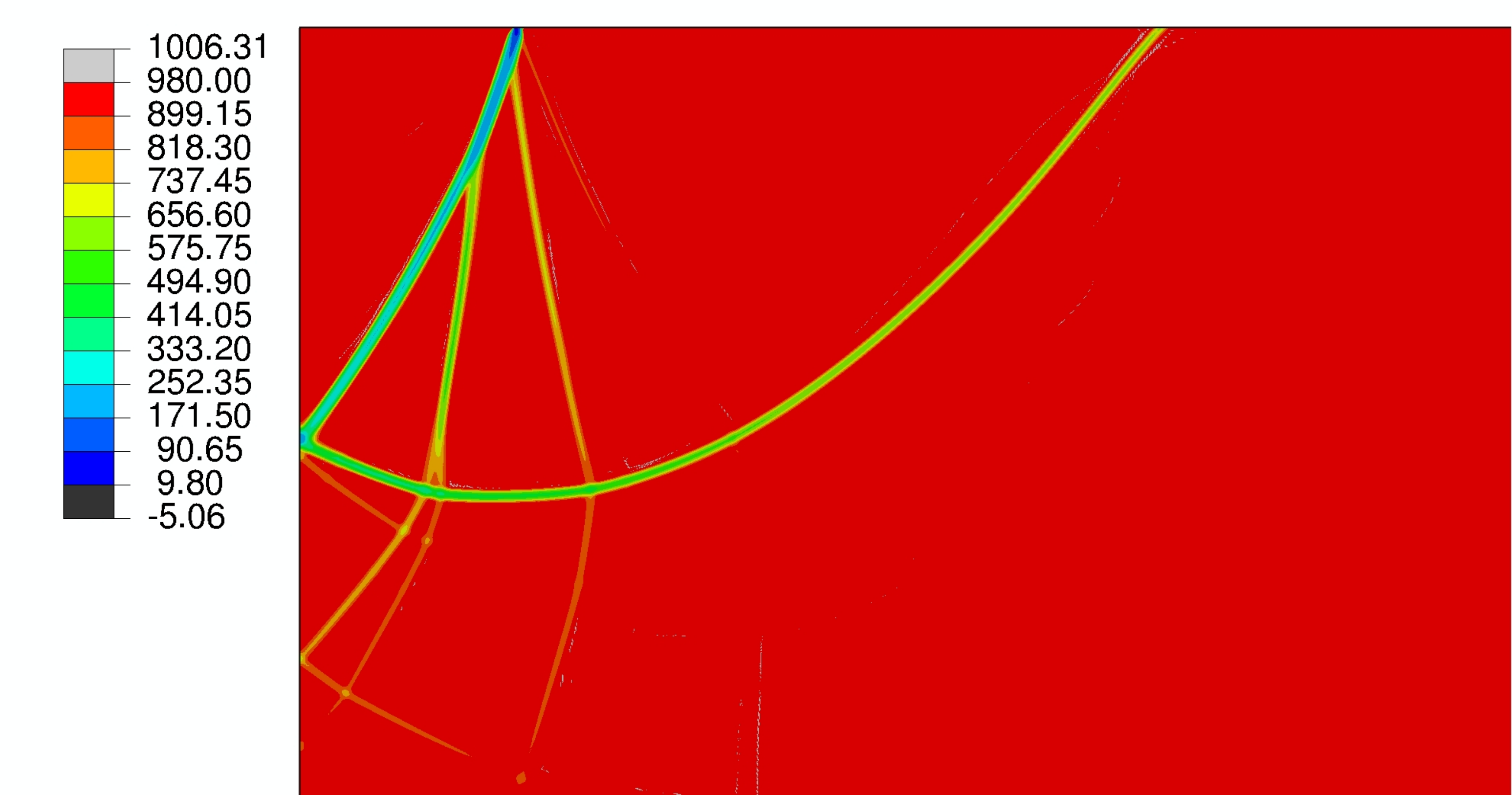} \\
		(c) & (d)\\
		\includegraphics[trim={2cm 2cm 10cm 0.5cm},clip,width=0.4\textwidth]{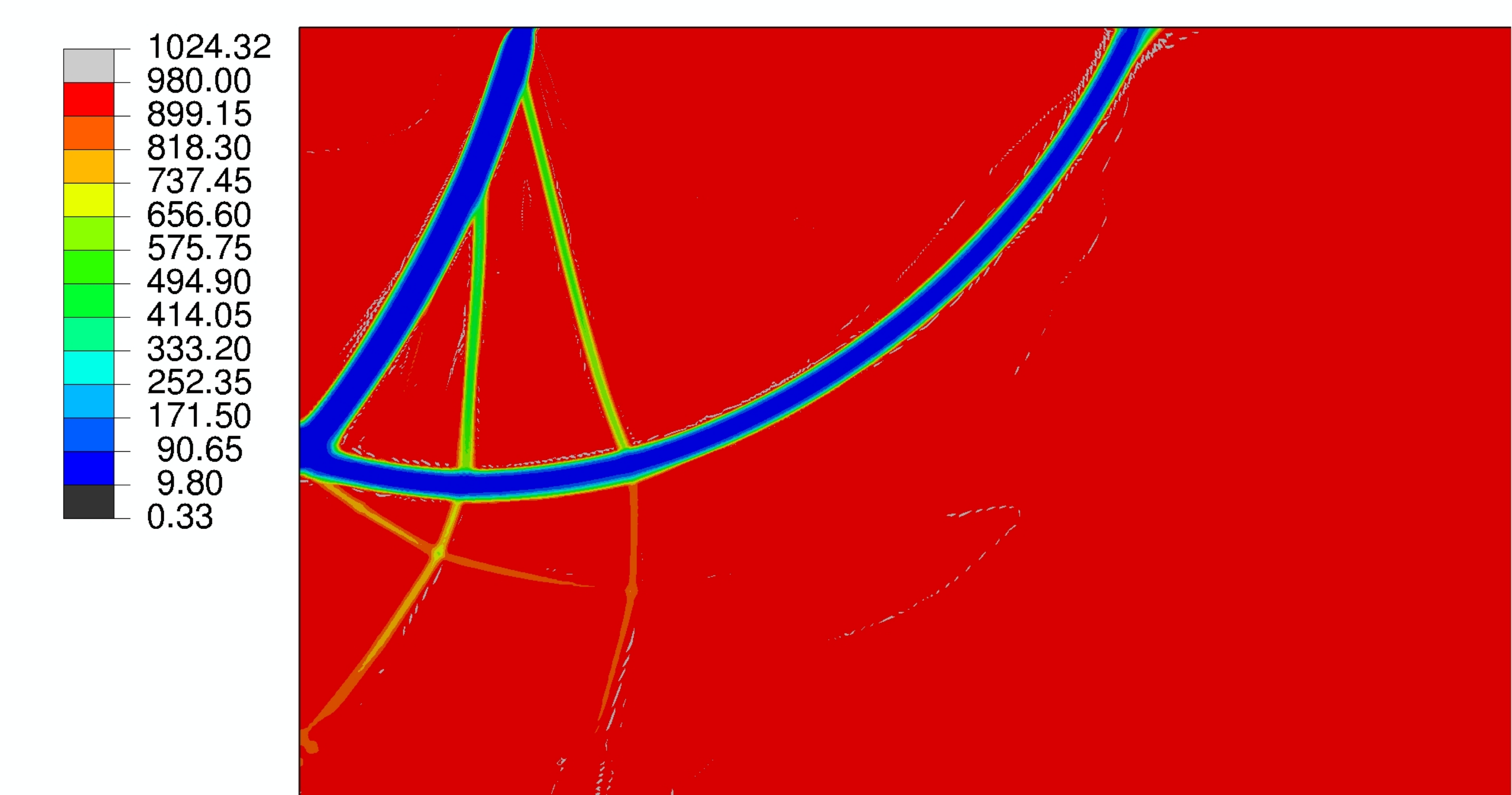} 
		&\includegraphics[trim={2cm 2cm 10cm 0.5cm},clip,width=0.4\textwidth]{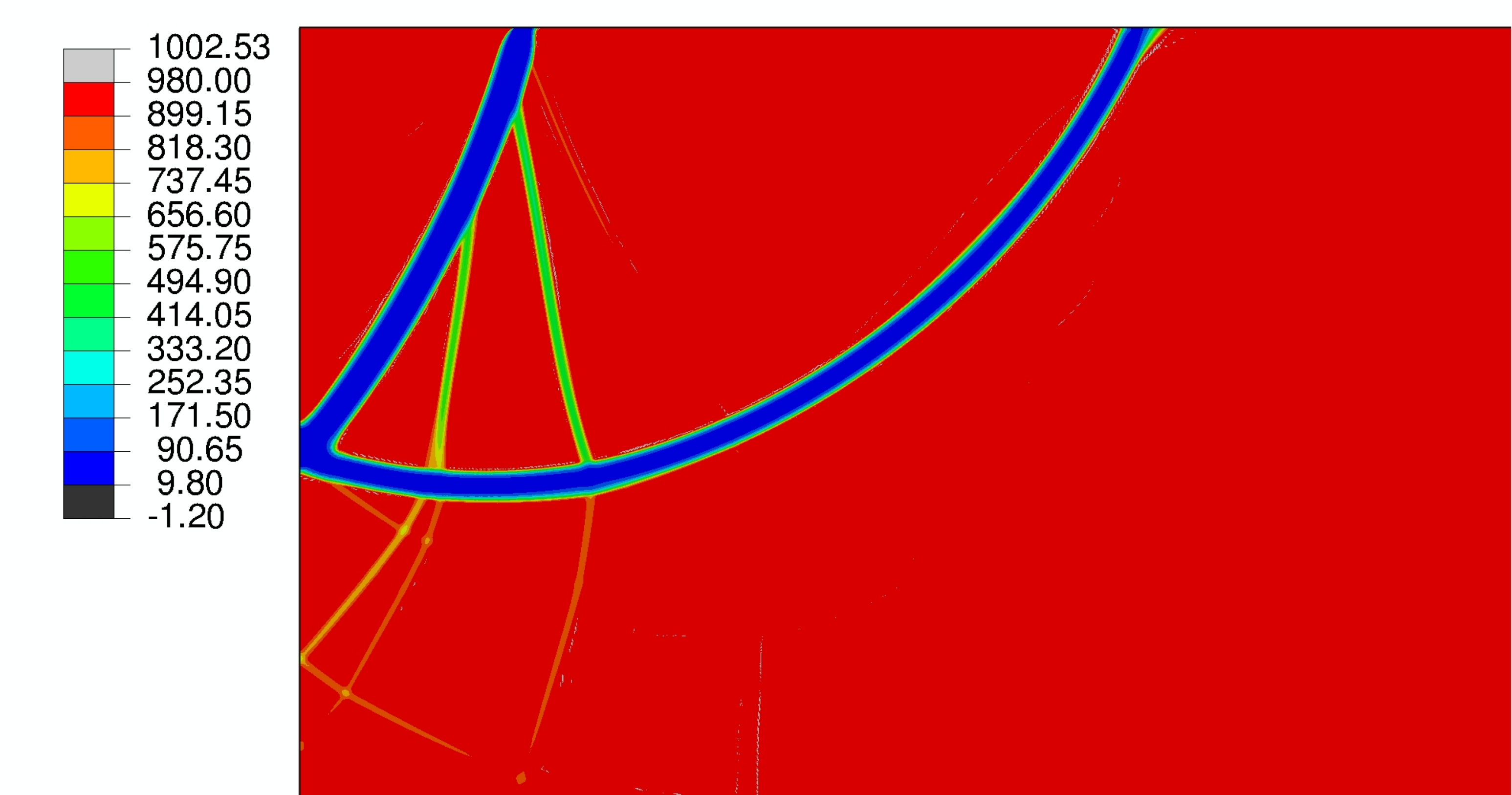} \\
		(e) & (f)
        \end{tabular}
\caption{Spatial distribution of $\kappa$ for the strip footing on Tresca soil with exponential softening, $r=6\cdot10^{-4}$. (a)--(b) Mesh~3 and Mesh~4 at peak load; (c)--(d) at $\bar u/B=0.02$; (e)--(f) at the end of the analysis ($\bar u/B=0.1$). The main shear bands are well captured and essentially coincident, whereas minor differences remain in the finer secondary structures, indicating that a finer discretization is required to fully resolve all localization features.}
\label{fig:tresca_conv_r6e-4}
\end{figure}

The localization patterns are analyzed by considering the spatial distribution of $\kappa$, which is a monotonic decreasing quantity and therefore allows a direct identification of the shear bands as regions of minimum $\kappa$. For $r=2\cdot10^{-3}$ (Fig.~\ref{fig:tresca_conv_r2e-3}), the pattern closely resembles the classical Prandtl mechanism, with two dominant shear bands and a slightly diffuse region between them. The results obtained with Mesh~2 and Mesh~3 are essentially coincident, indicating full mesh objectivity. For $r=10^{-3}$ (Fig.~\ref{fig:tresca_conv_r1e-3}), the localization pattern becomes significantly richer, with the development of a network of secondary shear bands. Both primary and secondary bands are almost perfectly coincident across meshes, showing that the formulation captures even fine-scale features in a mesh-independent manner. 
It is important to emphasize that $\ell$ is a \emph{material} parameter, not a \emph{smoothing} parameter. The value of $\ell$ for a given material depends on its material microstructure, which will determine the appropriate converged localization patterns.

For the smallest value $r=6\cdot10^{-4}$ (Fig.~\ref{fig:tresca_conv_r6e-4}), the pattern becomes even more intricate, with a dense network of secondary bands. The main shear bands are well captured and essentially coincident, whereas some differences persist in the finer secondary structures, indicating that an even finer mesh would be required to fully resolve all details. It is worth noting that the convergence of local fields is more demanding than that of global quantities such as load--displacement curves or dissipated energy, which converge faster with mesh refinement.

These results highlight that the internal length scale does not merely enforce mesh objectivity, but also governs the emergence of different localization mechanisms, leading to qualitatively different failure patterns as the ratio $r=\ell/B$ is varied.

\begin{table}[t]
 	\centering
 	\begin{tabular}{c|c c|c c|c c}
        \hline
 		\hline  		
        &	\multicolumn{2}{ c } {Mesh 2}  &
        \multicolumn{2}{ |c|} {Mesh 3} & \multicolumn{2}{ c} {Mesh 4} \\ 
        \cline{2-7}
 		$r$    &   loading &    total &loading & total &loading & total \\ 
 		&    steps     &  iterations   &    steps     &  iterations & steps     &  iterations\\
		\hline 		
 		$6\cdot10^{-4}$ 		& 250		& 1072	& 247	& 1130	& 261	& 1346\\	
 		$10^{-3}$ 	& 222		& 788	& 228	& 940 & 235	& 1070\\ 
		$2\cdot10^{-3}$ 	& 211		& 606	& 211	& 622 & $/$	& $/$	 \\
		\hline 
		\hline
 	\end{tabular} 
 	\caption{Numerical performance for the strip footing on Tresca soil with exponential softening, in terms of number of load increments and total iterations. Comparable values across meshes confirm the robustness of the formulation, with slightly higher computational effort for smaller values of $r=\ell/B$.}    \label{tab:trescaperform}
 \end{table}

The numerical performance of the simulations is summarized in Table~\ref{tab:trescaperform}. The number of load increments and total iterations remains comparable across different meshes, confirming the robustness of the formulation. The most demanding case corresponds to $r=6\cdot10^{-4}$, where the average time increment remains close to the maximum allowable value ($76.7\%$ of the maximum allowable time increment) and the average number of iterations per increment is approximately $5.2$, indicating a stable nonlinear solution process even in the presence of strong softening.

\begin{figure}[t]
	\centering
	\begin{tabular}{cc}
		\includegraphics[width=0.47\textwidth]{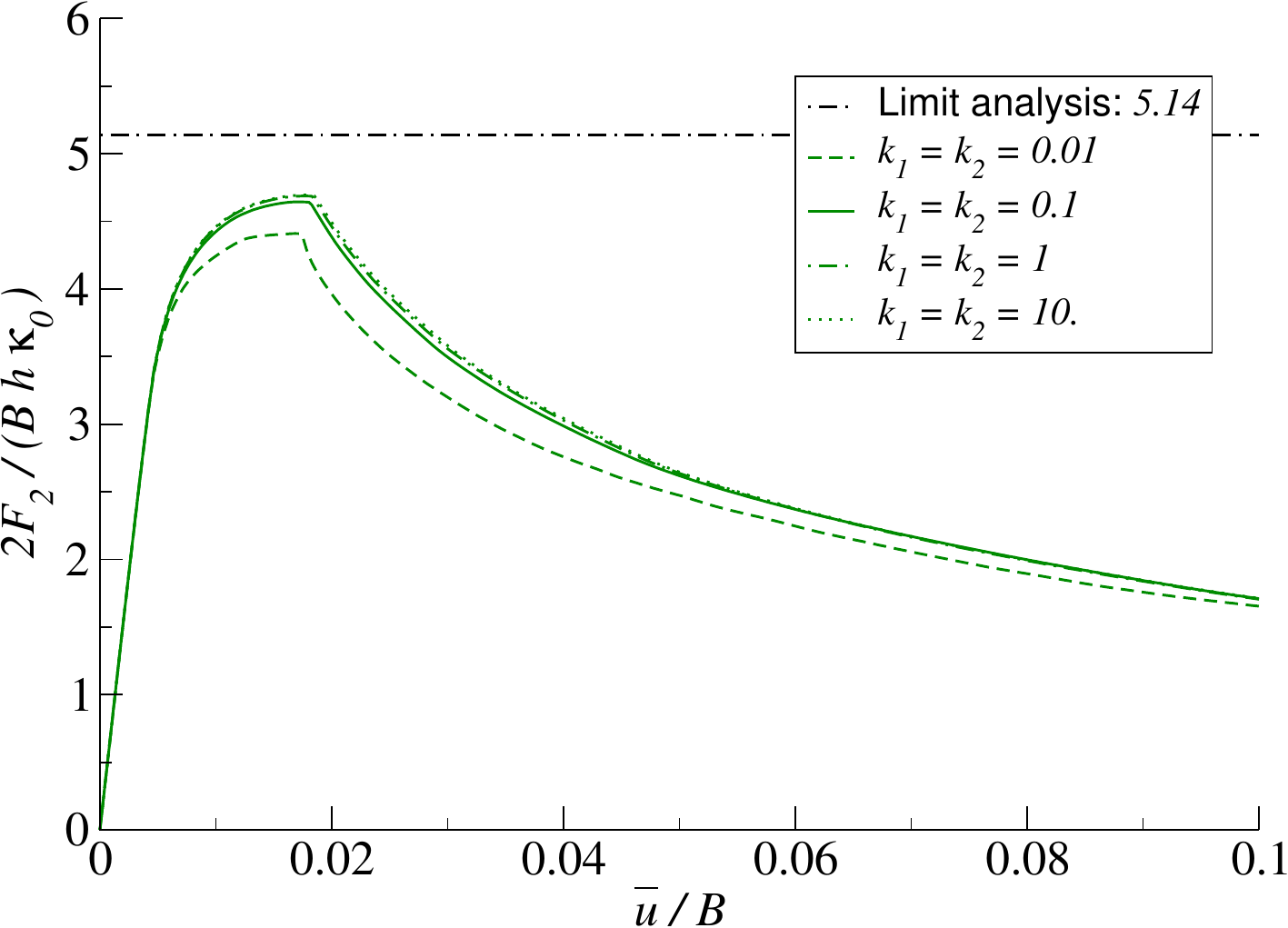} 
		&\includegraphics[width=0.47\textwidth]{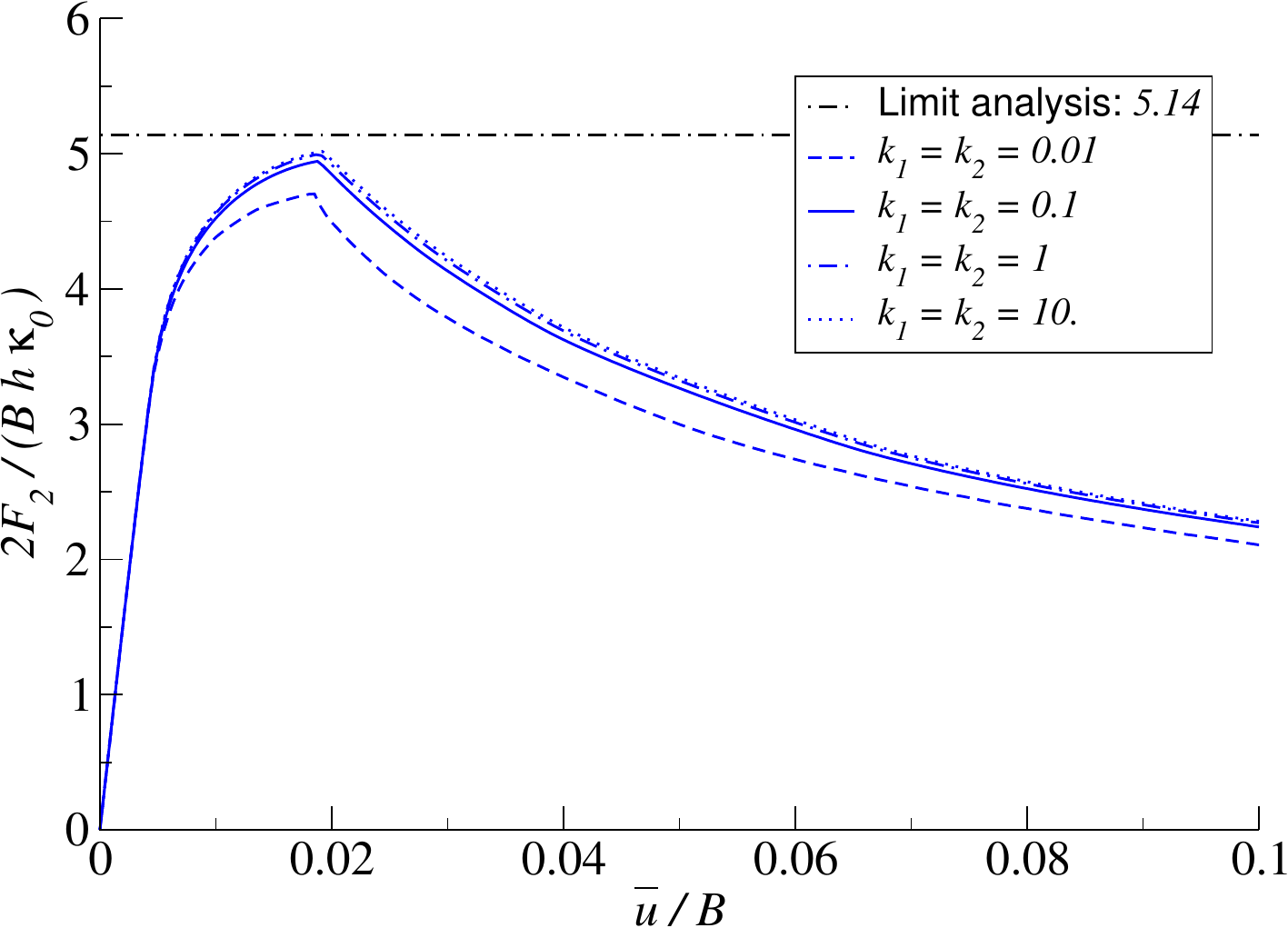}\\
		(a) & (b)
	\end{tabular}
\caption{Influence of the micro-elastic parameters $k_1$ and $k_2$ on the load--displacement response for the strip footing foundation on Tresca soil with exponential softening. Results are shown for Mesh 3, (a) $r=10^{-3}$ and (b) $r=2\cdot10^{-3}$. The response is only moderately affected and tends to saturate for increasing values of $k_1$ and $k_2$, while no significant changes in the localization patterns are observed.}
\label{fig:footing_tresca_inflk}
\end{figure}

Finally, Fig.~\ref{fig:footing_tresca_inflk} shows the influence of the micro-elastic parameters $k_1$ and $k_2$ for $r=10^{-3}$ and $r=2\cdot10^{-3}$. 
Notice that increasing $k_1 = k_2$ tends to increase the peak load but has little effect on the final load. 
The results indicate that the global response is only moderately affected by these parameters and tends to saturate for larger values. No significant changes in the localization patterns are observed, suggesting that the internal length plays a dominant role in governing the structural response.

\FloatBarrier

\subsection{Matsuoka--Nakai soil subject to self-weight}

These analyses are conducted using the Matsuoka--Nakai criterion \cite{matsuoka1974stress, PL2014}, a constitutive model commonly adopted to describe drained soil behavior. The elastic response is assumed to be linear, with shear modulus $G=41.67\,\mathrm{MPa}$ and bulk modulus $K_v=55.56\,\mathrm{MPa}$.

The soil is dry, with unit weight $\gamma_0=18\,\mathrm{kN}/\mathrm{m}^3$. Two sets of material parameters are considered, namely: (i) perfect plasticity, with no cohesion ($\kappa=0$) and constant angle of shearing resistance $\phi=\phi_0=30^\circ$, and (ii) a softening response, associated with a constant cohesion (modeled by assuming $\kappa=8\,\mathrm{kPa}$) and a decreasing angle of shearing resistance.
For the latter case, $\phi$ is assumed to evolve according to a piecewise linear law as a function of the accumulated plastic strain $\bar\varepsilon_p$. In particular, $\phi$ decreases linearly from $\phi_0=20^\circ$ to $\phi_u=10^\circ$ over the interval $0\le \bar\varepsilon_p \le \bar\varepsilon_p^{\,u}$=0.1, and remains constant thereafter:
\begin{equation}
\phi(\bar\varepsilon_p)=
\begin{cases}
\phi_0-\displaystyle\frac{\phi_0-\phi_u}{\bar\varepsilon_p^{\,u}}\;\bar\varepsilon_p,
& 0\le \bar\varepsilon_p \le \bar\varepsilon_p^{\,u},\\[8pt]
\phi_u,
& \bar\varepsilon_p > \bar\varepsilon_p^{\,u}.
\end{cases}
\end{equation}
The corresponding evolution of the material parameter $M$ is obtained from the current value of $\phi$ through Eq.~\eqref{eq:defM}.

Since volumetric plastic deformations may occur, a different definition of the accumulated plastic strain is adopted with respect to the Tresca case, i.e.,
\begin{equation}
    \bar \varepsilon_p=\int \dot \lambda \,\mathrm{d}t,
\end{equation}
where $\dot \lambda$ is the plastic multiplier, as defined in \cite{PL2018}.

A geostatic initial stress state is introduced to account for the soil self-weight. The vertical stress $T_{22}$ is in equilibrium with the body forces due to the soil weight, while the initial horizontal and out-of-plane stresses are assumed as:
\begin{equation}
    T_{11}=T_{33}= T_{22} \left(1-\sin \phi_0\right).
\end{equation}

The load is applied by prescribing the vertical displacement of the footing nodes. No lateral surcharge is applied, and a rough footing--soil interface is assumed. Unless otherwise specified, the micro-continuum parameters are set equal to $k_1=k_2=0.1$.

The perfect plastic case corresponds to the classical $N_\gamma$ problem. This benchmark is known to be particularly demanding from a numerical standpoint: in the absence of cohesion and lateral confinement, the yield condition is satisfied at the ground surface, while the elastic domain develops with depth due to the geostatic stress field. As a result, very steep stress gradients arise in the vicinity of the footing--soil interface.
The analyses of the $N_\gamma$ problem are performed using the standard nonlinear static solver of Abaqus \cite{ABA24}. The initial and maximum time increments are set equal to $10^{-4}\,t_0$ and $5\cdot10^{-3}\,t_0$, respectively, where $t_0$ is a dummy load time. The imposed normalized vertical displacement is $\bar u/B=0.1$.

\begin{figure}[t]
	\centering
	\begin{tabular}{cc}
		\includegraphics[width=0.47\textwidth]{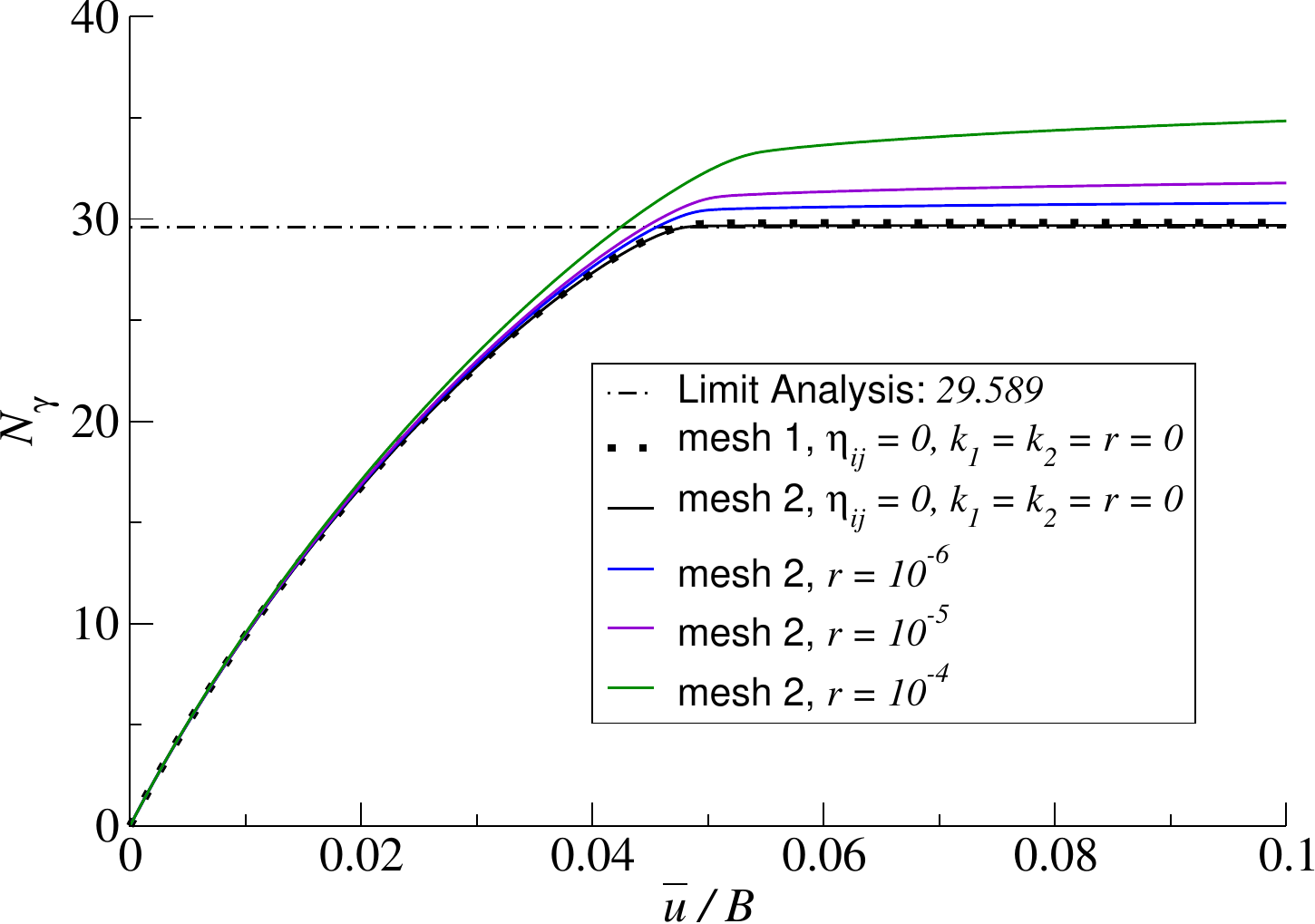} 
		&\includegraphics[width=0.47\textwidth]{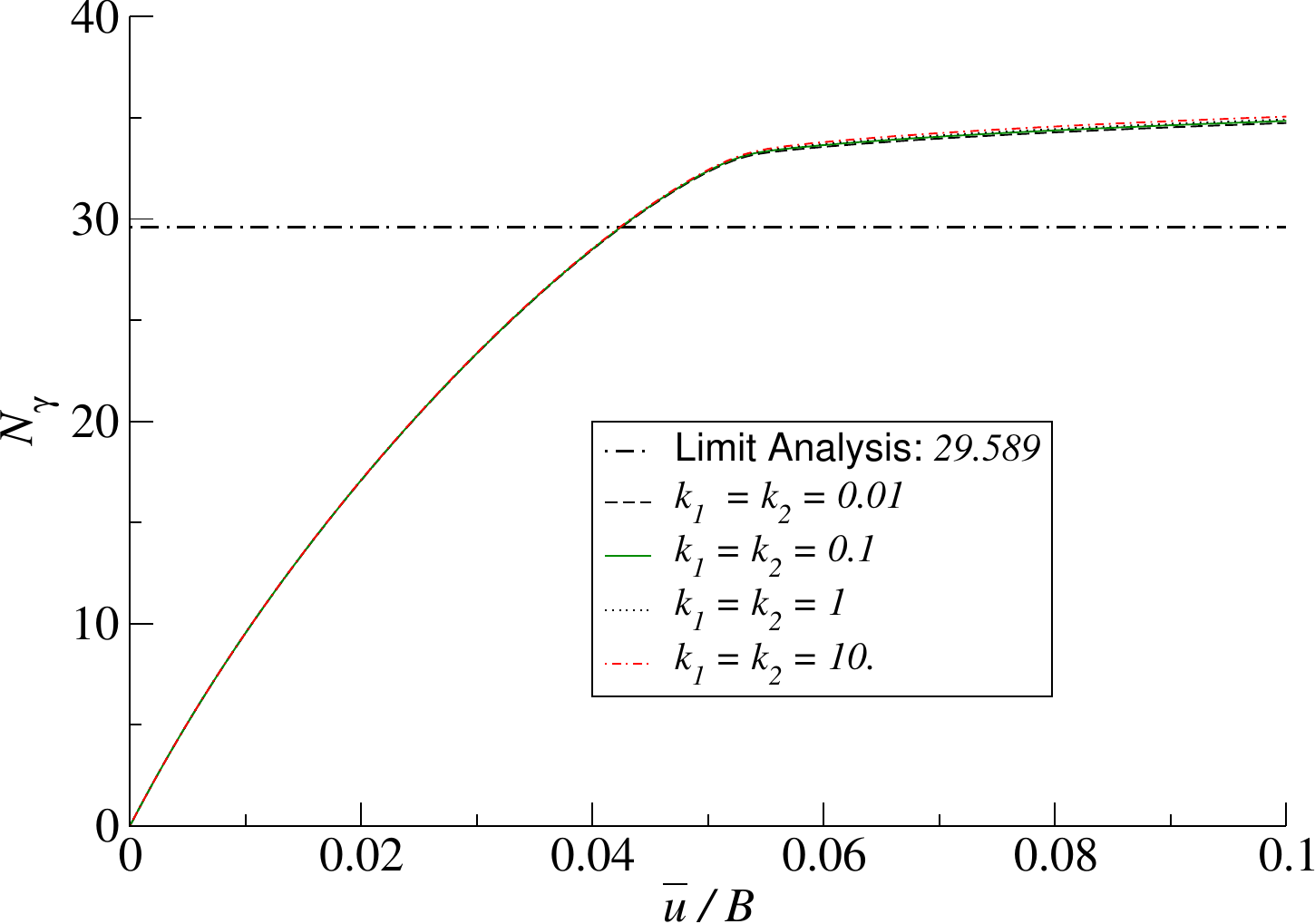}\\
		(a) & (b)
	\end{tabular}
\caption{
Strip footing on Matsuoka--Nakai soil under perfect plasticity ($N_\gamma$ problem).
(a) Normalized load--displacement curves for different values of $r=\ell/B$, showing convergence to the classical Cauchy solution.
(b) Influence of the micro-elastic parameters $k_1$ and $k_2$ for $r=10^{-4}$, Mesh 2: the curves are practically coincident, indicating a negligible effect of these parameters on the global response.
}
\label{fig:footing_ngamma}
\end{figure}

Fig.~\ref{fig:footing_ngamma}a shows the normalized load--displacement curves obtained under perfect plasticity. The dimensionless load is defined as
\begin{equation}
    N_\gamma=\frac{F_2}{B\,\gamma_0},
\end{equation}
where $F_2$ is the resultant vertical reaction per unit depth at the footing.

Only Mesh~1 and Mesh~2 are considered, since no mesh dependence is expected in the absence of softening. The classical Cauchy solution is recovered by enforcing $\boldeta=\bzero$ at all nodes and setting $k_1=k_2=\ell=0$. The reference solution is obtained from limit analysis using the method of characteristics, computed with the $ABC$ code developed by Martin \cite{Martin2003,Martin2004} and extended to general yield criteria following \cite{LP2017}.

The deformable Cosserat results are reported for $r=\ell/B=10^{-6}$, $10^{-5}$, and $10^{-4}$. The results show that the proposed formulation recovers the classical response in the limit $r \to 0$. Increasing $r$ leads to a slightly stiffer response and to a modest increase in the peak load, indicating the emergence of a size effect when the internal length becomes non-negligible compared to the footing width.

Fig.~\ref{fig:footing_ngamma}b shows the influence of the micro-elastic parameters $k_1$ and $k_2$ for $r=10^{-4}$. The curves are essentially coincident, indicating that the global response is weakly affected by these parameters under perfect plasticity.

\begin{figure}[t]
	\centering
	\begin{tabular}{cc}
		\includegraphics[trim={2cm 8cm 2cm 0.5cm},clip,width=0.474\textwidth]{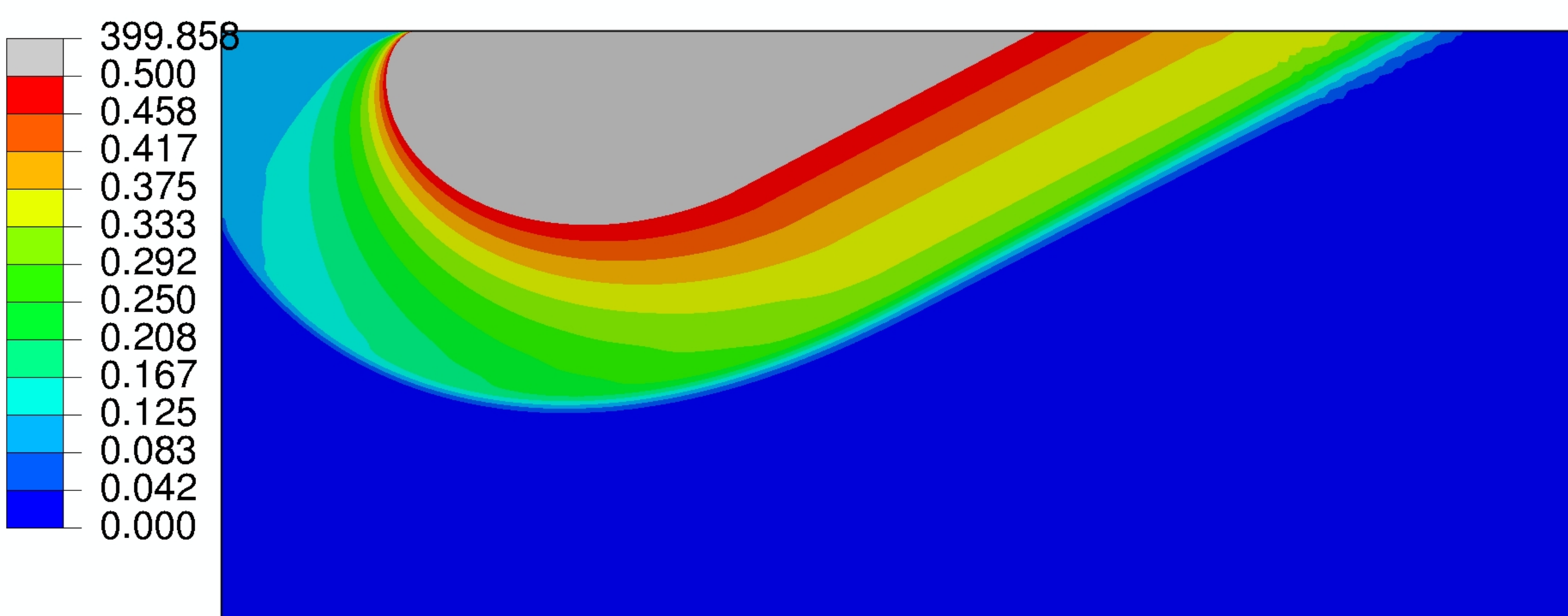} 
		&\includegraphics[trim={2cm 8cm 2cm 0.5cm},clip,width=0.474\textwidth]{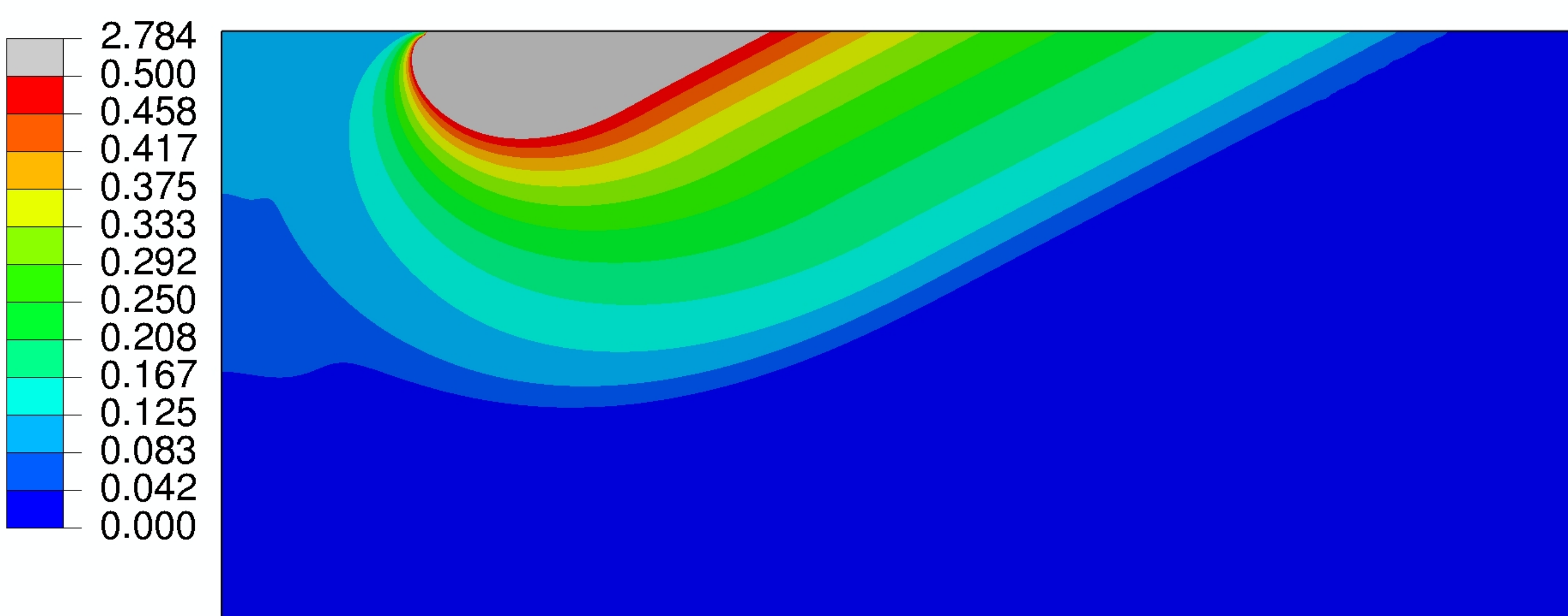} \\
		(a) & (b)\\
		\includegraphics[trim={2cm 8cm 2cm 0.5cm},clip,width=0.474\textwidth]{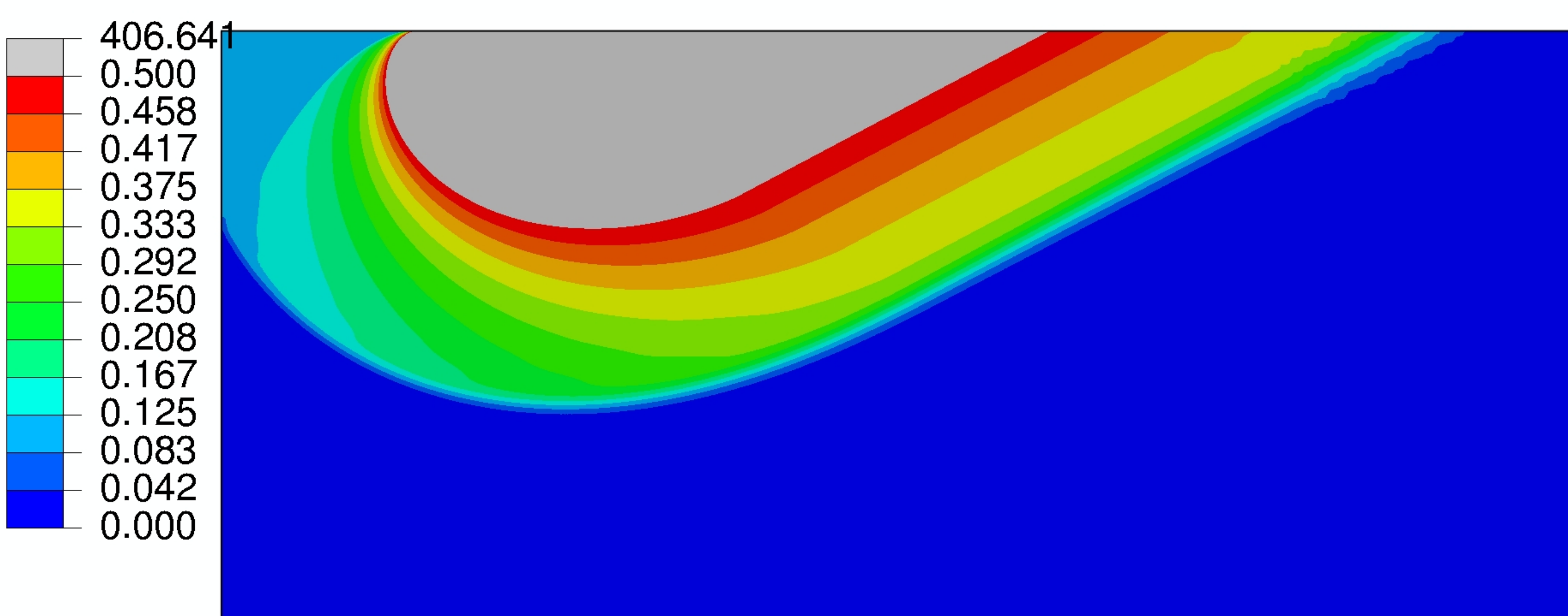} 
		&\includegraphics[trim={2cm 8cm 2cm 0.5cm},clip,width=0.474\textwidth]{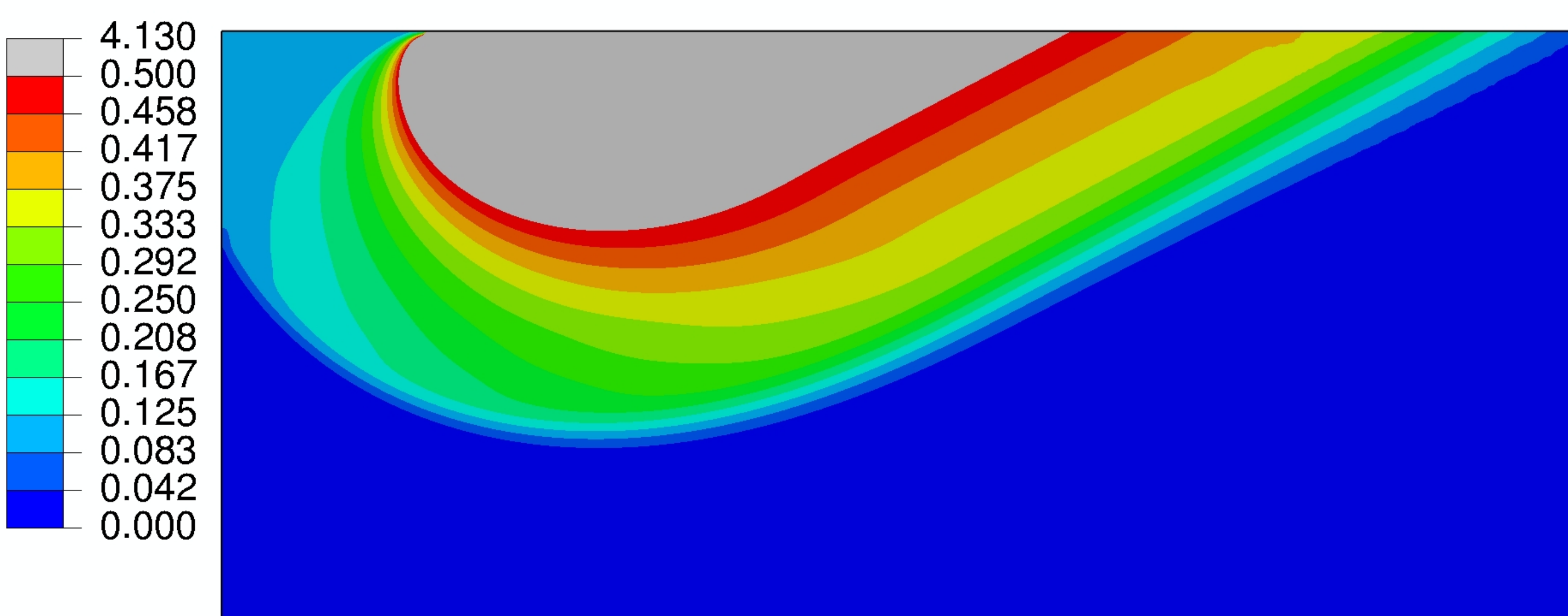} \\
		(c) & (d)\\
		\includegraphics[trim={2cm 8cm 2cm 0.5cm},clip,width=0.474\textwidth]{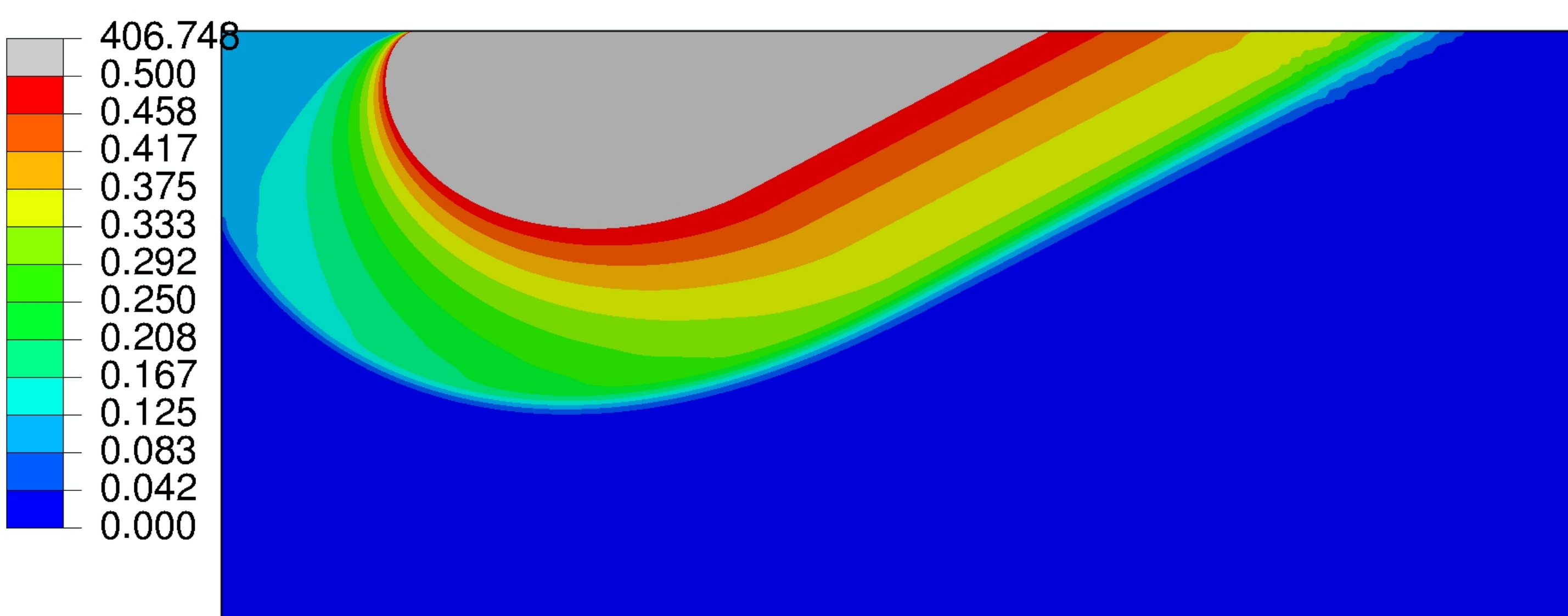} 
		&\includegraphics[trim={2cm 8cm 2cm 0.5cm},clip,width=0.474\textwidth]{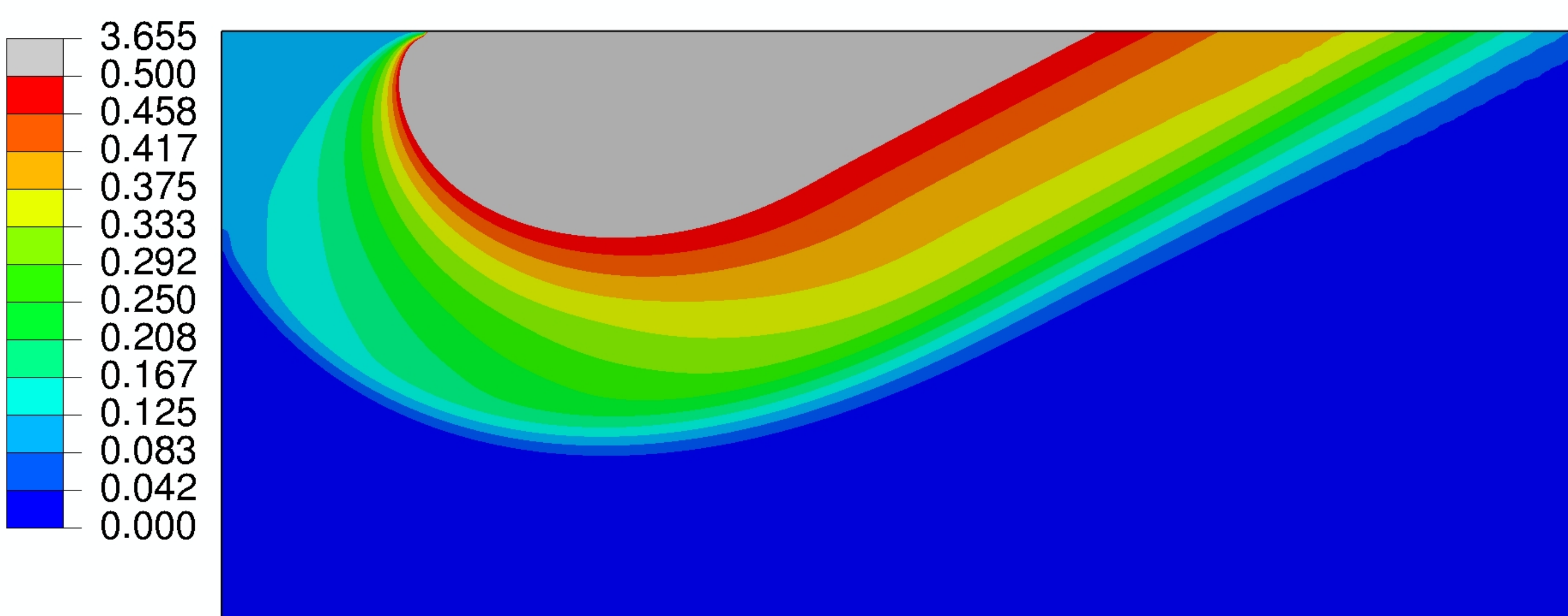} \\
		(e) & (f)
        \end{tabular}
\caption{Comparison between $|\dot{\bu}|$ obtained with Mesh 2 for the Cauchy Continuum (a-c-e) and for $r=10^{-4}$ (b-d-f) at different loading levels for the $N_\gamma$ problem. 
(a) Cauchy, $\bar u/B=0.05$, (b) $r=10^{-4}$, $\bar u/B=0.05$. (c) Cauchy, $\bar u/B=0.07$, (d) $r=10^{-4}$, $\bar u/B=0.07$.
(e) Cauchy, $\bar u/B=0.1$ (end of the analysis), (f) $r=10^{-4}$, $\bar u/B=0.1$ (end of the analysis)}
\label{fig:footing_velocity_field2}
\end{figure}

Fig.~\ref{fig:footing_velocity_field2} shows the spatial distribution of the velocity magnitude at different loading stages, for both the Cauchy continuum and the Cosserat model with $r=10^{-4}$.  Since the response is rate-independent the magnitude of the velocity
is a relative quantity. 
In the Cosserat case, the magnitude of the velocity is approximately two orders of magnitude smaller than in the classical solution, indicating a more distributed deformation process and a reduced tendency toward instantaneous localization.

\begin{figure}[t]
	\centering
	\begin{tabular}{cc}
		\includegraphics[width=0.47\textwidth]{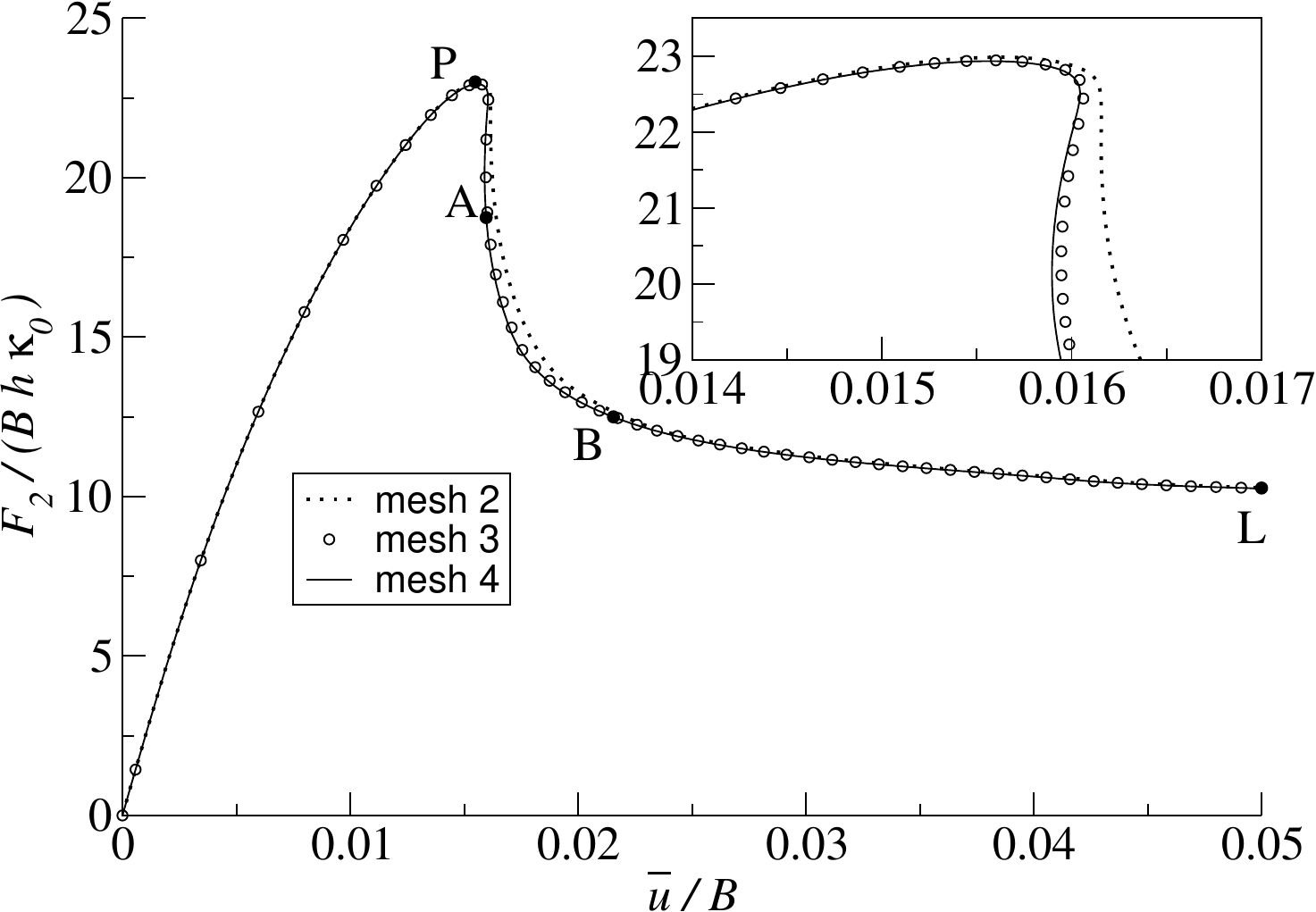} 
		&\includegraphics[width=0.47\textwidth]{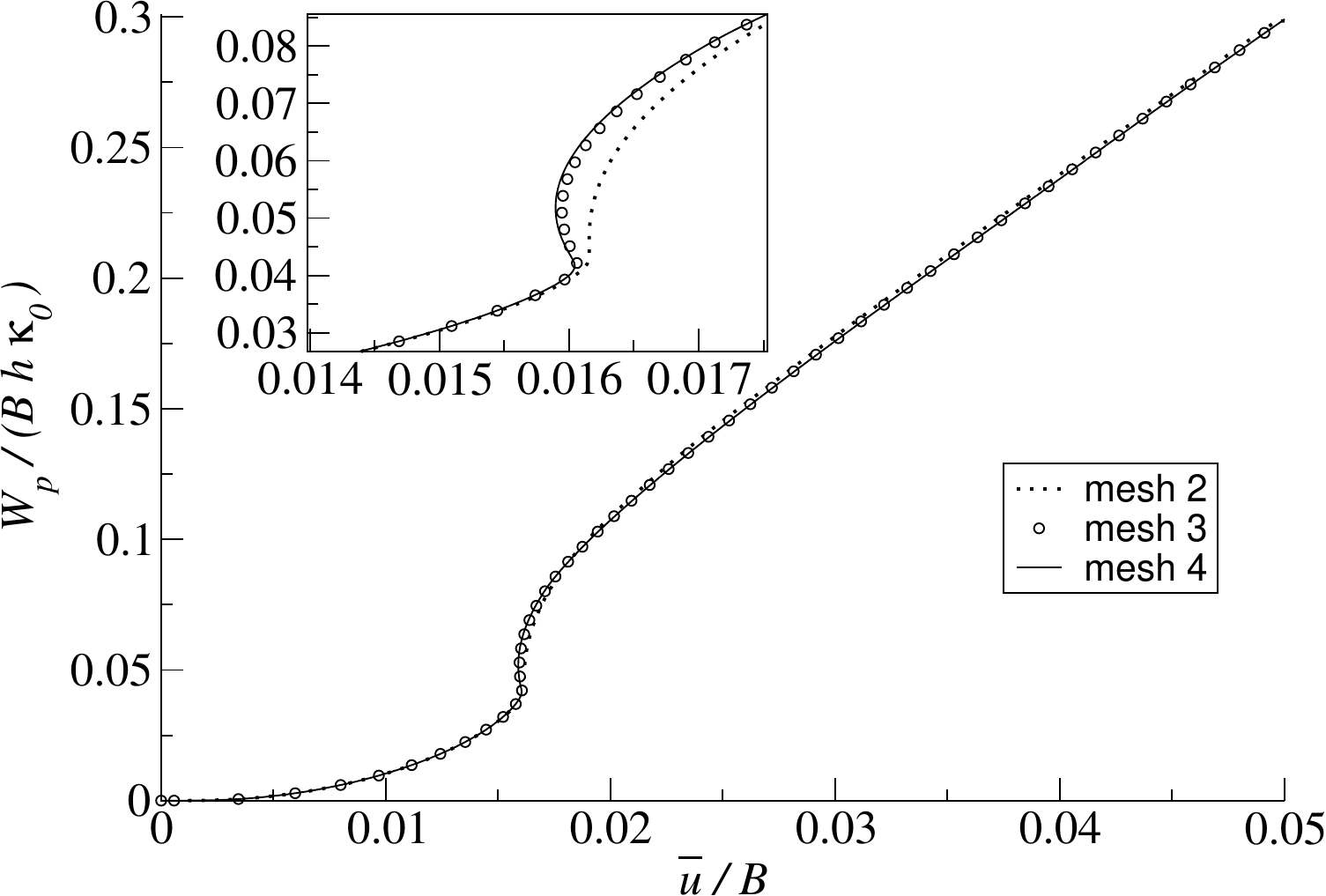}\\
		(a) & (b)
	\end{tabular}
\caption{
Convergence of the solution for the strip footing on Matsuoka--Nakai soil with softening, obtained using the Riks procedure for $r=4.5\cdot10^{-4}$.
(a) Load--displacement curves; (b) dissipated energy as a function of the normalized footing displacement. Both quantities exhibit a non-monotonic response with a pronounced post-peak regime associated with the unstable behavior, while showing clear convergence upon mesh refinement.
}
\label{fig:MN_soft_conv}
\end{figure}

The softening response is characterized by a pronounced loss of strength associated with the reduction of the angle of shearing resistance (see Fig.~\ref{fig:MN_soft_conv}). As a result, the global load--displacement response exhibits a non-monotonic behavior with a distinct post-peak regime, including an \emph{unstable} region. For this reason, the analyses cannot be performed using a standard displacement-controlled procedure and are instead carried out using an arc-length method \cite{deSouzaNeto}. In particular, the Riks algorithm available in Abaqus \cite{ABA24} is adopted. 
The initial arc-length increment is set equal to $1\cdot10^{-4}$, with a maximum increment equal to $5\cdot 10^{-3}$. The analyses have been stopped at the first increment in which $\bar u/B\geq0.05$.

The convergence properties of the proposed formulation are first assessed by considering the smallest value of the normalized internal length adopted in this study, namely $r=\ell/B=4.5\cdot10^{-4}$. This choice corresponds to the most demanding case in terms of localization, as smaller values of $r$ lead to increasingly complicated deformation patterns and sharper strain gradients.

It should be noted that the values of $r$ adopted in the present softening analyses are larger than those used for the $N_\gamma$ problem under perfect plasticity. This is intentional: in the absence of softening, very small values of $r$ can be employed to recover the classical Cauchy solution, whereas in the presence of strong strain softening, the choice of $r$ must ensure an adequate resolution of the localization mechanisms within a feasible mesh discretization.

Fig.~\ref{fig:MN_soft_conv} shows the convergence of the solution in terms of load--displacement curves (a) and dissipated energy as a function of the strip footing displacement (b). Both these curves exhibit a non-monotonic response with a pronounced post-peak regime due to the unstable behavior. 
Despite the presence of an unstable behavior, the load-displacement curves obtained with Mesh~3 and Mesh~4  are essentially coincident, indicating convergence of the global response. The same conclusion can be drawn from the dissipated energy, reported as a function of the  strip footing displacement
in Fig.~\ref{fig:MN_soft_conv}b. The entire dissipation history converges with mesh refinement, further confirming the consistency of the formulation.

\begin{figure}[t]
	\centering
	\begin{tabular}{cc}
		\includegraphics[trim={2cm 8cm 2cm 0.5cm},clip,width=0.475\textwidth]{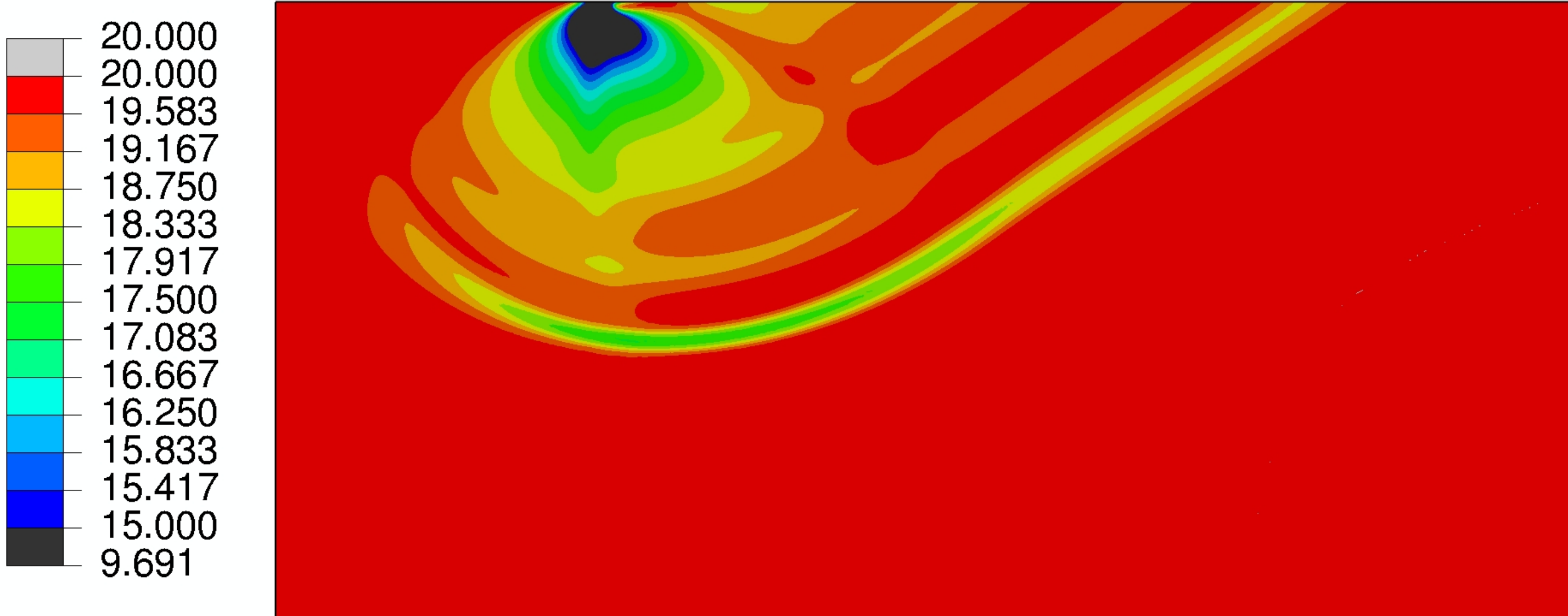} 
		&\includegraphics[trim={2cm 8cm 2cm 0.5cm},clip,width=0.475\textwidth]{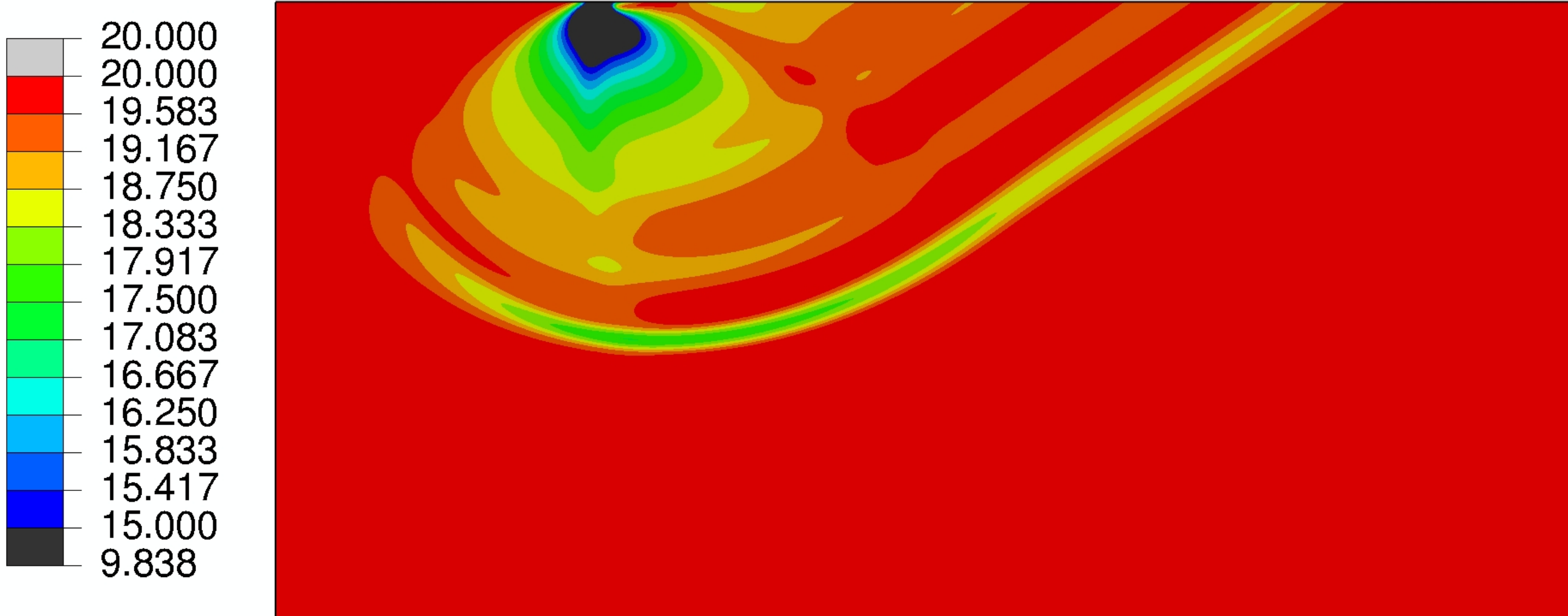} \\
		(a) & (b)\\
		\includegraphics[trim={2cm 8cm 2cm 0.5cm},clip,width=0.475\textwidth]{figures/riks_ell4,5e-4_mesh3_F150KPa_phi.pdf} 
		&\includegraphics[trim={2cm 8cm 2cm 0.5cm},clip,width=0.475\textwidth]{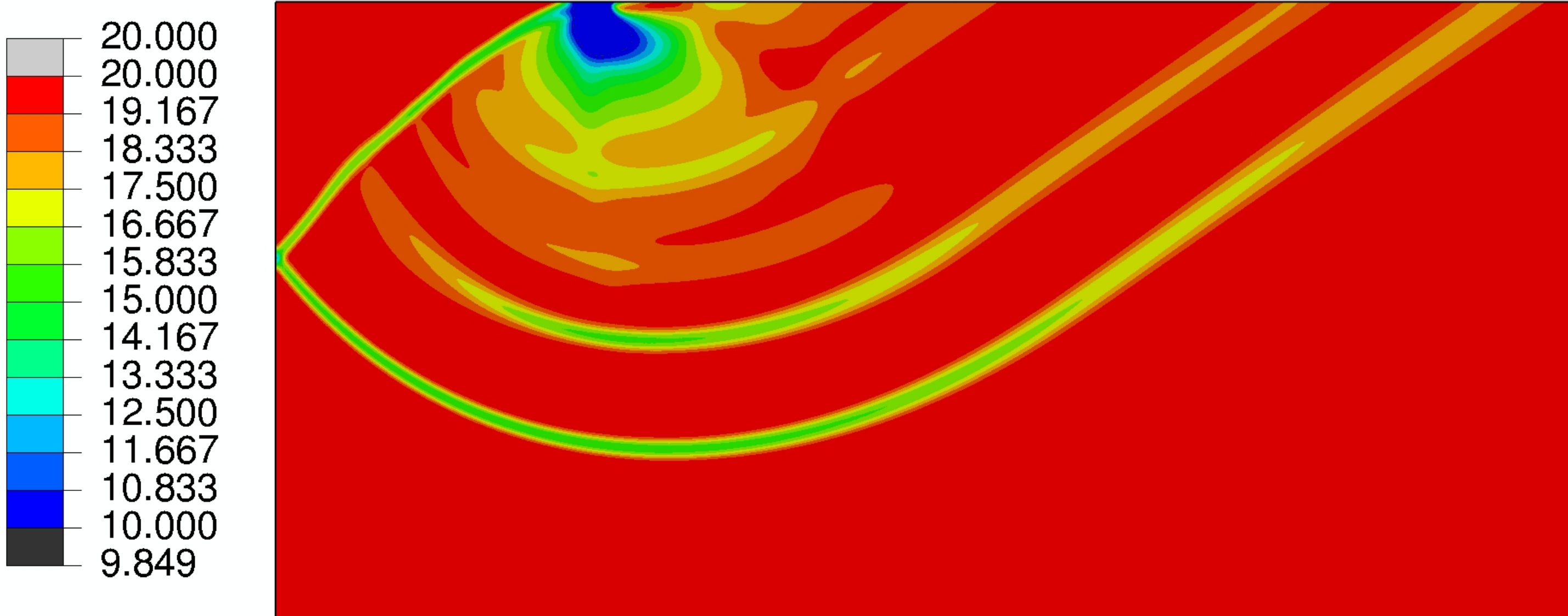} \\
		(c) & (d)  \\
        \includegraphics[trim={2cm 8cm 2cm 0.5cm},clip,width=0.475\textwidth]{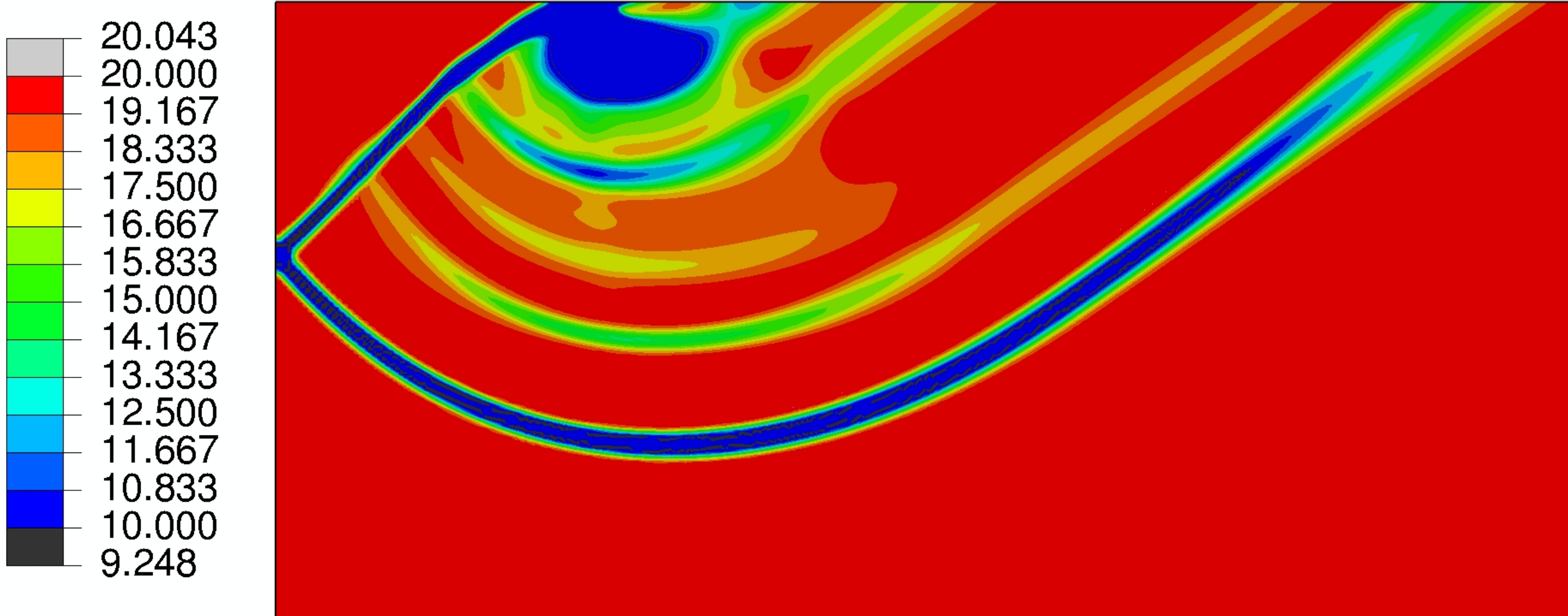} 
		&\includegraphics[trim={2cm 8cm 2cm 0.5cm},clip,width=0.475\textwidth]{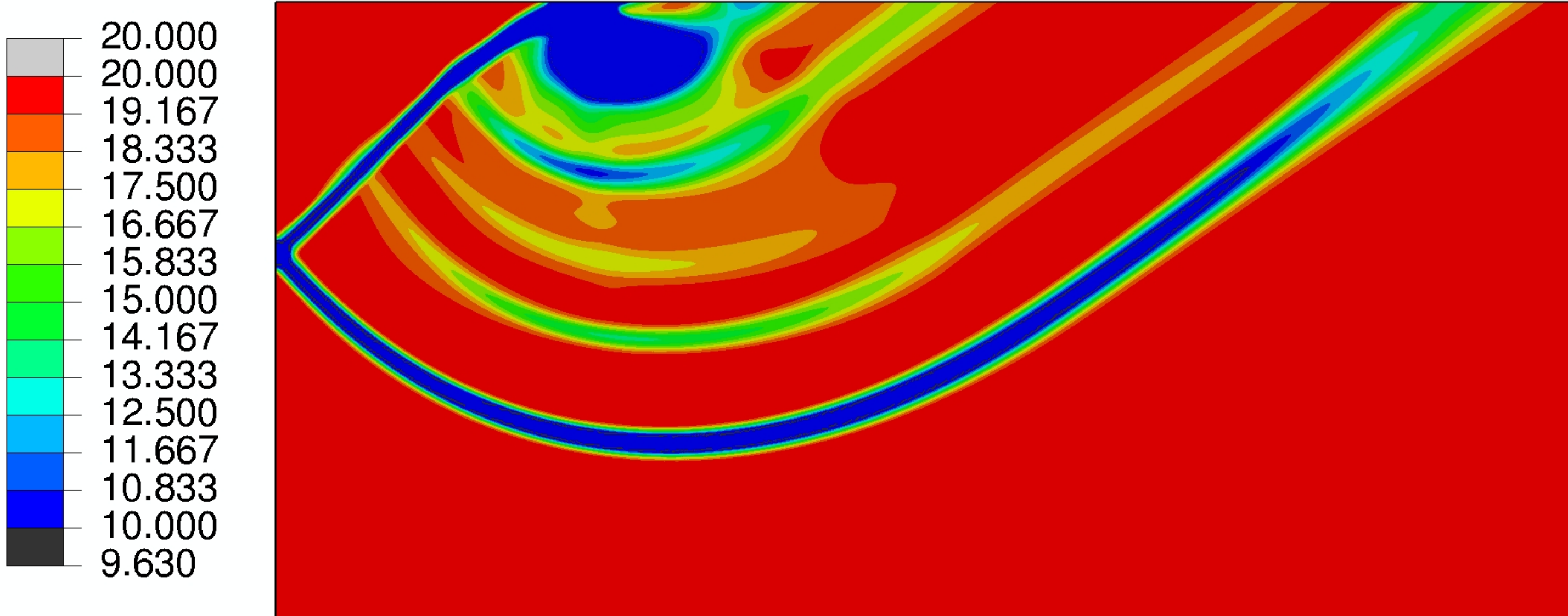} \\
		(e) & (f) \\
        \includegraphics[trim={2cm 8cm 2cm 0.5cm},clip,width=0.475\textwidth]{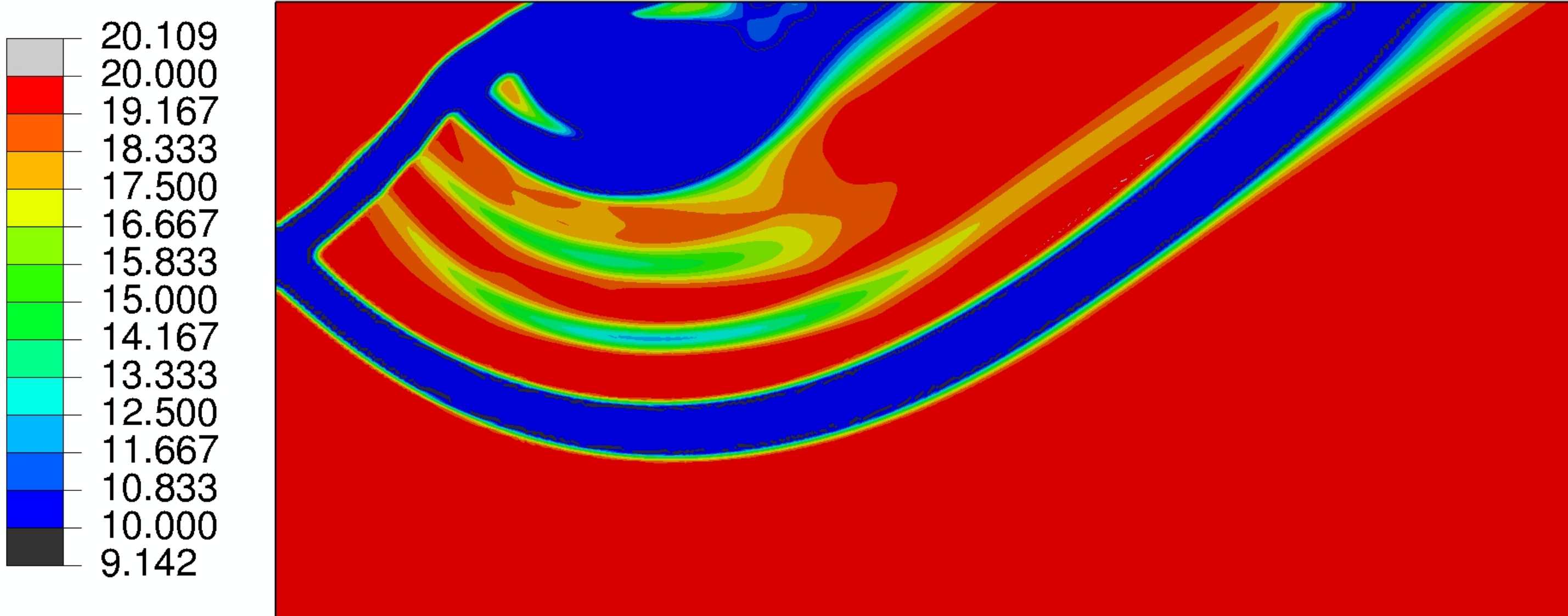} 
		&\includegraphics[trim={2cm 8cm 2cm 0.5cm},clip,width=0.475\textwidth]{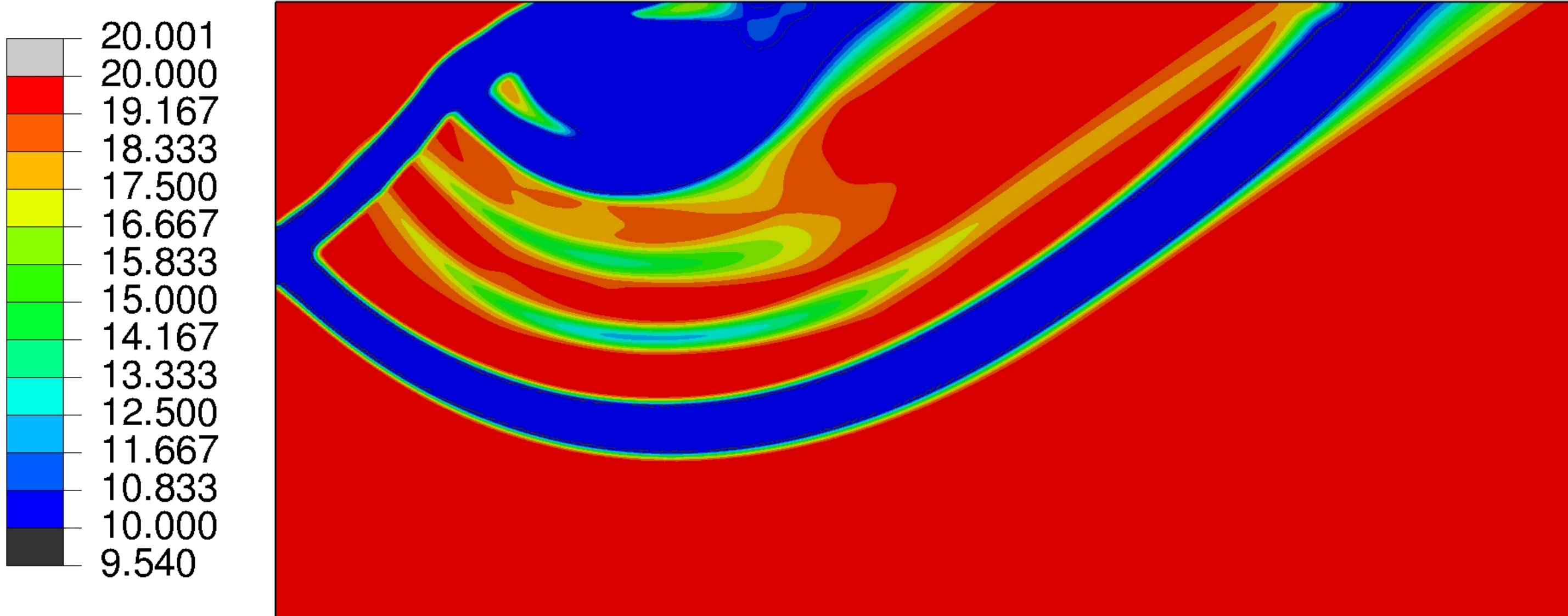} \\
		(g) & (h) 
	\end{tabular}
\caption{Convergence of the $\phi$ profiles for $r=4.5 \cdot 10^{-4}$.
(a) Mesh 3, stress peak (point $P$ of Fig. \ref{fig:MN_soft_conv}a), (b)  Mesh 4, stress peak (point $P$ of Fig. \ref{fig:MN_soft_conv}a). (c) Mesh 3, point $A$ of Fig. \ref{fig:MN_soft_conv}a, (d) Mesh 4,  point $A$ of Fig. \ref{fig:MN_soft_conv}a. (d) Mesh 3, point $B$ of Fig. \ref{fig:MN_soft_conv}a, (e) Mesh 4,  point $B$ of Fig. \ref{fig:MN_soft_conv}a. (c) Mesh 3, point $L$ of Fig. \ref{fig:MN_soft_conv}a, (d) Mesh 4,  point $L$ of Fig. \ref{fig:MN_soft_conv}a. 
}
\label{fig:footing_softening_convergence_shearbands}
\end{figure}

The spatial distribution of $\phi$, reported in Fig. \ref{fig:footing_softening_convergence_shearbands}, shows that the localization patterns are also fully converged between Mesh~3 and Mesh~4. In particular, both the primary and secondary shear bands are coincident during all their evolution, indicating that the proposed formulation is able to capture complex localization mechanisms in a mesh-independent manner even in the presence of strong softening and unstable structural response.

It is worth emphasizing that, compared to the Tresca case discussed in the previous section, the convergence behavior is improved, despite the more severe nonlinear response. This further highlights the robustness of the proposed approach.

\begin{table}[t]
 	\centering
 	\begin{tabular}{c|c c c}
        \hline
 		\hline  		
 		Mesh    &   loading &    total & total arc  \\ 
 		&           steps     &  iterations &length  \\
		\hline 		
 		$2$ 		& 598		& 1739 &	2.914 \\	
 		$3$ 	& 592		& 1750 & 2.916	\\ 
		$4$ 	& 592		& 1746 &2.916	 \\
		\hline 
		\hline
 	\end{tabular} 
   \caption{Numerical performance for the strip footing on Matsuoka--Nakai soil with linear softening of the angle of shearing resistance, in terms of number of load increments and total nonlinear iterations, for $r=4.5\cdot10^{-4}$ using the Riks procedure. The results are very similar across meshes, indicating mesh-independent computational cost and a robust solution procedure despite the presence of an unstable behavior.}
	    \label{tab:MNr4.5e-4}
 \end{table}

The computational cost of the analyses is summarized in Tab.~\ref{tab:MNr4.5e-4}. The average number of iterations per increment is $2.91$ for Mesh~2 and $2.95$ for both Mesh~3 and Mesh~4. 
The average arc-length increment is $4.87\cdot10^{-3}$ for Mesh~2 and $4.93\cdot10^{-3}$ for Mesh~3 and Mesh~4. It is worth noting that the average number of iterations per increment remains low, and that the arc-length increment is very close to the maximum allowable value, set equal to $5\cdot10^{-3}$.
These results indicate a very stable nonlinear solution process, even in the presence of a pronounced softening response related to an unstable behavior.
Moreover, the computational cost is essentially independent of the mesh refinement, further confirming the robustness of the proposed formulation.

\begin{figure}[t]
	\centering
		\includegraphics[width=0.47\textwidth]{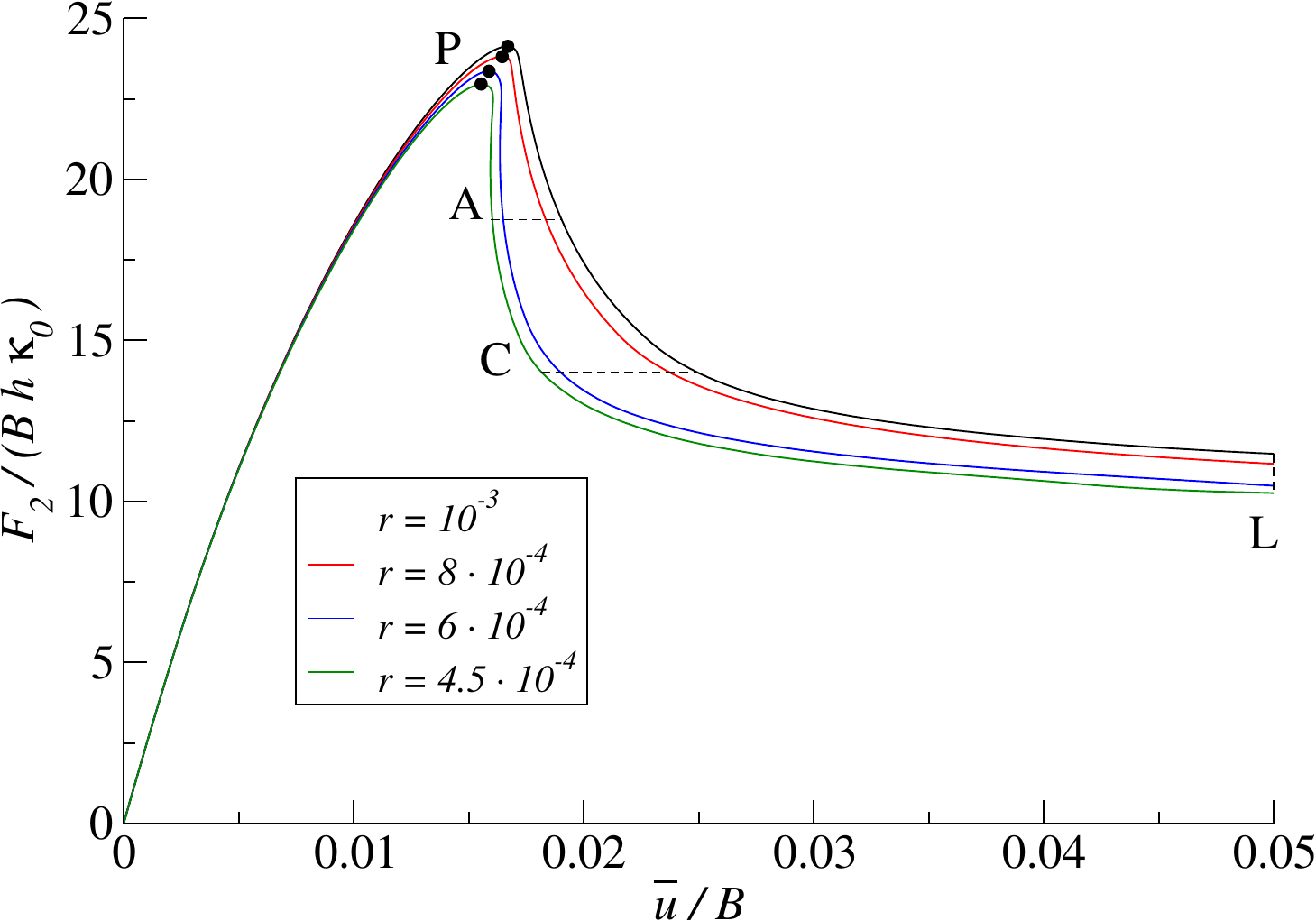} 
\caption{Influence of the internal length $r=\ell/B$ on the load--displacement response for the Matsuoka--Nakai model with softening (Mesh 3). Smaller values of $r$ lead to a more pronounced softening and to an unstable post-peak response, whereas larger values produce a smoother and more stable behavior.}
\label{fig:MN_soft_r}
\end{figure}

The influence of the internal length is further investigated by considering different values of the normalized parameter $r=\ell/B$, namely $r=4.5\cdot10^{-4}$, $6\cdot10^{-4}$, $8\cdot10^{-4}$, and $10^{-3}$. 
Fig.~\ref{fig:MN_soft_r} shows the corresponding load--displacement curves. A strong dependence of the global response on the internal length is observed. For the smallest values, $r=4.5\cdot10^{-4}$ and $r=6\cdot10^{-4}$, the response exhibits a pronounced softening behavior with an unstable post-peak regime. In contrast, for larger values of $r$, the response becomes progressively smoother, with reduced softening and more stable post-peak evolution.

\begin{figure}[t]
	\centering
	\begin{tabular}{cc}
		\includegraphics[trim={2cm 8cm 2cm 0.5cm},clip,width=0.475\textwidth]{figures/riks_ell4,5e-4_mesh3_peak_phi.pdf} 
		&\includegraphics[trim={2cm 8cm 2cm 0.5cm},clip,width=0.475\textwidth]{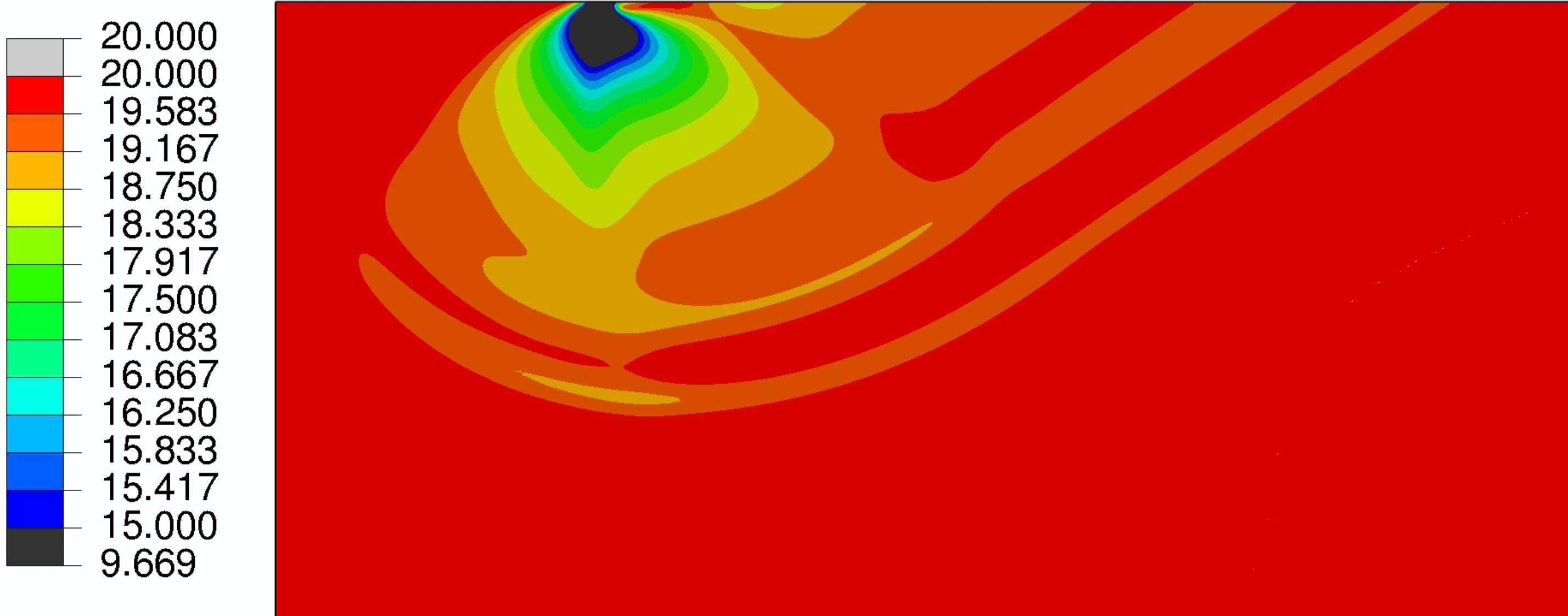} \\
		(a) & (b)\\
		\includegraphics[trim={2cm 8cm 2cm 0.5cm},clip,width=0.475\textwidth]{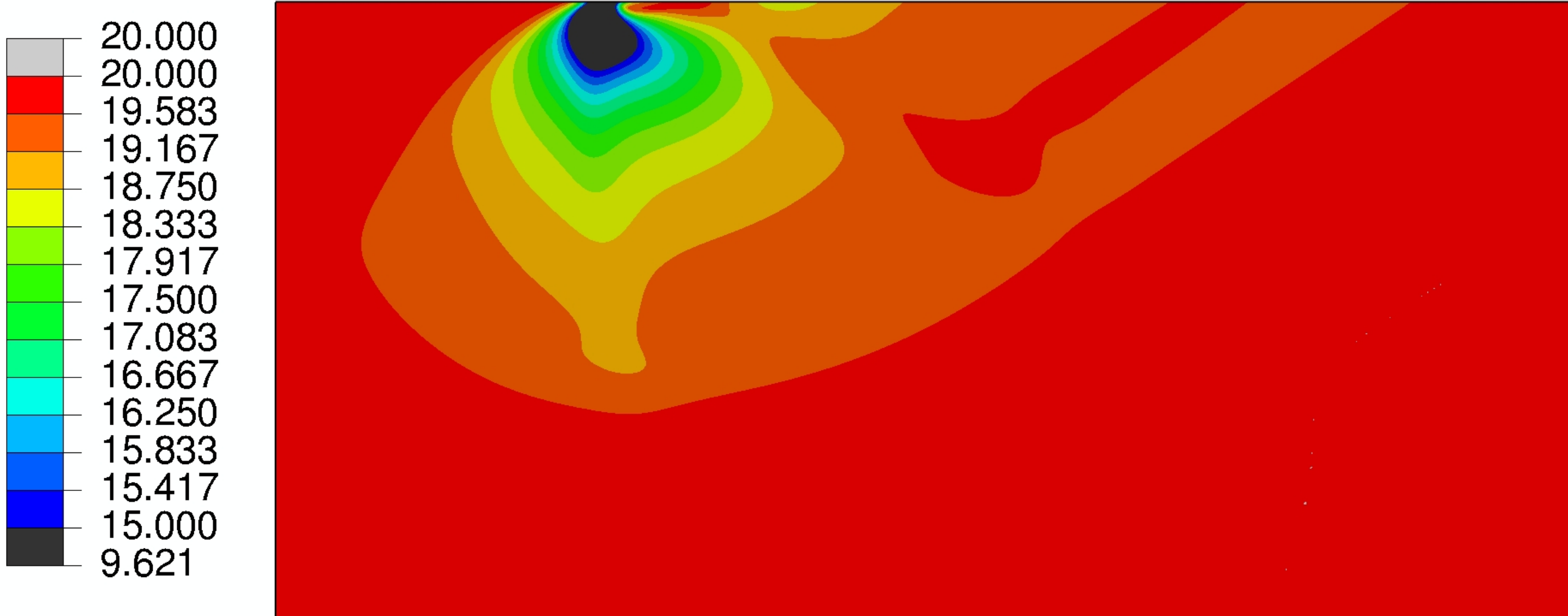} 
		&\includegraphics[trim={2cm 8cm 2cm 0.5cm},clip,width=0.475\textwidth]{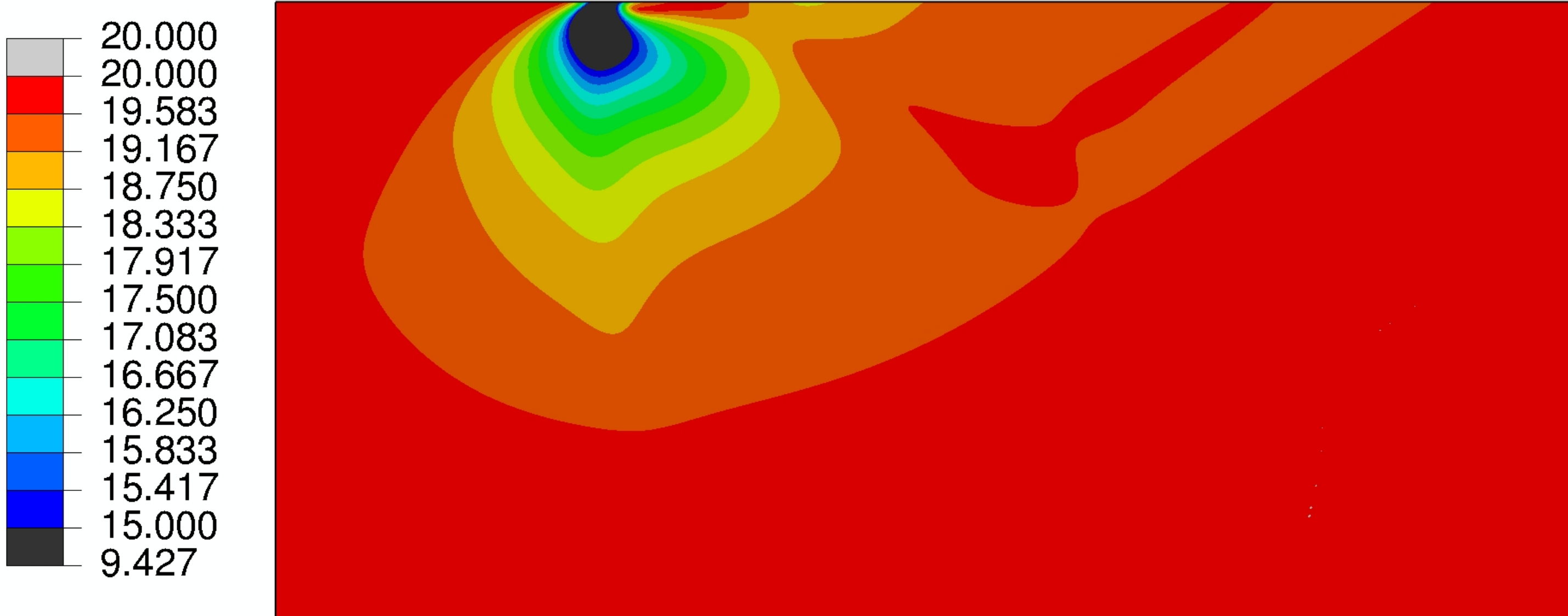} \\
		(c) & (d)  
	\end{tabular}
\caption{Spatial distribution of the angle of shearing resistance $\phi$ for different values of $r=\ell/B$, at peak load (point $P$ in Fig. \ref{fig:MN_soft_r}), for Mesh 3.
(a) $r=4.5\cdot 10^{-4}$, (b) $r=6\cdot 10^{-4}$, (c) $r=8\cdot 10^{-4}$, (d) $r= 10^{-3}$ Smaller values of $r$ lead to the development of multiple shear bands, approximately parallel to the limit analysis solution, while larger values promote simpler mechanisms with two dominant shear bands and a more diffuse intermediate region.}
\label{fig:MN_patterns1}
\end{figure}

\begin{figure}[t]
	\centering
	\begin{tabular}{cc}
		\includegraphics[trim={2cm 8cm 2cm 0.5cm},clip,width=0.475\textwidth]{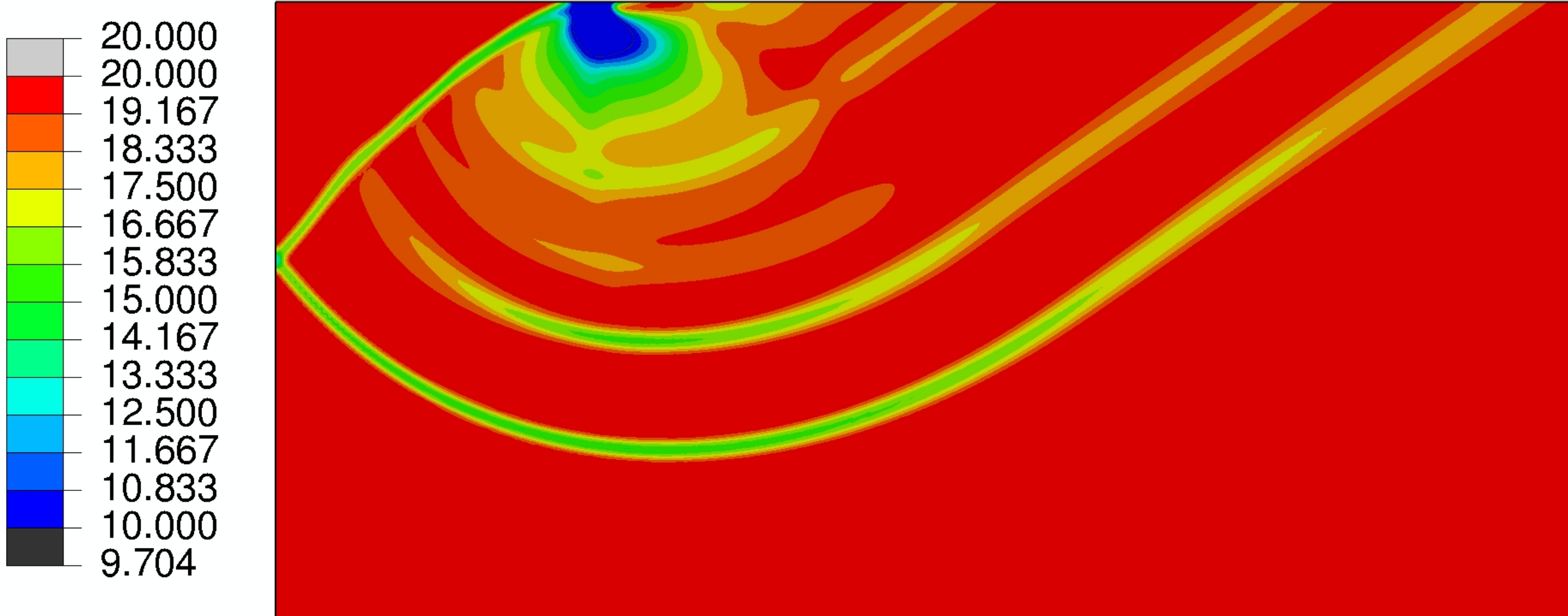} 
		&\includegraphics[trim={2cm 8cm 2cm 0.5cm},clip,width=0.475\textwidth]{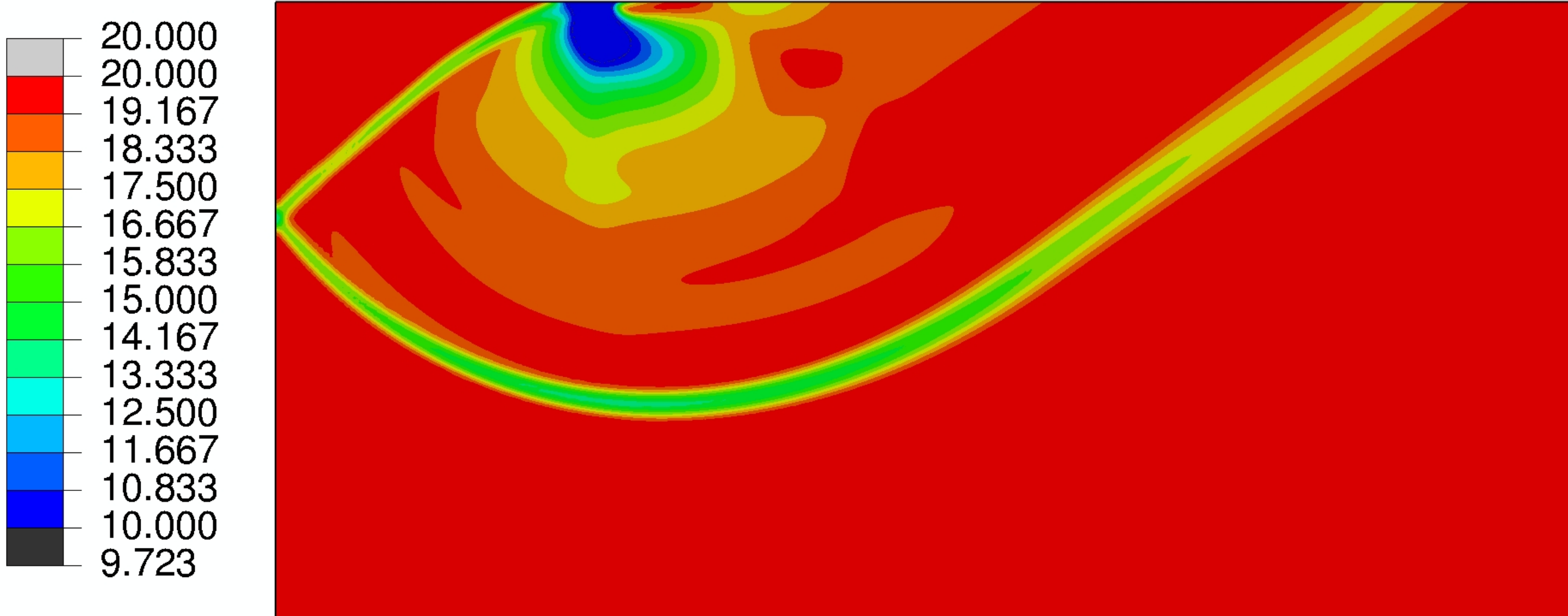} \\
		(a) & (b)\\
		\includegraphics[trim={2cm 8cm 2cm 0.5cm},clip,width=0.475\textwidth]{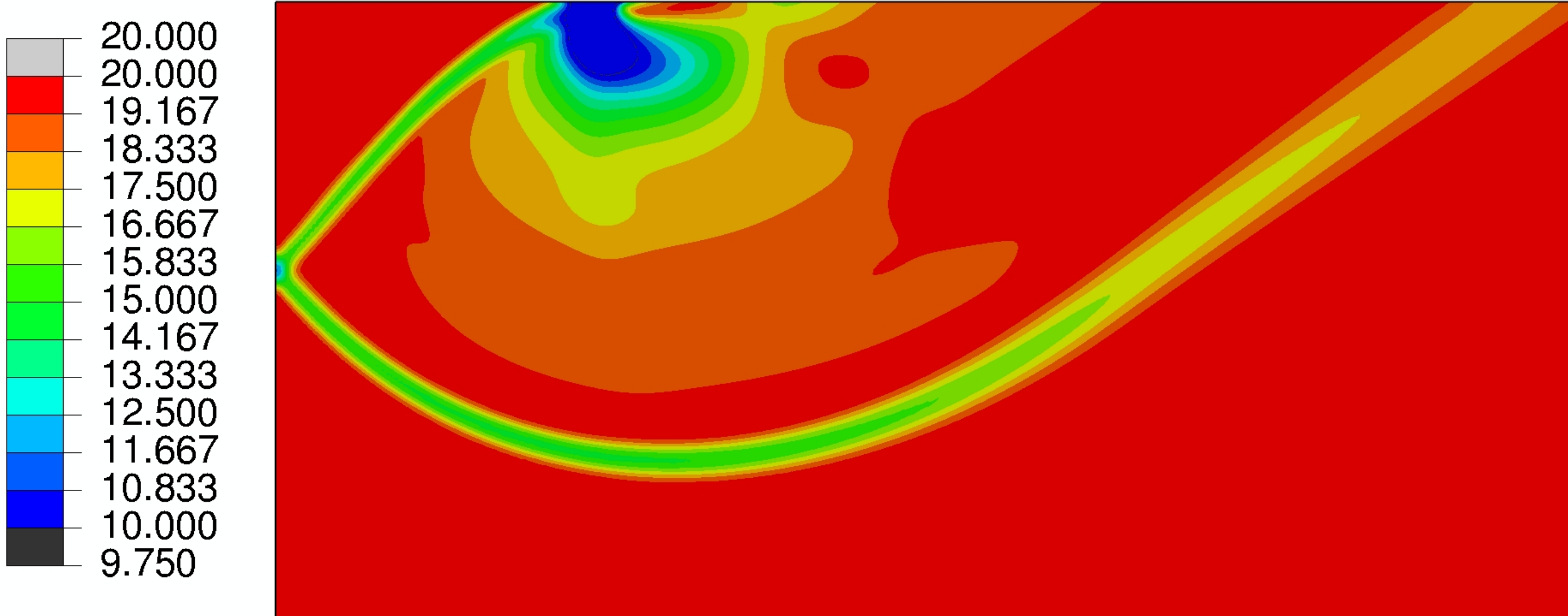} 
		&\includegraphics[trim={2cm 8cm 2cm 0.5cm},clip,width=0.475\textwidth]{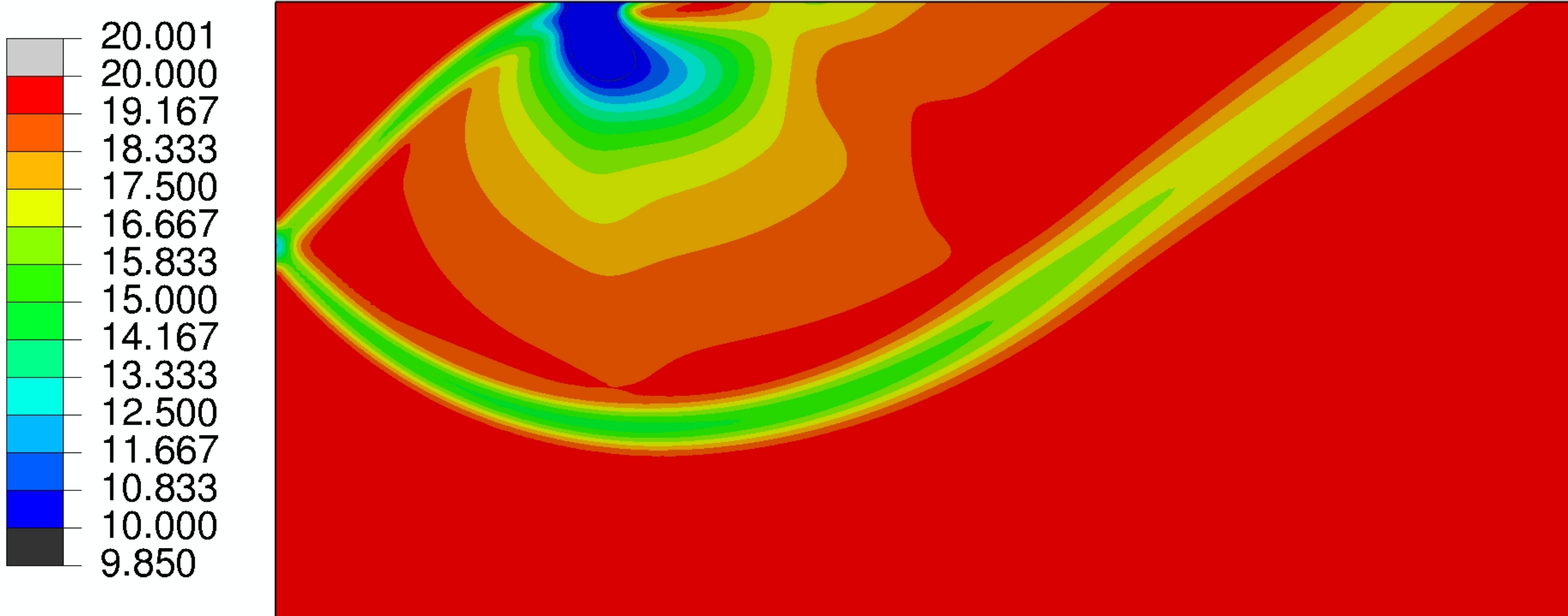} \\
		(c) & (d)  
	\end{tabular}
\caption{Spatial distribution of the angle of shearing resistance $\phi$ for different values of $r=\ell/B$, at intermediate load (point $A$ in Fig. \ref{fig:MN_soft_r}, (at $F_2/(B h \kappa_0) \approx 18.75$) for Mesh 3.
(a) $r=4.5\cdot 10^{-4}$, (b) $r=6\cdot 10^{-4}$, (c) $r=8\cdot 10^{-4}$, (d) $r= 10^{-3}$. 
}
\label{fig:MN_patterns2}
\end{figure}

\begin{figure}[t]
	\centering
	\begin{tabular}{cc}
		\includegraphics[trim={2cm 8cm 2cm 0.5cm},clip,width=0.475\textwidth]{figures/riks_ell4,5e-4_mesh3_F150kPa_phi.pdf} 
		&\includegraphics[trim={2cm 8cm 2cm 0.5cm},clip,width=0.475\textwidth]{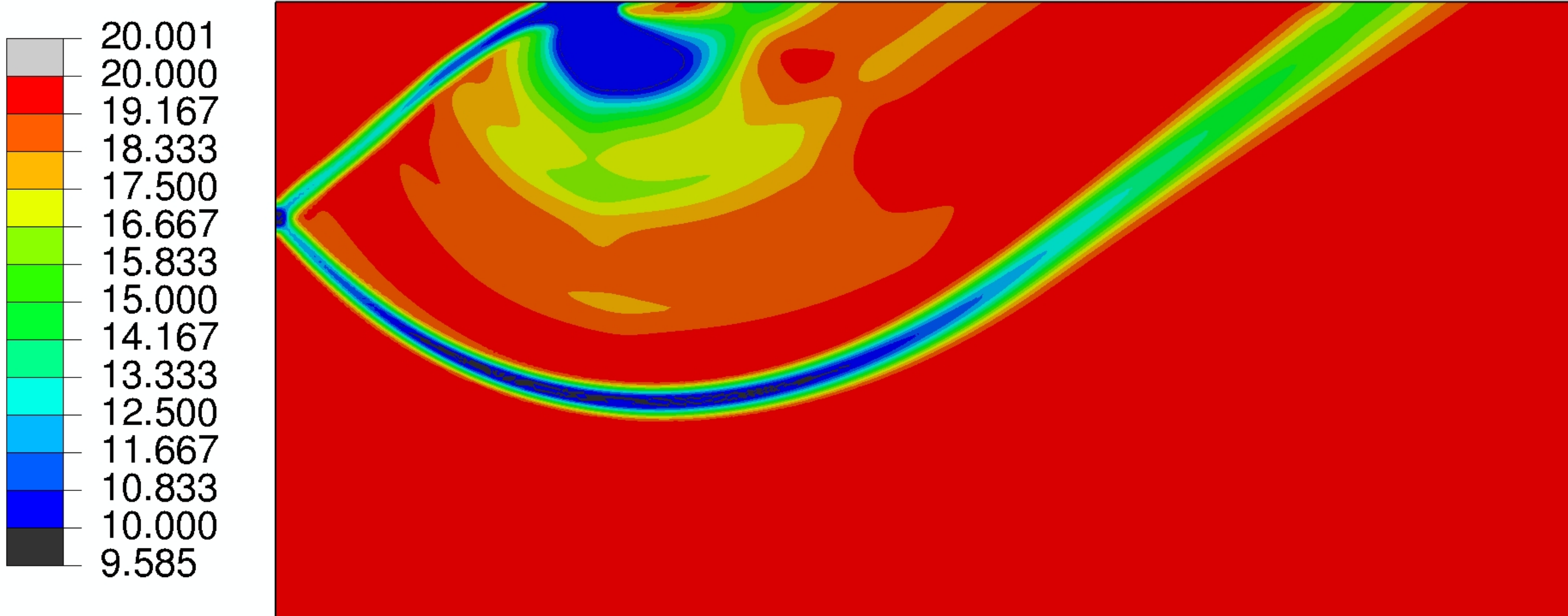} \\
		(a) & (b)\\
		\includegraphics[trim={2cm 8cm 2cm 0.5cm},clip,width=0.475\textwidth]{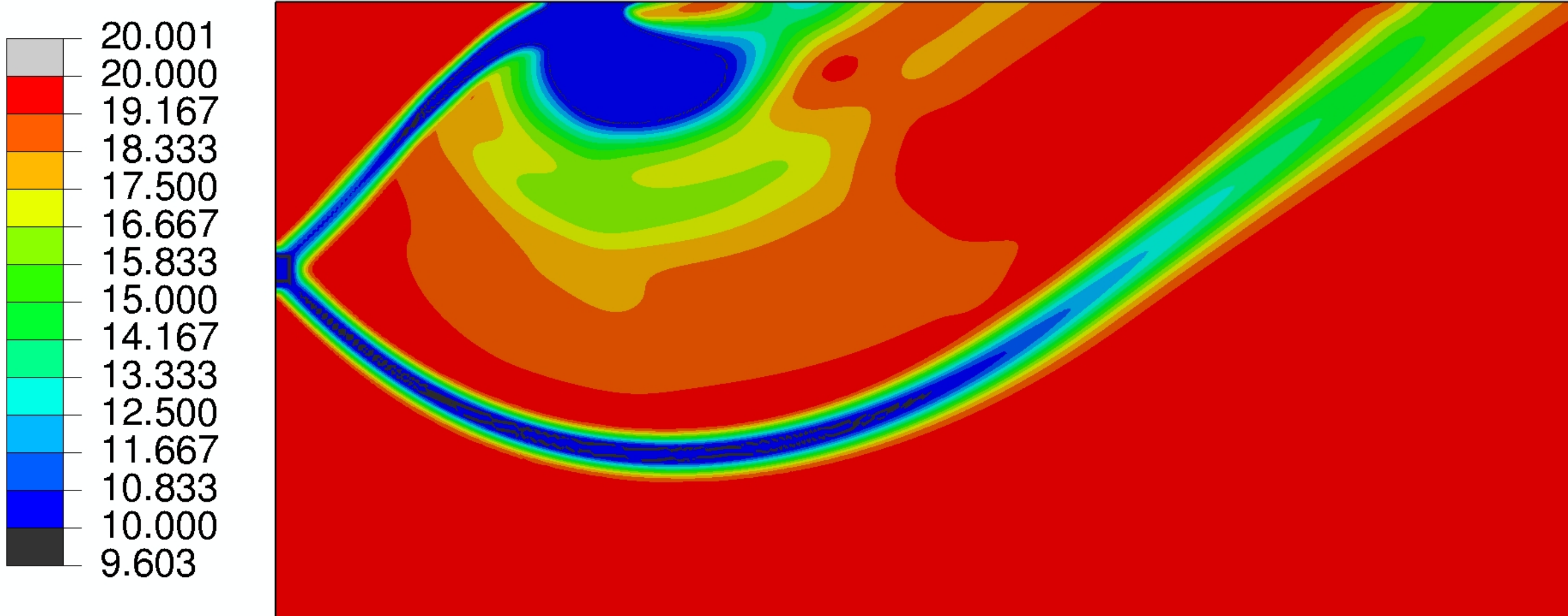} 
		&\includegraphics[trim={2cm 8cm 2cm 0.5cm},clip,width=0.475\textwidth]{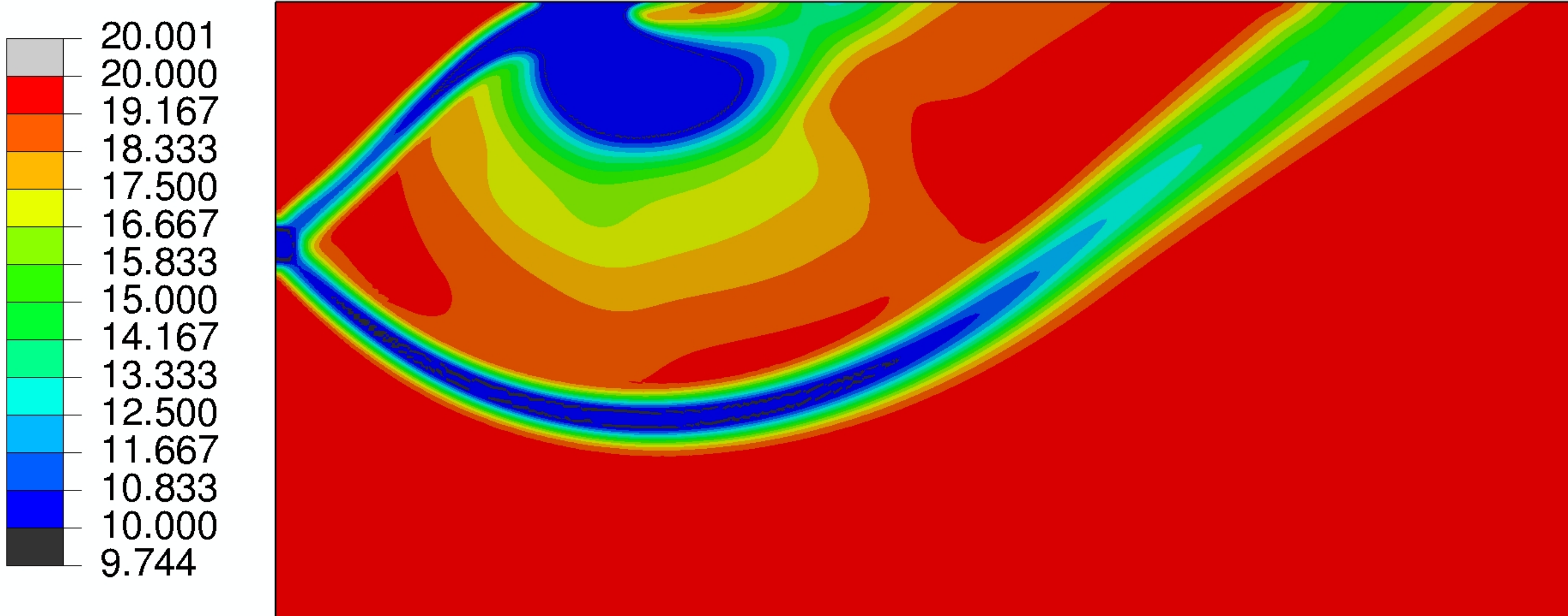} \\
		(c) & (d)  
	\end{tabular}
\caption{Spatial distribution of the angle of shearing resistance $\phi$ for different values of $r=\ell/B$, at intermediate load (point $C$ in Fig. \ref{fig:MN_soft_r}, at $F_2/(B h \kappa_0) \approx 14.$), for Mesh 3.
(a) $r=4.5\cdot 10^{-4}$, (b) $r=6\cdot 10^{-4}$, (c) $r=8\cdot 10^{-4}$, (d) $r= 10^{-3}$ 
}
\label{fig:MN_patterns3}
\end{figure}

\begin{figure}[t]
	\centering
	\begin{tabular}{cc}
		\includegraphics[trim={2cm 8cm 2cm 0.5cm},clip,width=0.475\textwidth]{figures/riks_ell4,5e-4_mesh3_last_phi.pdf} 
		&\includegraphics[trim={2cm 8cm 2cm 0.5cm},clip,width=0.475\textwidth]{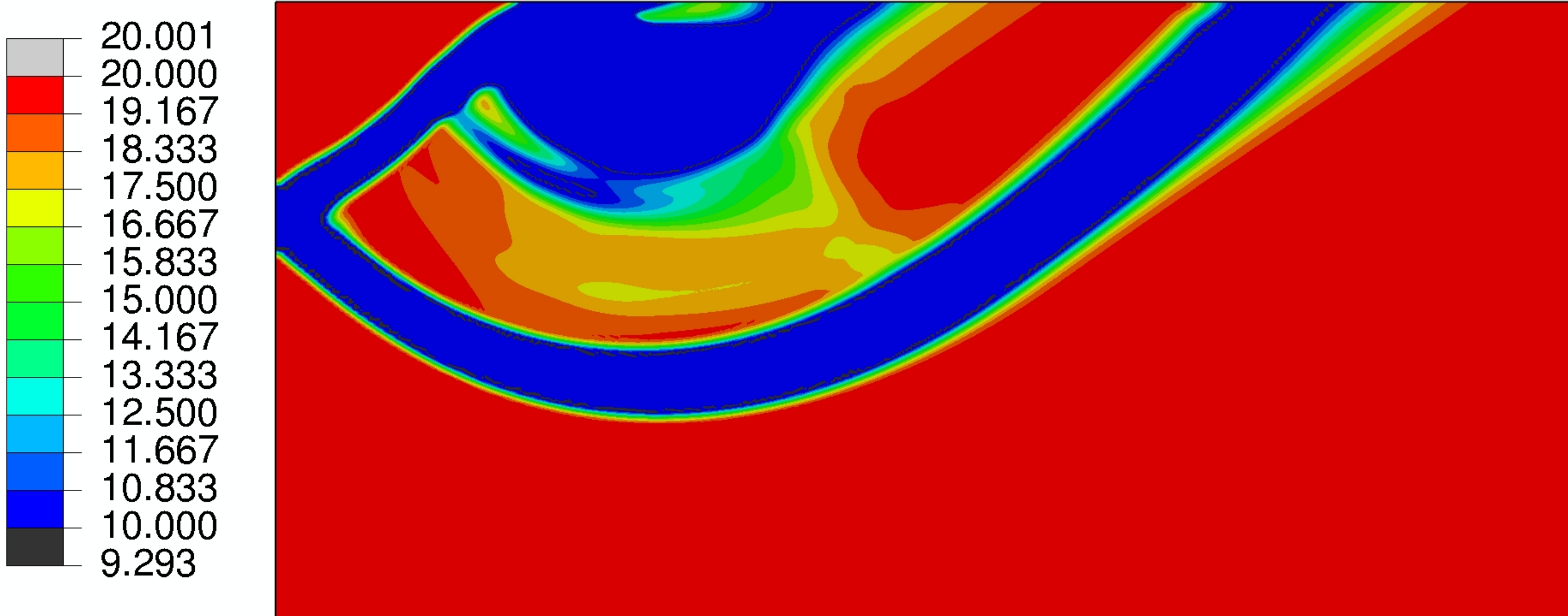} \\
		(a) & (b)\\
		\includegraphics[trim={2cm 8cm 2cm 0.5cm},clip,width=0.475\textwidth]{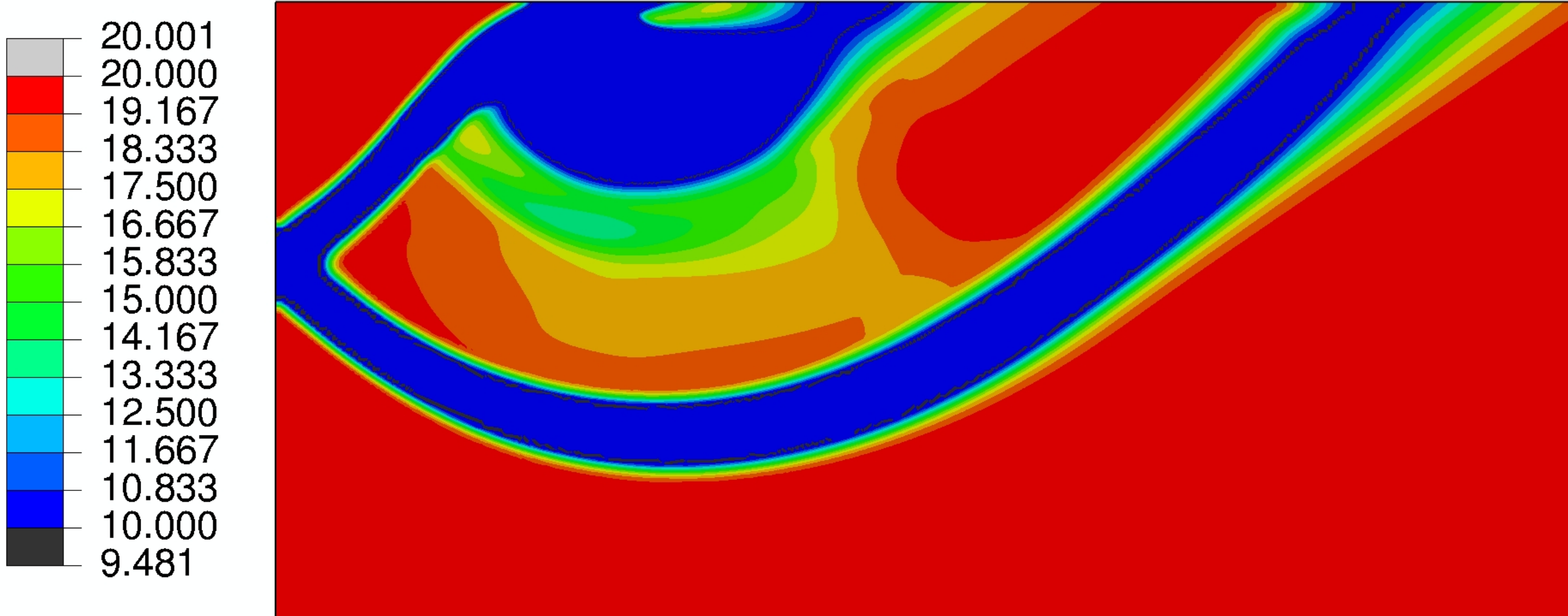} 
		&\includegraphics[trim={2cm 8cm 2cm 0.5cm},clip,width=0.475\textwidth]{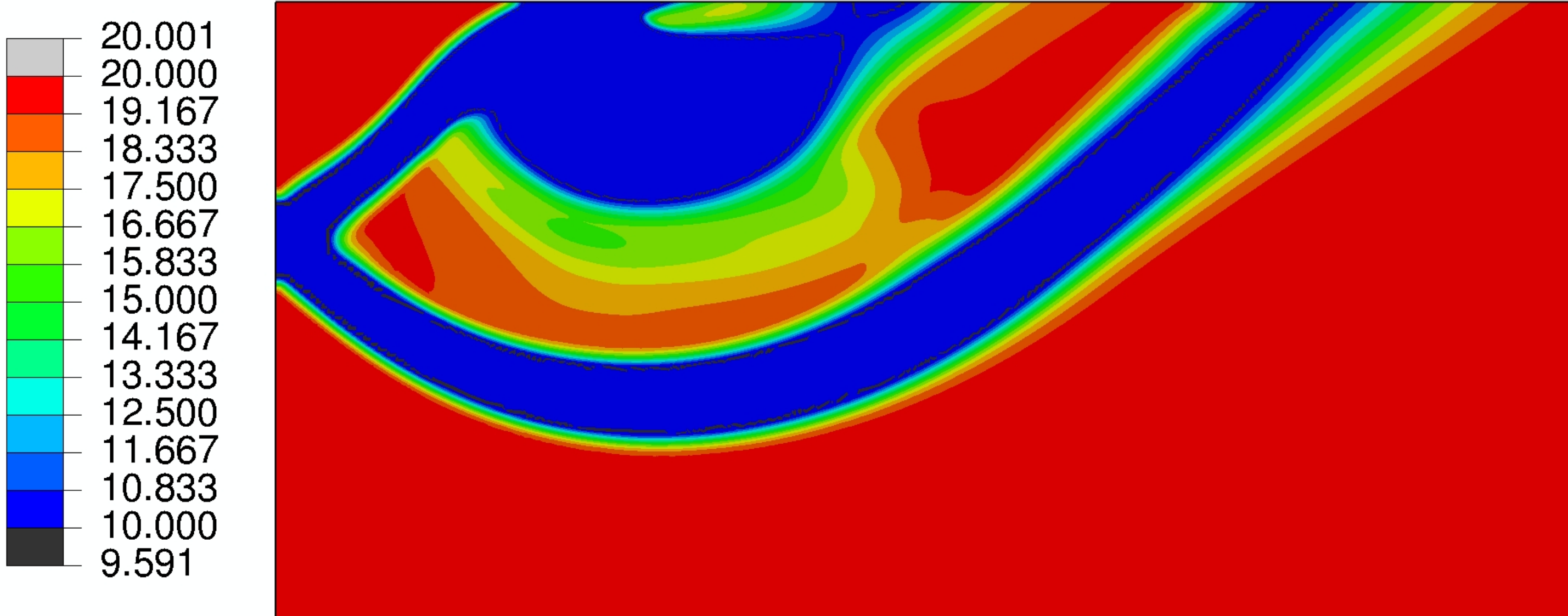} \\
		(c) & (d)  
	\end{tabular}
\caption{
Spatial distribution of the angle of shearing resistance $\phi$ for different values of $r=\ell/B$, at the last increment ($\bar u/H\approx0.05$, point $L$ in Fig. \ref{fig:MN_soft_r}, for Mesh 3.
(a) $r=4.5\cdot 10^{-4}$, (b) $r=6\cdot 10^{-4}$, (c) $r=8\cdot 10^{-4}$, (d) $r= 10^{-3}$.
}
\label{fig:MN_patterns4}
\end{figure}

The corresponding localization patterns are shown in Fig.~\ref{fig:MN_patterns1}-\ref{fig:MN_patterns4} in terms of the spatial distribution of the angle of shearing resistance $\phi$, at different loading stages (point $P$ peak load, 
point $A$, point $C$ and final step $L$).
A qualitative change in the deformation mechanisms is observed as the internal length varies. For the smallest values of $r$, the localization pattern is characterized by the development of multiple shear bands, which appear approximately parallel to those associated with the classical limit analysis solution. These bands form a complex network of interacting localization zones, indicating a highly localized deformation process.
For larger values of $r$, the deformation pattern becomes significantly simpler. The response is dominated by two main shear bands, consistent with the classical limit analysis mechanism, with a relatively diffuse deformation field developing in the region between them.
These results clearly indicate that the internal length does not merely affect the thickness of the shear bands, but fundamentally controls the selection of the deformation mechanisms. In particular, decreasing $r$ promotes the activation of multiple competing localization modes, whereas larger values of $r$ favor the emergence of simpler, more classical failure mechanisms.

\begin{table}[t]
 	\centering
 	\begin{tabular}{c|c c c}
        \hline
 		\hline  		
 		$r$    &   loading &    total  & total arc\\ 
 		&           steps     &  iterations  & length \\
		\hline 		
 		$4.5\cdot 10^{-4}$ 	& 592 		& 1750	& 2.916 \\ 
        $6\cdot 10^{-4}$ 	& 586		& 1686	& 2.886\\ 
        $8\cdot 10^{-4}$ 	& 539		& 1484	& 2.651\\ 
        $10^{-3}$ 	& 533		& 1264 &	2.621\\ 
		\hline 
		\hline
 	\end{tabular} 
 	\caption{Computational cost for the Matsuoka--Nakai model with softening of the angle of shearing resistance, obtained using the Riks procedure for Mesh~3 and different values of $r=\ell/B$. Decreasing $r$ increases the number of load increments, total iterations, and arc-length, reflecting the increased complexity of the nonlinear response. However, the solution procedure remains stable and very efficient.}
	    \label{tab:MNcost_r}
 \end{table}

The computational cost associated with different values of the internal length is reported in Tab.~\ref{tab:MNcost_r} for Mesh~3.
In contrast to the previous results, which highlighted the mesh-independent behavior, the present data show a clear dependence on the internal length. Decreasing $r$ leads to a progressively more demanding nonlinear problem, as indicated by the increase in the number of load increments, total iterations, and total arc-length.
This trend reflects the stronger localization and more pronounced softening response observed for smaller values of $r$, which result in a more complex and unstable equilibrium path.
The average number of iterations per increment increases from $2.37$ for $r=10^{-3}$ to $2.96$ for $r=4.5\cdot10^{-4}$, while the average arc-length increment remains close to the maximum allowable value of $5\cdot10^{-3}$ for all cases. This indicates that, although the problem becomes more demanding as $r$ decreases, the nonlinear solution procedure remains stable and very efficient throughout.

Overall, the results obtained for the Matsuoka--Nakai model with softening highlight the capability of the proposed formulation to robustly capture highly nonlinear responses, including pronounced softening and unstable behavior. The convergence of both global quantities and localization patterns confirms the mesh-objective nature of the approach, even in the presence of complex and evolving deformation mechanisms.
At the same time, the analysis of different values of the internal length $\ell$ shows that the model does not merely regularize the solution, but actively governs the selection of the deformation modes. In particular, smaller values of $\ell$ promote the activation of multiple interacting shear bands and a strongly unstable response, whereas larger values lead to simpler mechanisms, closer to the classical limit analysis solution.

These findings demonstrate that the proposed framework provides a consistent and efficient tool for the analysis of strain localization in softening plasticity, combining robustness, mesh objectivity, and the ability to reproduce a wide range of physically meaningful deformation patterns without modifying standard constitutive models.
\FloatBarrier

\section{Discussion: How the model works} 
The mechanism by which the deformable Cosserat model influences localization can be understood by examining when the micro-continuum is activated and how it affects the response.

The additional (i.e., non-standard) micro-stress $\bT_\mathrm{micro}$ does not depend on macroscopic deformation alone, but on the \emph{mismatch} between the directors $\bd_i$ and the material line elements $\ba_i$, governed respectively by the tensor $\boldeta$ and by the displacement gradient. This \emph{mismatch} is described by the second-order tensor $\bchi$, which measures both the corresponding relative deformation through its symmetric part
and the relative rotation through its skew-symmetric part:
\begin{equation}
\bchi = \frac{\partial \bu}{\partial \bx} - \boldeta^T .
\end{equation}
This feature shares a clear conceptual link with the standard rigid Cosserat model, where the micro-stress depends on the \emph{relative rotation} between a \emph{rigid} director triad \cite{PL2022a} and the macro-continuum. In the adopted formulation \cite{miles}, this idea is extended to a \emph{deformable} director field, 
so that the micro-continuum response is governed by a relative measure that includes not only rotations, but also deformations\footnote{It should be noted, however, that the deformable Cosserat model adopted here does not include the rigid Cosserat model as a special case, although both reduce to the classical Cauchy continuum in the limit as the material length $\ell \rightarrow 0$ \cite{PR2026}.}. A further analogy with the rigid Cosserat model concerns the higher-order kinematics. In the rigid setting, the curvature measure is defined by the gradient of the total director rotation. Similarly, in the present formulation, the curvature is determined by the gradient of the total director field $\boldeta$. Thus, in both cases, the micro-stress is governed by a \emph{relative} measure, whereas the micro-couples are determined by the gradients of the \emph{total} micro-kinematic fields. 

In spite of these analogies, the difference between the two models is fundamental. In both the rigid and deformable Cosserat models, the strain energy is typically assumed to admit an additive decomposition into a macro-part and a micro-part, so that the total stress is correspondingly split into the sum of macro- and micro-contributions. In the rigid Cosserat model the relative measure governing the micro-response is purely rotational and therefore it is influenced only by the skew-symmetric part of the displacement gradient; the associated micro-stress is correspondingly skew-symmetric. As a consequence, the symmetric part of the total stress is associated only with the macro-stress, whereas the skew-symmetric part is associated only with the micro-stress. 
In contrast, in the deformable Cosserat model $\bchi$ is a full second-order tensor. Hence, the symmetric part of the displacement gradient contributes to both the standard macro-constitutive response, through its rate, and to the tensor $\bchi$, which drives the micro-response. Accordingly, no sym/skw separation holds in the deformable case: the micro-stress generally has both symmetric and skew-symmetric parts, whereas the macro-stress remains symmetric, consistently with the assumed strain energy function.

If no explicit constitutive coupling is introduced, the interaction between the macro- and micro-continua is governed only by the  balance of linear momentum. In both models, the balance of director momentum couples the micro-stress to the divergence of the micro-couples, whereas the balance of linear momentum involves the total stress. In the rigid Cosserat case, since the micro-stress is purely skew-symmetric, the coupling enters the balance of linear momentum only through the skew-symmetric part of the total stress. 
In the deformable Cosserat case, by contrast, the micro-stress generally includes both symmetric and skew-symmetric parts. Therefore, the coupling induced by the balance equations is not confined to the skew-symmetric part of the total stress, but also acts through its symmetric part. The resulting stronger coupling through the balance equations is sufficient, in the present model, to make it possible to retain a completely standard macro-elastoplastic description, in which dissipation depends only on the symmetric macro-stress and on the corresponding symmetric plastic deformation rates, while still avoiding mesh-dependent localization.

In the rigid Cosserat case, the weaker coupling through the kinematics and the balance equations typically makes it necessary, in softening elastoplasticity, to formulate constitutive models in which the micro-stress, and especially the micro-couples, enter explicitly into the yield function and plastic potential \cite{muhlhaus1987thickness, PL2022a, russo2020thermomechanics} in order to activate gradient effects \cite{russo2020thermomechanics} and mitigate mesh dependence when localization occurs.
In contrast, the stronger coupling between macro- and micro- continua in the deformable Cosserat model can avoid mesh-dependent localization without plasticity in the micro-continuum.

A further important feature of the deformable Cosserat model concerns the behaviour of the micro-continuum under homogeneous deformation. In this case, the micro-continuum becomes inactive: the balance of director momentum \eqref{eq:dirmomentum} enforces $\bd_i=\ba_i$ and hence $\bchi=\mathbf{0}$, while the balance of linear momentum involves only the macro-stress. The response, therefore, coincides with that of the standard macro-continuum.

\begin{figure}[t]
\centering
\begin{tabular}{cc}
\includegraphics[trim={2cm 3cm 2cm 0.5cm},clip,width=0.475\textwidth]{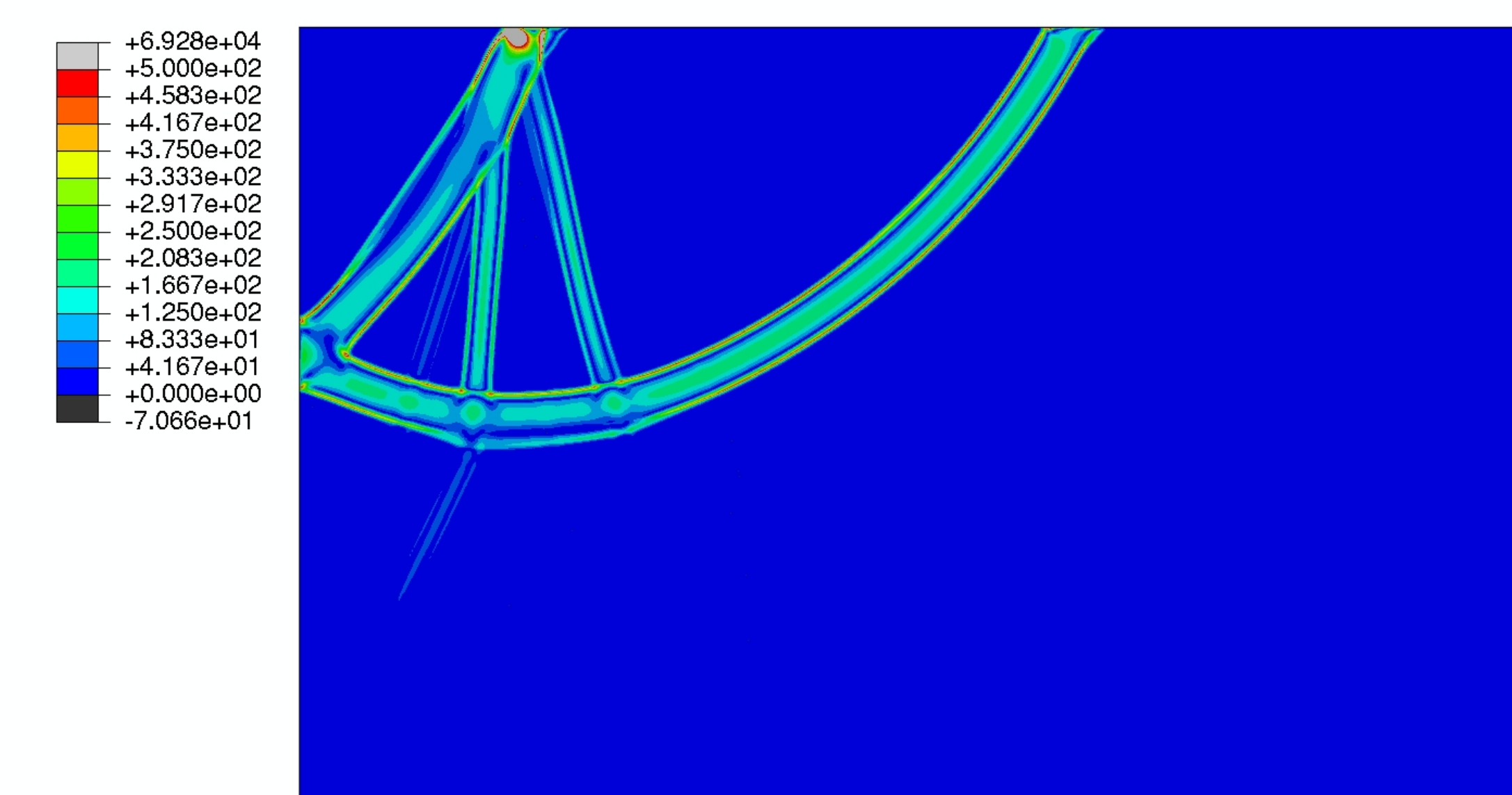}
&\includegraphics[trim={2cm 3cm 2cm 0.5cm},clip,width=0.475\textwidth]{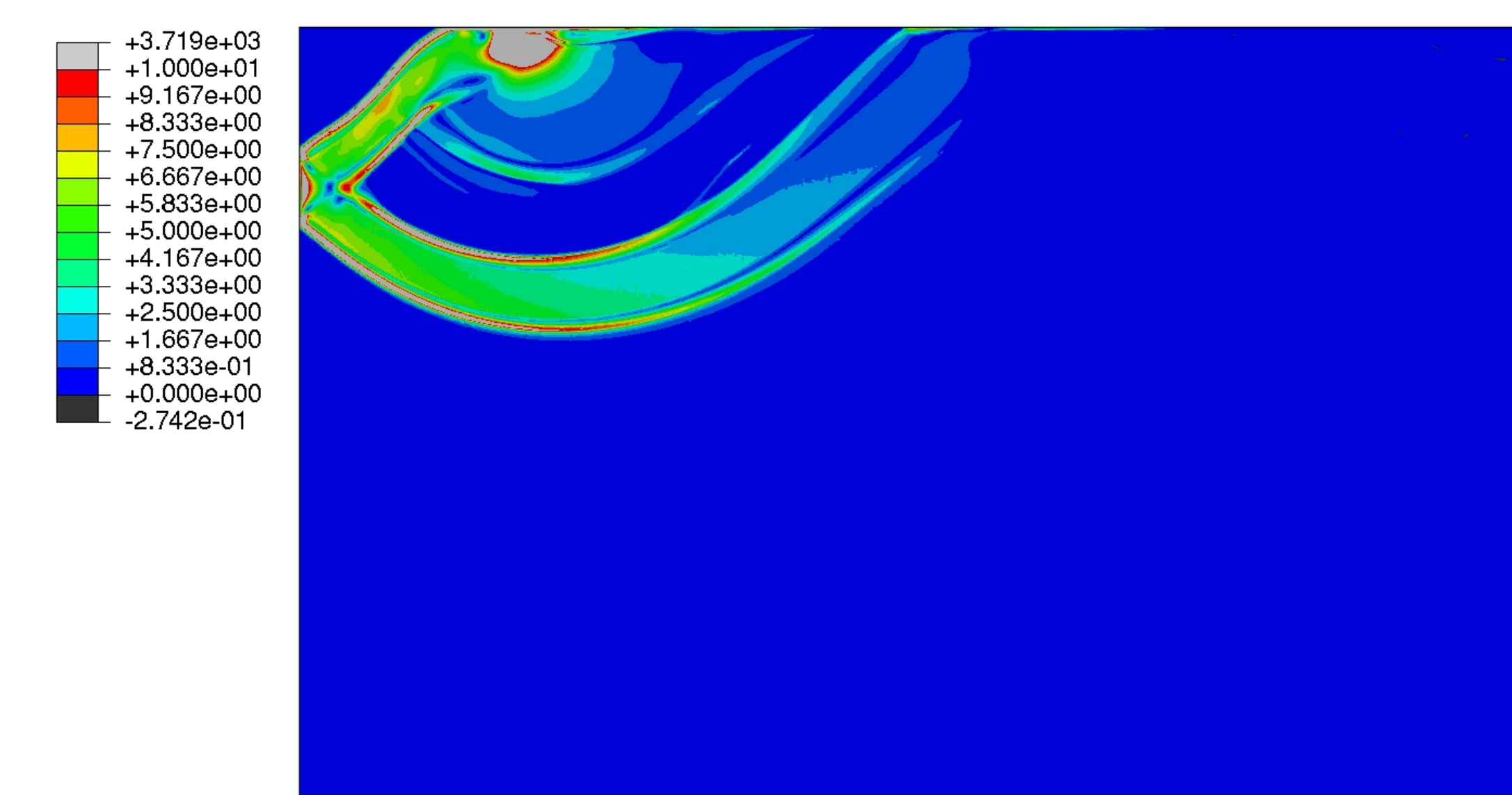} \\
(a) & (b)
\end{tabular}
\caption{
Spatial distribution of the modulus of the micro-stress $\sqrt{\bT_\mathrm{micro}:\bT_\mathrm{micro}}$ (in $\mathrm{kPa}$) at the end of loading, for $r=10^{-3}$. (a) Tresca soil with softening. (b) Matsuoka--Nakai soil with softening. In both cases, the micro-stress is negligible in the bulk of the domain and becomes significant only in narrow regions localized along the edges of the shear bands. Since $\bT_\mathrm{micro}$ is proportional to the tensor $\bchi = \partial \bu/\partial \bx - \boldeta^T$, the figure provides a direct indication of the regions where $\bchi \neq \mathbf{0}$, i.e. where the directors $\bd_i$ differ from the material line elements $\ba_i$.
}
\label{fig:Tmacro_patterns}
\end{figure}

This is exactly what is observed in Fig.~\ref{fig:Tmacro_patterns}: the micro-continuum becomes significant only in localized zones, where strong spatial variations of the displacement gradient induce non-zero values of $\bchi$, i.e., a mismatch between the directors $\bd_i$ and the material line elements $\ba_i$, driven by the displacement gradient. Accordingly, $\bT_\mathrm{micro}$ is negligible in the bulk and becomes significant only in narrow regions near the shear bands.

If the deformation is not homogeneous, the balance of director momentum \eqref{eq:dirmomentum} requires a non-vanishing micro-stress to be satisfied. The micro-couples entering Eq.~\eqref{eq:dirmomentum} are governed by the gradients of the director field, which follow from the kinematic relation
\begin{equation}
\boldeta = \left(\frac{\partial \bu}{\partial \bx} - \bchi\right)^T \,,  
\quad\mbox{ i.e.,}\quad\eta_{ij} = u_{j,i} - \chi_{ji} \,,
\end{equation}
and therefore satisfy
\begin{equation}
\eta_{ij,k} = u_{j,ik} - \chi_{ji,k}\,.
\end{equation}
Hence, the director gradients depend on both the gradient of the  tensor $\bchi$ and the second gradient of the displacement field. When localization starts to develop, the associated strong spatial variations induce non-zero curvatures of the director field and therefore generate micro-couples. Their divergence acts as the driving term in the balance of director momentum and must be balanced by a non-zero micro-stress $\bT_\mathrm{micro}$, which in turn implies $\bchi \neq \mathbf{0}$.

In the deformable Cosserat model, the energetic contributions associated with the micro-continuum become relevant only when strong spatial variations of the displacement gradient develop.
 The regularizing effect therefore arises naturally from the kinematic structure of the model, through the mismatch between $\bd_i$ and $\ba_i$, while the macro-response remains entirely standard. In a standard softening continuum, no comparable energetic mechanism opposes the collapse of localization into an arbitrarily narrow zone, so that the computed band thickness is ultimately governed by the discretization. Here, by contrast, sharply localized kinematics require large values of the tensor $\bchi$ and large director curvatures, each of which carry an energetic cost. The collapse into a zero-width zone therefore becomes energetically unfavorable, and the localization pattern develops over a finite width selected by the internal length $\ell$ rather than by the mesh.

It is emphasized that the present mechanism is fundamentally different from approaches in which an internal length is introduced at the constitutive level through a length-dependent effective measure governed by additional differential equations, such as phase field, gradient plasticity or other nonlocal models.

This also clarifies the role of the internal length $\ell$. Rather than acting as a numerical regularization parameter, $\ell$ governs the characteristic width and spatial extent of the regions where these gradients are energetically penalized. For a fixed value of $\ell$, the formulation yields mesh-independent solutions that converge towards well-defined localization patterns, while varying $\ell$ leads to qualitatively different responses in both the global behavior and the associated shear band structures. Thus, the role of $\ell$ is not merely to smooth the numerical solution, but to select a physically meaningful localization width and, consequently, a well-posed post-localization response. These observations support interpreting $\ell$ as a material parameter, to be identified from experimental measurements of shear-band thickness and spatial patterns.

\section*{Conclusions}

This paper has investigated strain localization in softening plasticity by means of the small-strain version of the deformable Cosserat model proposed in \cite{miles}, and has presented a finite element implementation of the resulting formulation.

The deformable Cosserat model consists of two interacting continua: a macro-continuum and a micro-continuum. Their constitutive responses are uncoupled, whereas their interaction is governed by the balances of linear and director momentum. The macro-continuum governs the entire dissipative response and is described by a standard elastoplastic constitutive law. The micro-continuum is non-dissipative and contributes only through additional energetic terms associated with the mismatch between the directors and  material line elements, as well as with the corresponding director curvatures. It therefore becomes relevant only when sufficiently strong spatial variations of the deformation field develop, and governs the structure of the resulting localization patterns.

The main advantage of the proposed approach is that, because of this constitutive structure, standard elastoplastic constitutive laws developed for the classical Cauchy continuum can be employed for the macro-continuum without any modification of either the stress update algorithm or the consistent tangent operator.

A key aspect of the proposed formulation is that the stronger kinematic coupling provided by the deformable-director setting is sufficient to influence localization without introducing plasticity at the micro-level. This distinguishes the present approach both from the rigid Cosserat model, where the coupling is weaker and typically requires a more direct constitutive involvement of micro-stresses and micro-couples, and from gradient-enhanced or nonlocal formulations, in which the internal length is introduced directly through the constitutive description.

The resulting finite element implementation is particularly simple. The macro-response is obtained from a standard constitutive update driven by the symmetric part of the displacement-gradient rate $\bD$, exactly as in a standard Cauchy continuum. All the quantities associated with the micro-continuum, namely the  relative deformation and rotation tensor $\bchi$, the director curvatures, the micro-stresses, and the micro-couples, are then determined explicitly from the nodal values of  $\boldeta$ and $\bu$ through the finite element operators and linear elastic relations. Hence, once the macro-stress and the corresponding consistent elastoplastic tangent have been obtained from a standard constitutive routine, the remaining part of the formulation reduces to the assembly of linear contributions. A general finite element algorithm has been presented, together with the corresponding FE operators for the plane strain case.
This opens the possibility of extending the approach to a broad class of standard constitutive models without reformulating their constitutive integration algorithms.

The performance of the proposed approach has been assessed through a boundary value problem known to be particularly demanding, namely the shallow foundation problem under plane strain conditions. Tresca and Matsuoka--Nakai plasticity models have been considered, including both perfectly plastic and softening responses. In the softening regime, for a fixed value of the material length $\ell$, the proposed formulation yields mesh-independent solutions with localization patterns that converge upon mesh refinement. In particular, for the Matsuoka--Nakai model, convergence is also obtained in cases exhibiting unstable post-peak behaviour. At the same time, varying $\ell$ leads to qualitatively different responses, both in terms of global behaviour and of the associated shear-band profiles.

These results support a clear interpretation of the role of the internal length $\ell$. It should not be regarded merely as a numerical regularization parameter, but rather as a material parameter governing the evolution of converged localization patterns. In this sense, the proposed model provides a physically meaningful framework for the description of strain localization in softening plasticity, while preserving the full simplicity and generality of standard elastoplastic constitutive updates.

\section*{Acknowledgements}
The finite element code Simulia Abaqus has been run at the Department of Civil, Environmental, Architectural Engineering and Mathematics, University of Brescia, Italy, under an academic license
\appendix

\section{Finite Element operators for the plane strain formulation}
\label{Apn_FEOper}

The matrix $\bH_r$ that contains the derivatives of the shape functions for displacements and geometry with respect to the intrinsic coordinate system $\br$, can be expressed as
\begin{equation}
\bH_r = \begin{bmatrix}
\der{N^{(1)}}{r} & .. &\der{N^{(n)}}{r} \\
\der{N^{(1)}}{s} & .. & \der{N^{(n)}}{s} \\
\end{bmatrix} \,.
\end{equation}
In the same way, the matrix $\bH_{\eta r}$ that contains the derivatives of the shape functions for the tensor components $\eta_{ij}$ with respect to the intrinsic coordinate system $\br$, can be expressed as
\begin{equation}
\bH_{\eta r} = \begin{bmatrix}
\der{N_{\eta}^{(1)}}{r} & .. &\der{N_{\eta}^{(n_{\eta})}}{r} \\
\der{N_{\eta}^{(1)}}{s} & .. & \der{N_{\eta}^{(n_\eta)}}{s} \\
\end{bmatrix}\,.
\end{equation}
The jacobian can be computed as
\begin{equation}
\bJ = \bH_r \hat{\bX}\,,
\end{equation}
where
\begin{equation}
\hat{\bX} = \begin{bmatrix}
\hat{x}^{(1)}_{1} & \hat{x}^{(1)}_{2} \\
.. &  ..\\
\hat{x}^{(n)}_{1} & \hat{x}^{(n)}_{2}  \\
\end{bmatrix}\,.
\end{equation}
The nodal displacements are ordered in the vector $ \hat \bu $ of $2n$ components
\begin{equation}
\hat \bu  =\begin{bmatrix}
\hat{u}^{(1)}_1 &
\hat{u}^{(1)}_2 &
.. &
\hat{u}^{(n)}_1 &
\hat{u}^{(n)}_2  
\end{bmatrix}^T \,.
\nonumber
\end{equation}
The components of the tensor $\eta_{ij}$ at the nodes are stored in an array $\hat \boldeta$ of $4 n_\eta$ components
\begin{equation}
\hat \boldeta  =\begin{bmatrix}
\hat{\eta}^{(1)}_{11} &
\hat{\eta}^{(1)}_{22} &
\hat{\eta}^{(1)}_{12} &
\hat{\eta}^{(1)}_{21} &
.. &
\hat{\eta}^{(n_\eta)}_{11} &
\hat{\eta}^{(n_\eta)}_{22} &
\hat{\eta}^{(n_\eta)}_{12} &
\hat{\eta}^{(n_\eta)}_{21} 
\end{bmatrix}^T \,.
\end{equation}
The second-order identity tensor, in array notation, has the form
\begin{equation}
\bI =
\begin{bmatrix}
1 & 1 & 1& 0 & 0 
\end{bmatrix}^T \,,
\end{equation}
while the operators $ \mathcal{I}_d$, to compute the deviatoric part, and 
$\mathcal{S}$, and $\mathcal{P}$, to convert arrays between compact symmetric, engineering and full tensorial notations are respectively equal to:
\begin{equation}
    \mathcal{I}_d=\begin{bmatrix}
        \frac{2}{3} & - \frac{1}{3}  &  -\frac{1}{3} & 0 & 0\\
        - \frac{1}{3} &  \frac{2}{3} &  - \frac{1}{3} & 0 & 0\\
        - \frac{1}{3} & - \frac{1}{3}&  \frac{2}{3} & 0 & 0\\
        0 & 0 &0 & 1 & 0\\
        0 & 0 &0 & 0 & 1\\
    \end{bmatrix},\quad
\mathcal{S} =
\begin{bmatrix}
1 & 0 & 0 & 0 & 0\\
0 & 1 & 0 & 0 & 0\\
0 & 0 & 1 & 0 & 0\\
0 & 0 & 0 & 1 & 1
\end{bmatrix}
\, \quad 
\mathcal{P} =
\begin{bmatrix}
1 & 0 & 0 & 0\\
0 & 1 & 0 & 0\\
0 & 0 & 1 & 0\\
0 & 0 & 0 & 1\\
0 & 0 & 0 & 1
\end{bmatrix}
\, .
\end{equation}
The $\mathcal{A}$ and $\mathcal{A}_\mathrm{sym}$ map the components of $\bH$  into the $5\times 2n$ matrices and have the form
\begin{equation}
\mathcal{A} =
\left[
\begin{matrix}
(\bH)_{11} & 0               &..               &(\bH)_{1n}  & 0 \\
0                & (\bH)_{21}& ..   &  	0            &(\bH)_{2n} \\
0                &0                & 	..                     &0                 & 0                \\
(\bH)_{21} &0&.. &	 (\bH)_{2n} &0 \\
0
&(\bH)_{11} &
.. 
& 0 &(\bH)_{1n}
\end{matrix}
\right], \quad
\mathcal{A}_\mathrm{sym}=
\left[
\begin{matrix}
(\bH)_{11} & 0               &..               &(\bH)_{1n}  & 0 \\
0                & (\bH)_{21}& ..   &  	0            &(\bH)_{2n} \\
0                &0                & 	..                     &0                 & 0                \\
\tfrac{1}{2}(\bH)_{21}
&\tfrac{1}{2}(\bH)_{11} &
.. 
& \tfrac{1}{2}(\bH)_{2n} &\tfrac{1}{2}(\bH)_{1n}\\
\tfrac{1}{2}(\bH)_{21} &\tfrac{1}{2}(\bH)_{11}&.. &	 \tfrac{1}{2}(\bH)_{2n} &\tfrac{1}{2}(\bH)_{1n} \\
\end{matrix}
\right]
\,.
\end{equation}
The operator 
 $\mathcal{N}_\tau$  maps the shape functions $\bN_{\eta}$ into a $5\times  4n_\eta$ matrix,  which is defined by
\begin{equation}
  \mathcal{N}_\tau=
  \begin{bmatrix}
      N_{\eta}^{(1)} & 0 & 0 & 0 & .. & N_{\eta}^{(n_\eta)} & 0 & 0 &0\\
      0 & N_{\eta}^{(1)} & 0 & 0  & .. & 0 &N_{\eta}^{(n_\eta)} & 0 &0\\
      0 & 0 & 0 & 0 & .. & 0 & 0 & 0 &0\\
      0 & 0 & 0 & N_{\eta}^{(1)}  & .. & 0 & 0 &0 &N_{\eta}^{(n_\eta)} \\
      0 & 0 & N_{\eta}^{(1)}& 0  & .. & 0 & 0 &N_{\eta}^{(n_\eta)}  &0\\
       \\
  \end{bmatrix} \,.
\end{equation}
The operator 
 $\mathcal{G}_\mathrm{sym}$ maps the components of  $\bH_{\eta}$ into a $8\times  4n_\eta$ matrix,  which is defined by
\begin{equation}
\mathcal{G}_\mathrm{sym}=
\begin{bmatrix}
(\bH_\eta)_{11} & 0 & 0 & 0 & \cdots & (\bH_\eta)_{1 n_\eta} & 0 & 0 & 0\\
0 & (\bH_\eta)_{11} & 0 & 0 & \cdots & 0 & (\bH_\eta)_{1 n_\eta} & 0 & 0\\
0 & 0 & \tfrac{1}{2}(\bH_\eta)_{11} & \tfrac{1}{2}(\bH_\eta)_{11} & \cdots & 0 & 0 & \tfrac{1}{2}(\bH_\eta)_{1 n_\eta} & \tfrac{1}{2}(\bH_\eta)_{1 n_\eta}\\
0 & 0 & \tfrac{1}{2}(\bH_\eta)_{11} & \tfrac{1}{2}(\bH_\eta)_{11} & \cdots & 0 & 0 & \tfrac{1}{2}(\bH_\eta)_{1 n_\eta} & \tfrac{1}{2}(\bH_\eta)_{1 n_\eta}\\
(\bH_\eta)_{21} & 0 & 0 & 0 & \cdots & (\bH_\eta)_{2 n_\eta} & 0 & 0 & 0\\
0 & (\bH_\eta)_{21} & 0 & 0 & \cdots & 0 & (\bH_\eta)_{2 n_\eta} & 0 & 0\\
0 & 0 & \tfrac{1}{2}(\bH_\eta)_{21} & \tfrac{1}{2}(\bH_\eta)_{21} & \cdots & 0 & 0 & \tfrac{1}{2}(\bH_\eta)_{2 n_\eta} & \tfrac{1}{2}(\bH_\eta)_{2 n_\eta}\\
0 & 0 & \tfrac{1}{2}(\bH_\eta)_{21} & \tfrac{1}{2}(\bH_\eta)_{21} & \cdots & 0 & 0 & \tfrac{1}{2}(\bH_\eta)_{2 n_\eta} & \tfrac{1}{2}(\bH_\eta)_{2 n_\eta}
\end{bmatrix}
\end{equation}

\bibliographystyle{elsarticle-num} 
\bibliography{References}

\end{document}